\documentclass[12pt]{amsart}
\usepackage{epsfig}
\usepackage{amscd,amsfonts,amsmath,amssymb}

\input xy
\xyoption{all}
\usepackage[all]{xy}

\begin{document}

\title{Vari\'et\'es rationnellement 
connexes sur un corps alg\'ebriquement clos}
\author{Laurent BONAVERO}

\date{}
\maketitle

\def\restriction{\string |}
\def\cad{c'est-\`a-dire}
\def\emptyset{\varnothing}
\def\moins{\mathop{\hbox{\vrule height 3pt depth -2pt
width 5pt}\,}}

\newcommand{\pp}{\rm ppcm}
\newcommand{\pg}{\rm pgcd}
\newcommand{\Ker}{\rm Ker}
\newcommand{\C}{{\mathbb C}}
\newcommand{\Q}{{\mathbb Q}}
\newcommand{\GL}{\rm GL}
\newcommand{\SL}{\rm SL}
\newcommand{\diag}{\rm diag}

\newcommand{\ev}{\operatorname{ev}}
\def\contentsname{Table des mati\`eres}
\def\refname{R\'ef\'erences}
\def\finpreuve
{\hskip 3pt \vrule height6pt width6pt depth 0pt}

\newtheorem{ex}{Exemple}
\newtheorem{theo}{Th\'eor\`eme}
\newtheorem{prop}[theo]{Proposition}
\newtheorem{lemm}[theo]{Lemme}
\newtheorem{lemmf}[theo]{Lemme fondamental}
\newtheorem{defi}[theo]{D\'efinition}
\newtheorem{exo}{Exercice}
\newtheorem{rem}{Remarque}
\newtheorem{cor}[theo]{Corollaire}
\newcommand{\CC}{{\mathbb C}}
\newcommand{\ZZ}{{\mathbb Z}}
\newcommand{\RR}{{\mathbb R}}
\newcommand{\QQ}{{\mathbb Q}}
\newcommand{\FF}{{\mathbb F}}
\newcommand{\PP}{{\mathbb P}}
\newcommand{\codim}{\operatorname{codim}}
\newcommand{\Exc}{\operatorname{Exc}}
\newcommand{\Ho}{\operatorname{Hom}}
\newcommand{\rg}{\operatorname{rg}}
\newcommand{\Pic}{\operatorname{Pic}}
\newcommand{\Aut}{\operatorname{Aut}}
\newcommand{\NE}{\operatorname{NE}}
\newcommand{\Nun}{\operatorname{N}}
\newcommand{\card}{\operatorname{card}}
\newcommand{\Hilb}{\operatorname{Hilb}}
\newcommand{\mult}{\operatorname{mult}}
\newcommand{\vol}{\operatorname{vol}}
\newcommand{\divi}{\operatorname{div}}
\newcommand{\pr}{\operatorname{pr}}
\newcommand{\con}{\operatorname{cont}}
\newcommand{\Spec}{\operatorname{Spec}}
\newcommand{\ima}{\operatorname{Im}}
\newcommand{\Chow}{\operatorname{Chow}}
\newcommand{\lieu}{\operatorname{lieu}}
\newcommand{\ol}[1]{\ensuremath{\mathcal{O}_{#1}}}
\newcommand{\Div}{\ensuremath{\operatorname{Div}}}
\newcommand{\Supp}{\operatorname{Supp}}
\newcommand{\Star}{\operatorname{Star}}
\newcommand{\ph}{\ensuremath{\varphi}}
\newcommand{\lis}[1]{\ensuremath{\{x_1,\dots,x_{#1}}\}}
\newcommand{\som}[1]{\ensuremath{x_1+\dots+x_{#1}}}
\newcommand{\fx}{\ensuremath{\Sigma_X}}
\newcommand{\fy}{\ensuremath{\Sigma_Y}}
\newcommand{\Z}{\ensuremath{\mathbb{Z}}}
\newcommand{\f}{\ensuremath{\Sigma}}
\newcommand{\Int}{\operatorname{Int}}
\newcommand{\N}{\ensuremath{\mathcal{N}_1}}
\newcommand{\A}{\ensuremath{\mathcal{A}_1}}
\newcommand{\PC}{\operatorname{PC}}

\newtheorem{teo}{Th\'eor\`eme}
\newtheorem{lemma}[teo]{Lemme}
\newtheorem{conj}[teo]{Conjecture}

\newcounter{subsub}[subsection]

\def\thesubsub{\thesubsection .\arabic{subsub}}
\def\subsub#1{\addtocounter{subsub}{1}\par\vspace{3mm}
\noindent{\bf \thesubsub ~ #1 }\par\vspace{2mm}}
\def\coker{\mathop{\rm coker}\nolimits}
\def\pr{\mathop{\rm pr}\nolimits}
\def\im{\mathop{\rm Im}\nolimits}
\def\hfl#1#2{\smash{\mathop{\hbox to 12mm{\rightarrowfill}}
\limits^{\scriptstyle#1}_{\scriptstyle#2}}}
\def\vfl#1#2{\llap{$\scriptstyle #1$}\big\downarrow
\big\uparrow
\rlap{$\scriptstyle #2$}}
\def\diagram#1{\def\normalbaselines{\baselineskip=0pt
\lineskip=10pt\lineskiplimit=1pt}   \matrix{#1}}
\def\limind{\mathop{\oalign{lim\cr
\hidewidth$\longrightarrow$\hidewidth\cr}}}

\long\def\InsertFig#1 #2 #3 #4\EndFig{
\hbox{\hskip #1 mm$\vbox to #2 mm{\vfil\includegraphics{#3}}#4$}}
\long\def\LabelTeX#1 #2 #3\ELTX{\rlap{\kern#1mm\raise#2mm\hbox{#3}}}

\newcommand{\num}{\stepcounter{subsection}{\bf\thesubsection.}}
\newcommand{\marque}{\stepcounter{subsection}{\thesubsection. }}
\newcommand{\hs}{\hspace*{5mm}}
\newcommand{\voir}{{\it cf.\ }}

\begin{center}
\begin{minipage}{130mm}
\scriptsize

{\bf R\'esum\'e.} 
Ce sont les notes d'un mini-cours sur les 
vari\'et\'es rationnellement connexes,
\'ecrit
pour les Etats de la Recherche de la Soci\'et\'e
Math\'ematique de France (Strasbourg, 2008)\footnote{
Je remercie tous les participants pour leurs commentaires
qui ont permis d'am\'eliorer grandement ce texte par rapport
\`a la version distribu\'ee le jour de la conf\'erence. St\'ephane Druel 
\`a Grenoble m'a apport\'e une aide inestimable lors de
la pr\'eparation de ce cours.}\footnote{
Institut Fourier, 
UFR de Math\'ematiques,
Universit\'e de Grenoble 1,
UMR 5582,
BP 74,
38402 Saint Martin d'H\`eres,
FRANCE.
e-mail : laurent.bonavero@ujf-grenoble.fr}. On met l'accent
sur les aspects g\'eom\'etriques. Ces notes sont aussi une invitation \`a lire 
le livre d'Olivier Debarre \cite{Deb01}, 
dont une grande partie de ce cours
est extraite. Ces notes doivent surtout permettre au lecteur de comprendre
l'\'enonc\'e suivant et l'un de ses fameux corollaires \cite{GHS03}.
\begin{theo}\label{theoGHS} {\bf (Graber, Harris et Starr)} 
Sur un corps alg\'ebriquement clos,
soient $X$ une vari\'et\'e projective lisse et $\varphi : X \to C$
un morphisme surjectif sur une courbe projective lisse $C$.
Si la fibre g\'en\'erale de $\varphi$ est s\'eparablement rationnellement
connexe, alors $\varphi$ poss\`ede une section.
\end{theo}

\begin{cor} {\bf (Graber, Harris et Starr)} 
Sur un corps alg\'ebriquement clos de caract\'eristique z\'ero,
soit $f~: X \to Y$ un morphisme dominant entre deux vari\'et\'es 
projectives. Si $Y$ et la fibre g\'en\'erale de 
$f$ sont rationnellement connexes,
alors $X$ est rationnellement connexe.
\end{cor}

Les preuves de ces deux r\'esultats sont donn\'ees 
dans le cours de Jason Starr,
le mat\'eriel pr\'eliminaire n\'ecessaire est pr\'esent\'e en d\'etail
dans ces notes.   
Ce cours est r\'edig\'e dans l'espoir de s'adresser \`a un public large,
\`a l'exception peut-\^etre du \S 7, o\`u
nous donnons les d\'etails de la preuve 
de la conjecture de connexit\'e rationnelle de Shokurov
par Hacon et ${\rm M^c}$Kernan,  
plus technique et o\`u les pr\'erequis sont un
peu plus importants.


\end{minipage}
\end{center}

\setcounter{tocdepth}{1}

\tableofcontents

\newpage

\centerline{\bf COURS 1}

\medskip

\begin{center}
\begin{minipage}{130mm}
\scriptsize

Ce premier cours est un survol sur l'importance et la pr\'esence
des courbes
rationnelles en g\'eom\'etrie alg\'ebrique. On traite le cas particulier
des hypersurfaces de l'espace projectif et on discute le lien
entre courbes rationnelles et g\'eom\'etrie birationnelle
classique (\'eclatements, lieux exceptionnel et d'ind\'etermination)
ou moderne (th\'eorie de Mori). On introduit enfin la notion
de connexit\'e rationnelle et des notions qui lui sont reli\'ees.

\end{minipage}
\end{center}

\medskip

\section*{Introduction}

Sauf mention explicite du contraire, toutes les vari\'et\'es alg\'ebriques
et les morphismes consid\'er\'es sont 
d\'efinis {\bf sur un corps $k$ alg\'ebriquement clos
de caract\'eristique arbitraire}, les vari\'et\'es sont irr\'eductibles
et r\'eduites. 
Certains \'enonc\'es ne sont valables 
qu'en caract\'eristique $0$ (ceux dont la preuve n\'ecessite le
th\'eor\`eme de 
lissit\'e g\'en\'erique ou ceux pour lesquels il faut utiliser
une r\'esolution des singularit\'es) ou sur un corps
non d\'enombrable (ceux pour lesquels il est important de
savoir qu'une vari\'et\'e n'est pas r\'eunion d\'enombrable 
de sous-vari\'et\'es propres), on le mentionnera explicitement.

Une courbe est une vari\'et\'e 
projective int\`egre (irr\'eductible et r\'eduite)
de dimension $1$. 

Si $X$ est une vari\'et\'e (quasi-) projective, on dira qu'un point est
en position g\'en\'erale, ou plus simplement
g\'en\'eral, s'il appartient \`a un ouvert non vide 
(non sp\'ecifi\'e) de
$X$, qu'il est en position tr\`es 
g\'en\'erale, ou plus simplement tr\`es g\'en\'eral, 
s'il appartient au compl\'ementaire
d'une r\'eunion d\'enombrable de 
ferm\'es stricts de $X$. Cette derni\`ere notion
n'est que tr\`es peu pertinente sur un 
corps d\'enombrable o\`u le compl\'ementaire
d'une r\'eunion d\'enombrable de ferm\'es 
stricts de $X$ peut \^etre vide.

\section{Quelques exemples de vari\'et\'es poss\'edant 
des courbes rationnelles.}

\subsection{G\'en\'eralit\'es}
Il y a trois grandes classes de courbes projectives lisses~: 
la droite projective $\PP^1$ (dont la topologie complexe
est celle 
d'une sph\`ere), les courbes elliptiques
(dont la topologie complexe est celle d'un tore), 
les courbes de genre $\geq 2$ (dont la topologie complexe 
est celle d'une bou\'ee
multi-places). Il y a \'enorm\'ement 
d'invariants ou outils
permettant de distinguer ces trois classes, celui qui nous sera le plus utile
sera le signe du fibr\'e canonique~: si $X$ est une vari\'et\'e
projective lisse de dimension $n$, 
on note classiquement 
$K_X := -\det T_X$\footnote{
On ne distinguera pas entre 
notation additive et multiplicative pour le groupe de
Picard. Ici, $K_X$ est le dual du fibr\'e $\det T_X$.}. 
C'est le fibr\'e
en droites dont les sections locales sont les $n$-formes r\'eguli\`eres
$f(z_1,\ldots,z_n) dz_1 \wedge \cdots \wedge dz_n$.
Pour une courbe lisse $C$, on a $K_C = T^*_C$ et 
les propri\'et\'es
suivantes sont satisfaites~: $-K_{\PP^1}$ est ample, $K_{E}$
est trivial si $E$ est une courbe elliptique et 
$K_C$ est ample si $C$ est de genre $\geq 2$.

\medskip

Le probl\`eme suivant est un 
probl\`eme classique en g\'eom\'etrie alg\'ebrique~:
si $X \subset \PP^N$ est une vari\'et\'e projective, existe-t-il une courbe
dans $X$ de degr\'e donn\'e et de genre donn\'e~? Si oui, en existe-t-il 
beaucoup et 
que peut-on dire du lieu dans $X$ couvert par ces courbes~? Une situation 
\'el\'ementaire o\`u l'on peut \'enoncer une r\'eponse compl\`ete est le cas 
des courbes planes~: 
si $C \subset \PP^2$ est une courbe lisse plane de genre $g$ et de degr\'e
$d$, alors $g= (d-1)(d-2)/2$. Si $C \subset \PP^2$ est une courbe de
degr\'e $d$ (ceci signifie que $C$ intersecte une droite
g\'en\'erale de $\mathbb P^2$ en $d$ points), 
alors le genre $g$ de sa d\'esingularis\'ee
v\'erifie $g \leq (d-1)(d-2)/2$. En particulier, si $C \subset \PP^2$
est isomorphe \`a $\PP^1$, alors $C$ est une droite ou une conique
dans $\PP^2$ et si $f~: \PP^1 \to \PP^2$ est un morphisme
dont l'image $f(\PP^1)$ est de degr\'e $\geq 3$, 
alors $f(\PP^1)$ est n\'ecessairement 
singuli\`ere. 

\begin{defi} Soit $X$ une vari\'et\'e projective. 
Une courbe rationnelle sur $X$ est un morphisme 
$f : \PP^1 \to X$ non constant\footnote{
Il pourra arriver que par abus de
langage (ou par inattention), on parle aussi de courbes rationnelles
pour d\'esigner l'image d'un morphisme non constant $f : \PP^1 \to X$.}.
Le degr\'e d'une courbe rationnelle $f : \PP^1 \to X$ est 
le degr\'e du morphisme $f$.
\end{defi} 

Le lecteur d\'ebutant doit commencer par se convaincre qu'il
n'y a pas de morphisme $f: \PP^1 \to C$ non constant si
$C$ est une courbe de genre $\geq 1$.\footnote{
Dans sa th\`ese \cite{Laz83}, 
R. Lazarsfeld montre le tr\`es joli r\'esultat
suivant : si $X$ est une vari\'et\'e projective 
complexe lisse de dimension $n$ et 
si $f: \PP^n \to X$ est un morphisme surjectif, alors $X \simeq \PP^n$.
La preuve utilise de fa\c{c}on essentielle 
la g\'eom\'etrie des courbes rationnelles
de $X$. Hwang et Mok ont depuis \'etendu ce r\'esultat \`a de nombreuses 
vari\'et\'es de Fano avec nombre de Picard \'egal \`a $1$.} 
   
\subsection{Courbes rationnelles contenues dans les hypersurfaces
de $\PP^{n}$}
Les courbes ration\-nelles les plus simples dans l'espace projectif
sont les droites. On se demande ici si une hypersurface
(g\'en\'erale) de degr\'e $d$ dans $\PP^{n}$ contient (au moins) 
une droite.

La vari\'et\'e des droites dans $\PP^{n}= \PP(V)$ (o\`u $V$ est un
espace vectoriel de dimension $n+1$\footnote{Dans tout ce cours,
nous suivons la convention na\"{\i}ve~: $\PP(V)$ d\'esigne
la vari\'et\'e des droites vectorielles de $V$.}) est la grassmannienne
$G(2,n+1)$ de dimension $2n-2$. Il y a sur $G(2,n+1)$
un sous-fibr\'e tautologique $E$ de rang $2$ du fibr\'e
trivial $G(2,n+1)\times V$ (la fibre $E_{[l]}$ de
$E$ au dessus de $[l] \in G(2,n+1)$ est le sous-espace vectoriel
de $V$ de dimension $2$ d\'efini par la droite $l \subset \PP(V)$).
Une hypersurface $X_d$ de degr\'e $d$ dans $\PP(V)$ est donn\'ee
par son \'equation, \`a savoir un polyn\^ome homog\`ene de
degr\'e $d$, autrement dit un \'el\'ement $s$ de $S^d(V^*)$. L'hypersurface
$X_d$ contient la droite $l$ si et seulement
si $s_{|E_{[l]}}$ est nul. La sous-vari\'et\'e $F_{X_d}(1,n,d)$ de $G(2,n+1)$
des droites contenues dans $X_d$ est donc le lieu
des z\'eros de $s$ vue comme section du fibr\'e $S^d(E^*)$ sur $G(2,n+1)$.
Il est alors bien connu (en caract\'eristique z\'ero seulement) 
que pour $s$ g\'en\'erale,
$\dim F_{X_d}(1,n,d) = \dim G(2,n+1) - \rg S^d(E^*)$ si cette quantit\'e
est positive ou nulle, et que $F_{X_d}(1,n,d)$ est vide sinon. 
Comme  $\rg S^d(E^*) = d+1$, on en d\'eduit l'\'enonc\'e suivant.

\begin{prop} En caract\'eristique nulle\footnote{Cet \'enonc\'e est en fait
encore vrai en caract\'eristique positive.},
une hypersurface g\'en\'erale de degr\'e $d$ dans $\PP^{n}$
contient une infinit\'e de droites si $d < 2n-3$, un
nombre fini de droites si $d= 2n-3$ et
ne contient pas de droites si $d > 2n-3$.
\end{prop}

Le lecteur int\'eress\'e consultera \cite{DM98} 
pour en savoir beaucoup
plus sur la vari\'et\'e des $r$-plans 
contenus dans une intersection compl\`ete. Il y trouvera
aussi des valeurs num\'eriques du nombre de droites contenues dans
une hypersurface g\'en\'erale de degr\'e $2n-3$ 
dans $\PP^{n}$, dont 
le c\'el\`ebre et classique~: il y a $27$ droites 
sur une cubique lisse de $\PP^{3}$. 

\medskip

Dans le cas o\`u $d < 2n -3$, il est possible de pr\'eciser un peu le lieu
de $X_d$ couvert par les droites contenues dans $X_d$.
Soit en effet $Z \subset X_d \times F_{X_d}(1,n,d)$ 
la vari\'et\'e d'incidence suivante :
$$ Z := \{ (x,[l]) \in X_d \times F_{X_d}(1,n,d) \mid x \in l\}.$$ 
Les fibres de la deuxi\`eme projection $Z \to F_{X_d}$ \'etant de dimension
$1$, on a $\dim Z = 2n-2 -d$. Soit $V= p_1(Z) \subset X_d$ l'image
de $Z$ dans $X_d$ par la premi\`ere projection : c'est le lieu
de $X_d$ couvert par les droites contenues dans $X_d$. 
Si $d\geq n$, 
on d\'eduit de ce qui pr\'ec\`ede que $V$ est une sous-vari\'et\'e stricte
de $X_d$, on montre aussi ais\'ement\footnote{Si $s$ est un polyn\^ome
homog\`ene de degr\'e $d$ 
d\'efinissant $X_d$ et si $x=[1:0:\cdots : 0]\in X_d$,
alors la droite passant par $x$ et un point $[0:x_1:\cdots : x_n]$
est contenue dans $X_d$ si et seulement si $s(t,x_1,\ldots,x_n)=0$
pour tout $t$. Comme $t\mapsto s(t,x_1,\ldots,x_n)$ est un polyn\^ome
en $t$ de degr\'e $\leq d-1$, ce polyn\^ome est nul si et
seulement si ses $d$ coefficients sont nuls, ce qui consiste \`a
r\'esoudre $d$ \'equations homog\`enes de degr\'es respectifs $1,2,\ldots,d$
en les 
variables $[x_1:\cdots : x_n]$.
Il y a au moins une solution si $d \leq n-1$. Dans le cas $d=n-1$, on remarque
qu'il y en a $(n-1)!$, autrement dit par un point g\'en\'eral
de $X_{n-1}\subset \mathbb P^n$ passent $(n-1)!$ droites. Cette 
observation \'el\'ementaire a inspir\'e
un tr\`es joli r\'esultat d\^u \`a J.M. Landsberg \cite{Lan03}.}    
que si $d \leq n-1$, alors $V=X_d$, autrement
dit $X_d$ est couverte par des droites. 
Finalement, dans le cas $d=n$, des arguments analogues
permettent de montrer que $X_n$ est couverte par des coniques\footnote{
{\em Merci \`a Laurent Manivel \`a qui je dois les 
r\'ef\'erences m'ayant permis d'\'ecrire cette ``footnote'', en r\'eponse \`a
une question de Jean-Louis Colliot-Th\'el\`ene.}
Par un point
g\'en\'eral de $X_n \subset \mathbb P ^n$ passent un nombre fini de
coniques. Il n'y a pas, \`a ma connaissance, de m\'ethodes \'el\'ementaires
pour d\'eterminer ce nombre en dimension quelconque. 
La formule g\'en\'erale 
d\'ecoule d'un calcul d'invariants
de Gromov-Witten \`a l'aide de la sym\'etrie 
miroir utilisant une \'equation diff\'erentielle ordinaire introduite
par Givental. Les lignes qui suivent sont issues de la lecture de 
\cite{JNS04} et \cite{Jin05}.

\begin{theo} \label{coniques} 
{\bf (Coates, Givental - Jinzenji, Nakamura, Suzuki)}
Sur le corps des nombres complexes, soit 
$X_n \subset \mathbb P ^n$ une hypersurface g\'en\'erale de degr\'e
$n$ dans $\mathbb P ^n$. Soit $N_n$ le nombre de coniques contenues dans
$X_n$ et passant
par un point g\'en\'eral de $X_n$. 
Alors 
$$ N_n = \frac{(2n)!}{2^{n+1}}-\frac{(n!)^2}{2}.$$ 
\end{theo}

On renvoie \`a \cite{Bea95} et \cite{BH08} pour d'autres r\'esultats sur les 
coniques contenues dans des hypersurfaces.

Expliquons bri\`evement comment on obtient ce r\'esultat. Si $a$, $b$,
$c$ et $d$ sont quatre entiers, notons 
$\langle \mathcal O_a \mathcal O_b \mathcal O_c\rangle_d$
l'invariant de Gromov-Witten comptant 
le nombre (\'eventuellement infini) de courbes 
rationnelles de degr\'e $d$ contenues dans   
$X_n$ et rencontrant $3$ sous-espaces projectifs de $\mathbb P ^n$, 
g\'en\'eraux, 
de codimension respective $a$, $b$ et
$c$. Lorsque $a$, $b$ ou $c$ valent $1$, chaque courbe rationnelle
de degr\'e $d$ contenue dans $X_n$ 
est compt\'ee $d$ fois puisque l'intersection d'une courbe de degr\'e $d$ 
et d'un hyperplan g\'en\'eral est constitu\'ee de $d$ points. 
Comme l'intersection d'une droite g\'en\'erale et de $X_n$ est 
constitu\'ee de $n$ points, le nombre $N_n$ cherch\'e vaut 
donc 
$N_n = \langle \mathcal O _1 \mathcal O _1 \mathcal O _{n-1}\rangle _2 /4n$. 
Dans \cite{Jin05} sont introduites des constantes $\tilde L _m ^{n+1,n,d}$,
dites ``constantes 
de structure de l'anneau de
cohomologie quantique de $X_n$''. Ces constantes satisfont aux formules
r\'ecurrentes suivantes (que l'on explicite uniquement dans les cas $d=1$
et $d=2$)~:
$$ \sum_{m=0}^{n-1}\tilde L _m ^{n+1,n,1}w^m=n\prod_{j=1}^{n-1}(jw+(n-j))$$
et 
$$ \sum_{m=0}^{n-2}\tilde L _m ^{n+1,n,2}w^m=
\sum_{j_2=0}^{n-2}\sum_{j_1=0}^{j_2}\sum_{j_0=0}^{j_1}
\tilde L _{j_1} ^{n+1,n,1}\tilde L _{j_2+1} ^{n+1,n,1}w^{j_1-j_0}\left(
\frac{1+w}{2}\right) ^{j_2-j_1}.$$
 
Il y est aussi montr\'e que pour tout entier $m$, $0 \leq m \leq n-2$,
on a $\tilde L _m ^{n+1,n,2}= 
\langle \mathcal O _1 \mathcal O _{n-1-m} \mathcal O _{m+1}\rangle _2 /n$. 

Le th\'eor\`eme d\'ecoule du calcul du coefficient de $w^{n-2}$ 
dans la deuxi\`eme
formule ci-dessus, de celui de $w^{n-1}$ dans la premi\`ere et
enfin de l'\'evaluation de cette derni\`ere en $w=2$.\finpreuve 

}.

\medskip

Si $X_d$ est une hypersurface lisse de degr\'e $d$ dans $\PP^n$,
son fibr\'e canonique est calcul\'e par la formule d'adjonction
et vaut $K_{X_d} = {\mathcal O}_{\PP^n}(d-n-1)_{|X_d}$. La proposition 
pr\'ec\'edente est une premi\`ere illustration d'un principe
g\'en\'eral : {\em ``plus le fibr\'e canonique $K_X$ est positif, moins il y a
de courbes rationnelles sur $X$, plus ce fibr\'e est n\'egatif, plus il
y a de courbes rationnelles sur $X$''}.
 
\medskip 

De nombreux auteurs (dont Clemens, Voisin et Pacienza)
ont \'etudi\'e l'existence de courbes rationnelles de degr\'e $\geq 2$
dans les hypersurfaces
de $\PP^n$.
Rassemblons leurs r\'esultats (voir \cite{Pac03} et sa bibliographie).

\begin{theo} Supposons le corps de base de caract\'eristique nulle.
Soit $X_d$ une hypersurface tr\`es g\'en\'erale de degr\'e $d$ dans $\PP^n$.
Alors :
\begin{enumerate}
\item {\bf (Clemens)} si $n\geq 3$ et $d \geq 2n-1$, $X_d$ ne contient 
pas de courbes rationnelles,
\item {\bf (Voisin)} si $n\geq 4$ et $d \geq 2n-2$, $X_d$ ne contient 
pas de courbes rationnelles,
\item {\bf (Voisin)} si $n\geq 5$ et si $\delta \geq 1$, $X_{2n-3}$ 
ne contient qu'un nombre fini
de courbes rationnelles de degr\'e $\delta$,
\item {\bf (Pacienza)} si $n\geq 6$ et si $\delta \geq 2$,  
$X_{2n-3}$ ne contient 
pas de courbes rationnelles de degr\'e~$\delta$.
\end{enumerate}
\end{theo}
 
\subsection{Courbes rationnelles provenant de la 
g\'eom\'etrie birationnelle classique}

La g\'eo\-m\'etrie birationnelle consiste \`a classifier les 
vari\'et\'es alg\'ebriques
en identifiant deux vari\'et\'es alg\'e\-briques si elles sont ``isomorphes
sur un ouvert (de Zariski) non vide''.

\begin{defi} Deux vari\'et\'es alg\'ebriques $X$ et $X'$ 
sont birationnellement \'equivalentes s'il existe une 
application rationnelle $\varphi : X \dashrightarrow X'$
et des ouverts non vides $U\subset X$ et $U'\subset X'$
tels que $\varphi_{|U}: U \to U'$ soit un isomorphisme.
Une telle $\varphi$ est une {\em application birationnelle}.
Si $V'$ est une-sous vari\'et\'e de $X'$ non contenue dans
$X'\setminus U'$, 
$V=:\overline{\varphi^{-1}(U'\cap V')}$ est la {\em transform\'ee stricte}
de $V'$ dans $X$.
 
\end{defi}

Un exemple fondamental d'application birationnelle est celui 
des \'eclatements
le long de sous-vari\'et\'es lisses~:
si $Y$ est une sous-vari\'et\'e ferm\'ee
lisse contenue dans le lieu non-singulier
d'une vari\'et\'e alg\'ebrique $X$, il
y a une vari\'et\'e alg\'ebrique $B_Y(X)$ 
et un morphisme birationnel
$\pi : B_Y(X) \to X$ qui se restreint en un 
isomorphisme
$\pi:  B_Y(X) \setminus \pi ^{-1}(Y)
\to X \setminus Y$ et telle que
$\pi ^{-1}(Y) \simeq \PP (N_{Y/X})$ o\`u $N_{Y/X}$
d\'esigne le fibr\'e normal de $Y$ dans $X$.
L'application birationnelle $\pi: B_Y(X) \to X$
s'appelle {\em l'\'eclatement de $X$ le long de $Y$, ou de centre $Y$}
et $E := \pi ^{-1}(Y)$ est le {\em diviseur exceptionnel} de $\pi$. 
Moralement, on remplace chaque point $y$ de $Y$ par
l'espace projectif des directions normales \`a $Y$ dans $X$
passant par $y$. Comme les espaces projectifs contiennent beaucoup
de courbes rationnelles, il y a en particulier des courbes rationnelles
dans $B_Y(X)$. 
Une premi\`ere utilisation des \'eclatements\footnote{
Les \'eclatements jouent un r\^ole
central en g\'eom\'etrie birationnelle comme le montrent
les trois \'enonc\'es difficiles et fondamentaux suivants (voir \cite{Hir64},
\cite{Hir75}, \cite{AKMW02} et le survol \cite{Bon02}).
Ils ne sont connus qu'en caract\'eristique z\'ero, leur preuve d\'epasse
largement le cadre de ce cours.

\begin{theo} {\bf (Hironaka - Th\'eor\`eme de d\'esingularisation)} 
Soit $X \subset \PP^N$ une vari\'et\'e projective.
Alors il existe une suite finie d'\'eclatements 
$Z_p \to Z_{p-1} \to \cdots \to Z_0 = \PP^N$ le long de sous-vari\'et\'es
lisses $Y_i \subset Z_i$ telle que, si $\pi_i$ d\'esigne
la compos\'ee $\pi_i: Z_{i}\to \PP^N$ et $X_i\subset Z_i$ la transform\'ee
stricte de $X = X_0$ sous $\pi_i$, alors 
\begin{enumerate}
\item chaque $Y_i$ est incluse dans le lieu singulier de
$X_i$,
\item la vari\'et\'e $X_p$ est une vari\'et\'e projective lisse.
On dit que $\pi_p : X_p \to X$ est une d\'esingularisation (plong\'ee) de
$X$.
\end{enumerate}
\end{theo} 

\begin{theo}\label{indet} {\bf (Hironaka - Lev\'ee des ind\'eterminations)} 
Soient $X \subset \PP^N$ une vari\'et\'e projective, $X'$ une vari\'et\'e
projective
et $\varphi : X \dashrightarrow X'$ une application rationnelle.
Alors il existe une suite finie d'\'eclatements 
$Z_p \to Z_{p-1} \to \cdots \to Z_0 = \PP^N$ le long de sous-vari\'et\'es
lisses $Y_i \subset Z_i$ telle que, si $\pi_i$ d\'esigne
la compos\'ee $\pi_i: Z_{i}\to \PP^N$ et $X_i\subset Z_i$ la transform\'ee
stricte de $X = X_0$ sous $\pi_i$, alors 
\begin{enumerate}
\item chaque $Y_i$ est incluse dans le lieu d'ind\'etermination de
$\varphi \circ \pi_i$,
\item l'application rationnelle $\varphi \circ \pi_p : X_p 
\dashrightarrow X'$
se prolonge en une application r\'eguli\`ere 
$\varphi \circ \pi_p : X_p \to X'$.
\end{enumerate}
\end{theo}

\begin{theo} {\bf (Abramovich, Karu, Matsuki
et W{\l}odarczyk - Th\'eor\`eme de factorisation)} 
Soit $\varphi:X \dashrightarrow X'$ 
une application birationnelle entre deux vari\'et\'es 
projectives lisses
$X$ et $X'$. 
Alors, $\varphi$ se factorise en 
une suite d'\'eclatements et de contractions de centres lisses.
Autrement dit, il y a une suite 
d'ap\-pli\-ca\-tions birationnelles entre vari\'et\'es projectives lisses
$$ X_1 = V_0 \stackrel{\varphi_0}{\dashrightarrow} V_1
\stackrel{\varphi_1}{\dashrightarrow} \cdots
\stackrel{\varphi_{i-1}}{\dashrightarrow} V_i
\stackrel{\varphi_{i}}{\dashrightarrow} V_{i+1}
\stackrel{\varphi_{i+1}}{\dashrightarrow}
\cdots \stackrel{\varphi_{l-2}}{\dashrightarrow}
V_{l-1} \stackrel{\varphi_{l-1}}{\dashrightarrow} V_l = X_2
$$
de sorte que 
$\varphi = \varphi_{l-1} \circ \varphi_{l-2} \circ \cdots \varphi_1 \circ
\varphi_0$ et pour tout $i$,
$\varphi_i:V_{i} \dashrightarrow V_{i+1}$ ou $\varphi_i^{-1}:V_{i+1}
\dashrightarrow V_{i}$
est un \'eclatement le long d'une sous-vari\'et\'e lisse. 
\end{theo}
}
permet de d\'emontrer les r\'esultats suivants, le premier est connu
sous le nom de lemme d'Abhyankar.

\begin{prop}\label{prol} 
Soient $X$ et $Y$ des vari\'et\'es projectives
et $\pi : X \to Y$ un morphisme birationnel.
On suppose que $Y$ est lisse.
Alors, par un point g\'en\'eral de toute composante irr\'eductible 
de ${\rm Exc}(\pi)$\footnote{${\rm Exc}(\pi)$
d\'esigne le ferm\'e de $X$ form\'e des points au voisinage desquels
$\pi$ n'est pas un isomorphisme.} passe une 
courbe rationnelle contract\'ee par $\pi$. 
\end{prop}

La preuve consiste \`a se ramener au cas o\`u ${\rm Exc}(\pi)$
est irr\'eductible lisse et d'image lisse dans $Y$,
puis \`a faire des \'eclatements successifs de $Y$ le long de 
sous-vari\'et\'es
lisses, le premier \'eclatement se faisant le long de $\pi ({\rm Exc}(\pi))$,
et enfin \`a montrer que ${\rm Exc}(\pi)$ est birationnellement
\'equivalent \`a l'un des diviseurs exceptionnels de cette suite 
d'\'eclatements\footnote{En caract\'eristique z\'ero, on peut, 
en levant les ind\'eterminations de $\pi^{-1}$ par 
une suite d'\'eclatements 
de centres lisses, montrer que toute fibre de $\pi$ est rationnellement
connexe par cha\^{\i}nes. Cette notion sera introduite 
plus loin et ce r\'esultat
sera \'etendu au cas o\`u $X$ est peu singuli\`ere, c'est l'objet de la 
toute derni\`ere partie de ce cours.}.

\begin{cor}\label{prolon} Soient $X$ et $Y$ des vari\'et\'es projectives.
On suppose que $X$ est lisse et que $Y$ ne contient pas de courbe
rationnelle. Alors toute application rationnelle
$\varphi : X \dashrightarrow Y$ se prolonge en une application
r\'eguli\`ere $X \to Y$. 
\end{cor}

{\bf D\'emonstration.} Soit $G_{\varphi} \subset X \times Y$ l'adh\'erence
du graphe de $\varphi$. La premi\`ere projection $p: G_{\varphi} \to X$
est une application birationnelle. Comme $X$ est lisse,
d'apr\`es la proposition pr\'ec\'edente,
toute composante irr\'eductible de ${\rm Exc}(p)$ contient
une courbe rationnelle $C \subset G_{\varphi} \subset
X\times Y$ contract\'ee par $p$. 
Comme $Y$ ne contient pas de courbes
rationnelles, c'est que la deuxi\`eme projection $q : G_{\varphi} \to Y$ 
contracte aussi $C$, ce qui est absurde, une courbe dans $X\times Y$
ne pouvant
\^etre \`a la fois ``horizontale'' et ``verticale''. 
Ainsi ${\rm Exc}(p)= \emptyset$ ce qui implique le r\'esultat.
\finpreuve
 
\medskip

Le r\'esultat pr\'ec\'edent implique de suite. 

\begin{cor} Soit $X$ une vari\'et\'e projective lisse ne contenant
pas de courbes rationnelles.
Alors toute application $\varphi : X \dashrightarrow X$ 
birationnelle se prolonge en un isomorphisme (bi-r\'egulier) de $X$.
\end{cor}

{\bf D\'emonstration.} Par le corollaire pr\'ec\'edent, $\varphi$
et $\varphi^{-1}$ se prolongent en des morphismes r\'eguliers.\finpreuve

\medskip

A l'inverse, les vari\'et\'es $X$ poss\'edant 
beaucoup de courbes rationnelles
ont en g\'en\'eral des applications birationnelles $X \dashrightarrow X$ 
qui ne sont pas des isomorphismes. C'est en particulier le cas 
des espaces projectifs $\PP^n$, $n\geq 2$, dont 
le groupe des transformations birationnelles est beaucoup plus gros
(et compliqu\'e) que son groupe d'isomorphismes r\'eguliers qui n'est autre
que 
le groupe ${\rm PGL}_{n+1}$. Nous verrons plus loin que ce
principe est \`a prendre avec prudence~: on montre en effet parfois 
que certaines 
vari\'et\'es de Fano\footnote{Une vari\'et\'e projective est de Fano
si $-K_X$ est un diviseur de Cartier ample.} 
(poss\'edant beaucoup de courbes rationnelles,
voir \S 6) ne sont pas rationnelles ({\em i.e.} birationnellement
\'equivalentes \`a $\PP^n$) en prouvant que leur groupe de transformations 
birationnelles co\"{\i}ncide avec leur groupe d'isomorphismes r\'eguliers
(c'est le cas de la quartique complexe lisse dans $\PP^4$ par exemple).

\subsection{Courbes rationnelles et th\'eorie de Mori}

Ce paragraphe est une premi\`e\-re introduction
aux liens entre courbes rationnelles et signe du fibr\'e
canonique tels qu'ils apparaissent dans la g\'eom\'etrie birationnelle
``moderne'', \`a savoir la th\'eorie de Mori (ou MMP pour ``Minimal
Model Program'').

Quelques notations : si $f : X \to Y$ est un morphisme entre
deux vari\'et\'es projectives normales, on note 
${\rm N}_1 (X/Y)$ l'espace vectoriel r\'eel engendr\'e par les 
courbes contract\'ees par $f$, modulo l'\'equivalence num\'erique.
Le sous-c\^one convexe ferm\'e de ${\rm N}_1 (X/Y)$ engendr\'e par les 
classes des $1$-cycles effectifs est not\'e $\overline{{\rm NE}}(X/Y)$. 
 
Le r\'esultat \`a la base de la th\'eorie de Mori est le
suivant. Il est le fruit des travaux de nombreux
auteurs dont les principaux sont Benveniste, Kawamata, Koll\'ar, Mori, Reid
et Shokurov. On renvoie au tr\`es r\'ecent et lumineux texte de Druel 
\cite{Dru08} pour une pr\'esentation
des toutes derni\`eres avanc\'ees du MMP par Birkar, Cascini, 
Hacon et ${\rm M^c}$Ker\-nan \cite{BCHM06}.

\begin{theo}\label{cone} {\bf (Th\'eor\`eme du c\^one)}
On suppose que le corps de base est de ca\-ract\'eristi\-que z\'ero.
Soient $X$ une vari\'et\'e projective \`a singularit\'es 
terminales\footnote{
Le lecteur d\'ebutant peut supposer $X$ lisse.}
et $f: X \to Z$ un morphisme sur une vari\'et\'e projective $Z$.
\begin{enumerate}
\item Il existe une famille au plus d\'enombrable $(\Gamma_i)_{i\in I}$
de courbes rationnelles telle que pour tout $i\in I$, 
\begin{enumerate}
\item[(i)] $\dim (f(\Gamma _i))=0$,
\item[(ii)] $0 < -K_X \cdot \Gamma _i \leq 2 \dim (X)$, 
\item[(iii)]  $R_i:= {\mathbb R}^+ [\Gamma_i ]$ est une ar\^ete du c\^one 
$\overline{{\rm NE}}(X/Z)$, 
\item[(iv)] $ \overline{{\rm NE}}(X/Z) =  \overline{{\rm NE}}(X/Z)_{K_X \geq 0}
+ \sum_{i\in I}R_i.$
\end{enumerate} 
\item Soit $i \in I$. Il existe un unique morphisme \`a fibres connexes
$c_i : X/Z \to X_i /Z$ sur une vari\'et\'e projective normale $X_i$ 
tel que, pour toute courbe $C \subset X$, $\dim (c_i(C))=0$
si et seulement si $[C]\in R_i$ ; le morphisme $c_i$ est appel\'e 
la contraction
de $R_i$. 
\end{enumerate}
\end{theo}

Ce th\'eor\`eme dit en particulier que 
si le fibr\'e canonique $K_X$ n'est pas ($f$-) nef\footnote{
On rappelle qu'un fibr\'e en droites $L$ (ou un diviseur de Cartier)
est nef, pour num\'eriquement 
effectif, si
$L\cdot C \geq 0$ pour toute courbe $C$,
qu'il est $f$-nef si $L\cdot C \geq 0$ pour toute courbe
$C$ contenue dans une fibre de $f$.}, alors il existe
une courbe rationnelle $\Gamma$ telle que $K_X \cdot \Gamma <0$.
 
Soit $c_i : X/Z \to X_i /Z$ la contraction
d'une ar\^ete $R_i$ avec $-K_X \cdot R_i >0$.
Deux cas se pr\'esentent~: 
\begin{enumerate}
\item[$\bullet$] si $\dim(X_i) =\dim(X)$,
$c_i$ est une application birationnelle (on distingue alors en g\'en\'eral
deux sous-cas suivant que ${\rm Exc}(c_i)$ est de codimension $1$ - on
dit que $c_i$ est {\em divisorielle} - ou que 
${\rm Exc}(c_i)$ est de codimension $\geq 2$ - on
dit que $c_i$ est {\em petite}), 
\item[$\bullet$] sinon $\dim(X_i) < \dim(X)$, 
$c_i$ est alors une {\em fibration de Mori} dont la fibre
g\'en\'erale est une vari\'et\'e de Fano (en g\'en\'eral singuli\`ere).
\end{enumerate}

Le bilan est donc le suivant sur un corps de caract\'eristique z\'ero~: 
{\em si $X$ est une vari\'et\'e projective 
(disons lisse ou \`a singularit\'es terminales), soit $K_X$ est nef, soit
il y a une courbe rationnelle $\Gamma$ telle que $-K_X\cdot \Gamma >0$
contract\'ee par une application birationnelle ou par une fibration
dont la fibre g\'en\'erale est de Fano}.

\medskip

Comme application du th\'eor\`eme du c\^one,
d\'emontrons la proposition suivante.

\begin{prop}\label{term}
Soient $X$ et $Y$ deux vari\'et\'es projectives \`a singularit\'es
terminales, avec $Y$ $\mathbb Q$-factorielle, 
sur un corps de caract\'eristique z\'ero
et $\pi : X \to Y$ un morphisme birationnel.
Alors, si $\pi$ n'est pas un isomorphisme, 
il existe une courbe rationnelle $C$ contract\'ee par $\pi$ 
telle que $-K_X \cdot C >0$.
\end{prop}

{\bf D\'emonstration\footnote{Voir \cite{Deb01} 
pour une preuve \`a la main quand
$X$ et $Y$ sont lisses.}.}
L'hypoth\`ese que $Y$ est $\mathbb Q$-factorielle assure que
le lieu exceptionnel de $\pi$ est de codimension pure $1$.

Ecrivons 
$$ K_X = \pi ^* K_Y + \sum a_i E_i$$ 
o\`u les $E_i$ sont les composantes irr\'eductibles de 
${\rm Exc}(\pi)$ et les $a_i$ sont tous $>0$ (c'est ici que
l'on utilise l'hypoth\`ese que $Y$ est \`a singularit\'es terminales).

Par l'absurde, si $K_X$ est $\geq 0$ sur toute courbe rationnelle
contract\'ee par $\pi$, c'est que $K_X$ est $\pi$-nef par le th\'eor\`eme
du c\^one (que l'on peut appliquer puisque $X$ est \`a singularit\'es
terminales). 
On en d\'eduit que $\sum a_i E_i = K_X-\pi ^* K_Y $ est aussi
$\pi$-nef. 
Comme $\pi_*(-\sum a_i E_i)=0$, le lemme de n\'egativit\'e\footnote{
Le lemme de n\'egativit\'e est le lemme suivant. 
\begin{lemm}\label{negat}{\bf (Lemme de n\'egativit\'e)}
Soient $\pi: X \to Y$ un morphisme entre deux vari\'et\'es projectives normales
et $B$ un ${\mathbb Q}$-diviseur de Cartier sur $X$. 
On suppose que $-B$ est $\pi$-nef. Alors
$B$ est effectif si et seulement si $\pi_* B$ est effectif.
\end{lemm} 
} implique que les $a_i$ sont tous $\leq 0$ 
ce qui fournit la contradiction.\finpreuve 

\medskip

\section{Cinq notions mesurant la pr\'esence 
de courbes rationnelles.}

On donne ici cinq notions permettant de mesurer la pr\'esence de 
courbes rationnelles. Rappelons \`a nouveau que les vari\'et\'es  
et morphismes consid\'er\'es sont 
d\'efinis sur un corps alg\'ebriquement clos
de caract\'eristique arbitraire.

\begin{defi} Soit $X$ une vari\'et\'e projective de dimension $n$.
On dit que $X$ est 
\begin{enumerate}
\item rationnelle s'il existe une application 
birationnelle $\varphi : \mathbb P ^n \dashrightarrow X$,
\item unirationnelle s'il existe une application rationnelle
dominante $\varphi : \mathbb P ^n \dashrightarrow X$,
\item r\'egl\'ee s'il existe une vari\'et\'e projective $Y$ de dimension
$n-1$ et une application birationnelle  
$\varphi : \mathbb P ^1 \times Y  \dashrightarrow X$,
\item unir\'egl\'ee s'il existe une vari\'et\'e projective $Y$ de dimension
$n-1$ et une application rationnelle dominante 
$\varphi : \mathbb P ^1 \times Y  \dashrightarrow X$,
\item rationnellement connexe s'il existe une vari\'et\'e 
quasi-projective $T$ 
et un morphisme $F : \mathbb P ^1 \times T \to X$
tels 
que le morphisme $\mathbb P ^1 \times \mathbb P ^1 \times T \to X \times X$
qui \`a $(u,u',t) \in \mathbb P ^1 \times \mathbb P ^1 \times T$ associe
$(F(u,t), F(u',t)) \in X \times X$ soit dominant. 

\end{enumerate}
\end{defi} 

Si $X$ v\'erifie l'une des propri\'et\'es ci-dessus, alors $X$
est couverte par des courbes rationnelles au sens o\`u par tout
point de $X$ en position g\'en\'erale passe une courbe rationnelle.
M\^eme si ce cours traite par la suite de la connexit\'e rationnelle,
donnons quelques exemples et les liens entre ces cinq notions.

\medskip

{\bf Remarques.}
\'Evidemment, $X$ rationnelle implique $X$ unirationnelle,
$X$ r\'egl\'ee implique $X$ uni\-r\'egl\'ee. Comme $\mathbb P^{n-1} \times 
\mathbb P ^1$ est une vari\'et\'e rationnelle, $X$ rationnelle implique $X$
r\'egl\'ee et $X$ unirationnelle implique $X$ unir\'egl\'ee (sauf si 
$\dim (X)=0$).
Enfin, $X$ unirationnelle implique $X$ rationnellement connexe
et $X$ rationnellement connexe implique $X$ unir\'egl\'ee (sauf
si $\dim (X)=0$).

\medskip

Les vari\'et\'es de la forme $\mathbb P^1 \times Y$
o\`u $Y$ ne contient pas de courbes rationnelles sont r\'egl\'ees
mais ne sont pas unirationnelles. Plus g\'en\'eralement,  
toutes les implications mentionn\'ees ci-dessus 
sont connues pour \^etre strictes \`a l'exception
de la suivante~:
{\em on ne connait pas de vari\'et\'e rationnellement connexe
non unirationnelle}. Mentionnons que toute cubique complexe lisse 
dans $\mathbb P ^4$ est unirationnelle\footnote{
Soient $X_3 \subset \mathbb P ^4$ une cubique complexe lisse et $L$ une droite
contenue dans $X_3$, il en existe, on a m\^eme d\'ej\`a vu
que les droites contenues dans $X_3$ couvrent $X_3$. 
On note $\Sigma$ l'ensemble
des droites $l$ de $\mathbb P ^4$ telles qu'il existe $p\in L$
avec $p\in l \subset T_p( X_3)$. Si $q \in \mathbb P ^4$ est 
un point g\'en\'eral
fix\'e et si $\Pi \subset \mathbb P ^4$ est un $2$-plan g\'en\'eral fix\'e,
l'application $\varphi : L \times \Pi \dashrightarrow \Sigma$ d\'efinie
par $\varphi (p,r) =  T_p( X_3) \cap \overline{pqr}$ est une application 
birationnelle, donc $\Sigma$ est une 
vari\'et\'e ration\-nelle de dimension $3$.
Par ailleurs, si $l \in \Sigma$ est une droite tangente
\`a $X_3$ passant par un point $p$ de $L$, son intersection avec 
$X_3$ contient
un deuxi\`eme point, permettant de d\'efinir ainsi une 
application rationnelle $\psi : \Sigma \dashrightarrow X_3$. Il est facile
de voir que $\psi$ est dominante, de degr\'e $2$. 
}
et non rationnelle (il s'agit d'un r\'esultat d\^u \`a
Clemens et Griffiths). Mentionnons aussi (voir par exemple
\cite{Mar00} et \cite{Mar06}) qu'il existe des quartiques
lisses unirationnelles 
dans $\mathbb P ^4$ et que 
les
quartiques complexes lisses  
dans $\mathbb P ^4$ ne sont pas  
rationnelles (la non rationalit\'e est due \`a Iskovskikh et Manin). 

Avant de mentionner le joli r\'esultat de Koll\'ar \cite{Kol95}, 
rappelons qu'une 
hypersurface lisse 
dans $\PP^{n}$ de degr\'e $d \leq n$ est de Fano, donc rationnellement
connexe en caract\'eristique z\'ero comme nous le verrons au \S 6.

\begin{theo} {\bf (Koll\'ar)} Sur $\mathbb C$, une hypersurface 
tr\`es g\'en\'erale
dans $\PP^{n}$ de degr\'e $d$ v\'erifiant
$$\frac{2}{3}(n + 2) \leq d \leq n$$ 
n'est pas r\'egl\'ee (et n'est {\em a fortiori} pas rationnelle).
\end{theo}

La preuve de ce r\'esultat se fait en passant en caract\'eristique
positive, Koll\'ar
construit dans le m\^eme article, {\em en caract\'eristique positive}, 
des exemples de vari\'et\'es de Fano non s\'eparablement rationnellement
connexes (voir plus loin pour cette derni\`ere notion).

\medskip

Les cinq notions ci-dessus sont des notions birationnellement invariantes~:
si $X$ est birationnellement \'equivalente \`a $X'$, alors
$X$ est rationnelle (resp. unirationnelle, resp. r\'egl\'ee,
resp. unir\'egl\'ee, resp. rationnellement connexe) si et seulement si $X'$
l'est.

\medskip

D'une certaine fa\c{c}on, les vari\'et\'es qui nous int\'eressent dans ce cours
ont toutes un fibr\'e canonique ``n\'egatif'' en vertu du
r\'esultat profond suivant, d\^u \`a Boucksom, Demailly, P\u aun et Peternell :
une vari\'et\'e est unir\'egl\'ee si et seulement si son fibr\'e canonique
n'est pas ``limite de diviseurs effectifs'' \cite{BDPP04} (voir aussi
le survol \cite{Deb04}).

\begin{theo}{\bf (Boucksom, Demailly, P\u aun et Peternell)}\label{carunireg}
Soit $X$ une vari\'et\'e projective complexe\footnote{Les auteurs
d\'emontrent ce r\'esultat en utilisant des techniques transcendantes. 
Ce th\'eor\`eme est maintenant cons\'equence des travaux de Birkar, Cascini, 
Hacon et ${\rm M^c}$Kernan \cite{BCHM06} et est donc valable sur un corps
alg\'ebriquement clos de caract\'eristique z\'ero.} 
lisse. Alors
$X$ est unir\'egl\'ee si et seulement si $K_X$ n'est pas 
pseudo-effectif\footnote{Un diviseur de Cartier sur $X$
est pseudo-effectif si sa classe
dans ${\rm N}^1(X)$ est dans l'adh\'erence du c\^one engendr\'e par les 
classes de diviseurs effectifs.}.
\end{theo}

Il n'est pas question de d\'emontrer ce r\'esultat ici, mentionnons 
simplement une cons\'equence imm\'ediate, beaucoup plus \'el\'ementaire,
que nous d\'emontrerons plus loin \`a l'aide des courbes rationnelles. 

\begin{prop} En caract\'eristique nulle, les plurigenres d'une vari\'et\'e
projective lisse et uni\-r\'e\-gl\'ee $X$ sont tous nuls :
$$ \forall m >0 \,\,\,\,\,\,\, H^0(X,mK_X) =0.$$
\end{prop}

La r\'eciproque est conjectur\'ee.

\begin{conj} Soit $X$ une vari\'et\'e projective lisse sur un 
corps de caract\'eristique z\'ero.
Alors $X$ est unir\'egl\'ee si et seulement si 
les plurigenres de $X$ sont tous nuls :
$$ \forall m >0 \,\,\,\,\,\,\, H^0(X,mK_X) =0.$$ 
\end{conj}

Cette conjecture est cons\'equence de la ``Conjecture d'abondance'',
qui est, apr\`es les toutes derni\`eres avanc\'ees du programme de
Mori par Birkar, Cascini, 
Hacon et ${\rm M^c}$Kernan, 
la conjecture majeure encore ouverte
dans le programme de Mori.

\newpage

\centerline{\bf COURS 2}

\medskip

\begin{center}
\begin{minipage}{130mm}
\scriptsize

Dans ce deuxi\`eme cours, nous introduisons les notions de courbes rationnelles
libres et tr\`es libres. L'\'etude de leurs d\'eformations permet de
caract\'eriser les vari\'et\'es (s\'eparablement) rationnellement connexes
en terme d'existence de telles courbes.

\end{minipage}
\end{center}

\medskip

\section{Connexit\'e rationnelle. Courbes rationnelles libres et 
tr\`es libres.}

Les r\'esultats de cette section sont dus \`a 
Koll\'ar, Miyaoka et Mori \cite{kmm}. Je me suis beaucoup appuy\'e
sur \cite{Deb01} et \cite{AK03}.

\subsection{Retour sur la d\'efinition de connexit\'e rationnelle}

\begin{defi} Soit $X$ une vari\'et\'e projective de dimension $n$.
On dit que $X$ est 
rationnellement connexe s'il existe une vari\'et\'e quasi-projective
$T$ et un morphisme $F : \mathbb P ^1 \times T \to X$
tels 
que le morphisme $\mathbb P ^1 \times \mathbb P ^1 \times T \to X \times X$
qui \`a $(u,u',t) \in \mathbb P ^1 \times \mathbb P ^1 \times T$ associe
$(F(u,t), F(u',t)) \in X \times X$ soit dominant. 
\end{defi} 

Si $X$ est lisse et le corps de base de caract\'eristique nulle, 
l'application tangente
d'un morphisme dominant est g\'en\'eriquement de rang maximal. 
Ceci n'est plus vrai en caract\'eristique positive\footnote{En 
caract\'eristique positive, un morphisme dominant
peut \^etre de diff\'erentielle 
identiquement nulle : le morphisme
de Frobenius $x\mapsto x^p$ est 
l'exemple le plus c\'el\`ebre.} et il s'av\`ere n\'ecessaire
d'\'etendre un peu la notion de connexit\'e rationnelle.

\begin{defi} Soit $X$ une vari\'et\'e projective de dimension $n$.
On dit que $X$ est s\'epara\-blement
rationnellement connexe s'il existe une vari\'et\'e quasi-projective
$T$ et un morphisme $F : \mathbb P ^1 \times T \to X$
tels 
que le morphisme $\mathbb P ^1 \times \mathbb P ^1 \times T \to X \times X$
qui \`a $(u,u',t) \in \mathbb P ^1 \times \mathbb P ^1 \times T$ associe
$(F(u,t), F(u',t)) \in X \times X$ soit dominant et g\'en\'eriquement lisse. 
\end{defi}

Soit $X$ une vari\'et\'e projective. Pour chaque entier $d$, il
y a d'apr\`es Grothendieck un sch\'ema quasi-projectif 
de type fini, not\'e ${\rm Mor}_d(\PP^1,X)$, param\'etrant
les morphismes de degr\'e $d$ de $\PP^1$
\`a valeurs dans $X$. 
Ce sch\'ema n'est en g\'en\'eral ni r\'eduit, ni irr\'eductible.  
Si $f : \PP^1 \to X$ est un morphisme
de degr\'e $d$, on notera $[f]\in {\rm Mor}_d(\PP^1,X)$ le 
point correspondant. Si $X$ est lisse le long de l'image de $f$,
l'espace tangent de Zariski de ${\rm Mor}_d(\PP^1,X)$ 
au point $[f]$ est isomorphe \`a $H^0(\mathbb P ^1,f^* T_X)$ et 
${\rm Mor}_d(\PP^1,X)$ est lisse  
au point $[f]$ d\`es que $H^1(\mathbb P ^1,f^* T_X)=0$.

\medskip

Comme le degr\'e d'un morphisme est constant en famille,
on en d\'eduit que si $X$ est une vari\'et\'e rationnellement
connexe, il existe un entier $d >0$ tel que le morphisme naturel
$$ \PP^1 \times \PP ^1 \times {\rm Mor}_d(\PP^1,X) \to X \times X$$
soit dominant : autrement dit, 
si $x$ et $y$ sont g\'en\'eraux dans $X$ (au sens o\`u $(x,y)$
est g\'en\'eral dans $X\times X$),
il y a une courbe rationnelle de degr\'e $d$ joignant $x$ \`a $y$. 

En particulier, si $x$ est g\'en\'eral
dans $X$, alors par $y$ g\'en\'eral dans $X$ passe une courbe rationnelle
de degr\'e $d$ issue de $x$. Autrement dit, le morphisme
$\PP^1 \times {\rm Mor}_d(\PP^1,X, 0 \mapsto x) \to X$
qui \`a $(u,[f])$ associe $f(u)$ est dominant\footnote{
${\rm Mor}_d(\PP^1,X, 0 \mapsto x)$ d\'esigne le sous-sch\'ema ferm\'e de 
${\rm Mor}_d(\PP^1,X)$ des morphismes $[f]$ v\'erifiant de plus
$f(0)=x$. Plus g\'en\'eralement,
on notera ${\rm Mor}_d(\PP^1,X,\forall i\, p_i \mapsto x_i)$ 
le sous-sch\'ema ferm\'e de 
${\rm Mor}_d(\PP^1,X)$ des morphismes $[f]$ v\'erifiant de 
plus $f(p_i)=x_i$ pour tout $1\leq i\leq r$. 
Son espace tangent au point $[f]$
est $H^0(\PP^1,f^*T_X\otimes {\mathcal O}_{\mathbb P ^1}(-\sum_{i=1}^r p_i))$,
il est lisse au point $[f]$ si  
$H^1(\PP^1,f^*T_X\otimes {\mathcal O}_{\mathbb P ^1}(-\sum_{i=1}^r p_i))=0$
et sa dimension au point $[f]$ est toujours minor\'ee 
par $-K_X\cdot f_*(\mathbb P ^1) + (1-r)\dim X$.}. 

\medskip

Si le corps de base n'est pas d\'enombrable, on peut montrer la r\'eciproque 
suivante. {\em
Si par deux points g\'en\'eraux d'une vari\'et\'e projective $X$ 
passe une courbe rationnelle,
alors il existe un entier $d$ tel que le morphisme naturel 
$$ \PP^1 \times \PP ^1 \times {\rm Mor}_d(\PP^1,X) \to X \times X$$
est dominant\footnote{Il suffit de consid\'erer le sch\'ema 
localement noeth\'erien de type fini ${\rm Mor}(\PP^1,X)$ 
union {\em d\'enombrable} 
des ${\rm Mor}_d(\PP^1,X)$ pour $d \in {\mathbb N}$ et le 
morphisme associ\'e.}~; $X$ est donc rationnellement connexe.}

\medskip

Il est important de noter que la vari\'et\'e ${\rm Mor}_d(\PP^1,X)$
n'est en g\'en\'eral pas projective mais seulement quasi-projective :
pour $X =\PP^2$ et $0 \neq t \in \mathbb C$, la
famille de morphismes 
$$f_t([u:v])=[t(u^2-v^2):2tuv:u^2+v^2]$$ ne
se compactifie pas dans ${\rm Mor}_2(\PP^1,\PP^2)$
lorsque $t$ tend vers $0$ ou l'infini. Il y a un ph\'enom\`ene
de ``cassage'' qui nous am\`enera \`a consid\'erer 
des ``cha\^{\i}nes de courbes
rationnelles''.

\subsection{Courbes rationnelles libres et tr\`es libres}

Rappelons que le groupe
de Picard de $\mathbb P^1$, et
plus g\'en\'eralement celui de $\mathbb P^m$,
est isomorphe \`a $\mathbb Z$ et que l'on note
$\mathcal O _{\mathbb P^m} (1)$ le g\'en\'erateur
ample de ce groupe. En vertu d'un th\'eor\`eme d\^u 
\`a Grothendieck, tout fibr\'e vectoriel $E$ de rang $r \geq 1$
sur $\mathbb P^1$ s'\'ecrit de fa\c{c}on unique 
$\oplus_{i=1}^r \mathcal O _{\mathbb P^1} (a_i)$ pour des
entiers $a_1\geq \cdots \geq a_r$. Attention, l'\'enonc\'e
correspondant sur $\mathbb P^m$, $m\geq 2$, est faux.

\begin{ex}
Si $m\geq 2$, le fibr\'e tangent $T_{\mathbb P^m}$ n'est 
pas une somme directe de fibr\'es en droites. En revanche,
pour toute droite $l \subset {\mathbb P^m}$,
sa restriction \`a $l$ l'est et la d\'ecomposition ne d\'epend
pas de~$l$~:  
$$  (T_{\mathbb P^m})_{|l} \simeq \mathcal O _{\mathbb P^1} (2) 
\oplus \mathcal O _{\mathbb P^1} (1)^{\oplus (m-1)}.$$
De m\^eme, $T_{\mathbb P^1} \simeq \mathcal O _{\mathbb P^1} (2)$.
\end{ex}

Pour d\'emontrer les r\'esultats principaux concernant les vari\'et\'es
rationnellement connexes, nous allons devoir d\'eformer
(et lisser) les (cha\^{\i}nes de) courbes rationnelles. 
Les courbes rationnelles se d\'eforment
d'autant plus facilement que leur fibr\'e normal a tendance \`a \^etre
positif. Formalisons ceci \`a l'aide d'une d\'efinition maintenant
classique.

\begin{defi} 
Soient $X$ une vari\'et\'e projective de dimension $n$ 
et $f~: \mathbb P ^1 \to X$
une courbe rationnelle. On suppose que $f(\mathbb P^1)$ est
contenu dans le lieu lisse $X_{\rm reg}$ de $X$.
Ecrivons alors
$$ f^*T_X \simeq \oplus_{i=1}^n \mathcal O _{\mathbb P^1} (a_i)
\mbox{ avec } a_1\geq \cdots \geq a_n.$$     
On dit que la courbe rationnelle $f$ est 
\begin{enumerate}
\item libre si $a_n \geq 0$,
\item tr\`es libre si $a_n \geq 1$,
\item $r$-libre (avec $r\geq 0$) si $a_n \geq r$.  
\end{enumerate}
\end{defi}

Un point important \`a retenir est que si $f$ est libre, alors 
le sch\'ema ${\rm Mor}(\PP^1,X)$ est lisse au point $[f]$, de
dimension $h^0(\PP^1,f^*T_X)= \dim H^0(\PP^1,f^*T_X)= n+\sum_{i=1}^n a_i
= -K_X\cdot f_*(\mathbb P^1)+n$.
On parlera alors en particulier de {\em la} 
composante de ${\rm Mor}(\PP^1,X)$
contenant $[f]$.  
De plus, l'ensemble des morphismes $r$-libres 
est un ouvert de ${\rm Mor}(\PP^1,X)$\footnote{Ceci vient
du fait que $f$ est $r$-libre si et seulement si 
$H^1(\PP^1,f^*T_X\otimes  \mathcal O _{\mathbb P^1}(-r-1))=0$
et du th\'eor\`eme de semi-continuit\'e de la cohomologie.}. Nous allons voir
que les courbes libres se d\'eforment 
beaucoup, nous permettant de caract\'eriser
les vari\'et\'es lisses rationnellement connexes \`a l'aide des courbes
tr\`es libres. 

\begin{theo}\label{caracrconn}
Soit $X$ une vari\'et\'e projective lisse de dimension $n$. 
\begin{enumerate}
\item 
Si $X$ contient une courbe rationnelle tr\`es libre $f~: \mathbb P ^1
\to X$, alors
pour tout sous-ensemble fini g\'en\'eral $\{x_1,\ldots,x_m\}$ de $X$,
il existe une courbe rationnelle tr\`es libre sur $X$
passant par tous les $x_i$ dont le degr\'e ne d\'epend
que de celui de $f$ et de $m$.  
En particulier, $X$ est rationnellement
connexe (et m\^eme s\'eparablement rationnellement connexe). 

\item En caract\'eristique z\'ero, si $X$ est rationnellement
connexe, alors, par un point g\'en\'eral de $X$ passe une courbe tr\`es libre.

\item En caract\'eristique quelconque, si $X$ est s\'eparablement 
rationnellement connexe, alors, par un point g\'en\'eral de $X$ 
passe une courbe tr\`es libre.  
\end{enumerate}
\end{theo}

{\bf D\'emonstration.}
Elle se fait en plusieurs \'etapes.

\medskip

{\em Etape 1.} Soient $f~: \mathbb P ^1
\to X$ une courbe $r$-libre ($r\geq 0$), $s$ un entier $\geq 1$
et $${\rm ev}^s : 
(\mathbb P ^1)^s \times {\rm Mor}(\PP^1,X)^{\rm red}
\to X^s $$ le morphisme qui \`a $(u_1,u_2,\ldots,u_s,[g]) 
\in (\mathbb P ^1)^s \times {\rm Mor}(\PP^1,X)$ associe
$(g(u_1),\ldots,g(u_s)) \in X^s$. 
Montrons que si $s \leq r+1$, alors ${\rm ev}^s$
est un morphisme lisse au voisinage du point $(u_1,u_2,\ldots,u_s,[f])$
pour tout $(u_1,u_2,\ldots,u_s)\in (\mathbb P ^1)^s$.

\medskip

En effet, soit $T$ la composante de ${\rm Mor}(\PP^1,X)$
passant par $[f]$. On a vu que $T$ est lisse au voisinage de
$[f]$, de dimension $n+\sum_{i=1}^n a_i \geq n(r+1)$
(avec $f^*T_X \simeq \oplus_{i=1}^n \mathcal O _{\mathbb P^1} (a_i)$).

Il suffit donc de montrer que pour tout 
$(u_1,u_2,\ldots,u_s)\in (\mathbb P ^1)^s$,
la diff\'erentielle de ${\rm ev}^s$ est surjective
au point $(u_1,u_2,\ldots,u_s,[f])$. Cette diff\'erentielle
n'est autre que l'application naturelle
$$ \bigoplus_{i=1}^s T_{u_i} \mathbb P ^1 \oplus H^0 (\mathbb P ^1,f^*T_X)
\to \bigoplus_{i=1}^s T_{f(u_i)}X = \bigoplus_{i=1}^s (f^*T_X)_{f(u_i)}$$
qui \`a $(v_1,\ldots,v_s,\sigma) \in \bigoplus_{i=1}^s T_{u_i} 
\mathbb P ^1 \oplus H^0 (\mathbb P ^1,f^*T_X)$ associe 
$$((Tf)_{u_1}(v_1)+\sigma(u_1),\ldots,(Tf)_{u_s}(v_s)+\sigma(u_s)),$$
elle est donc surjective d\`es que pour tout $i=1,\ldots,n$, 
l'application d'\'evaluation 
$$H^0(\mathbb P ^1,\mathcal O _{\mathbb P^1} (a_i)) \to \bigoplus_{j=1}^s
(\mathcal O _{\mathbb P^1} (a_i))_{u_j}$$ l'est. C'est
le cas si (et seulement si lorsque les $u_j$ sont deux \`a deux distincts) 
$a_i\geq s-1$.\finpreuve 

\medskip

{\em Etape 2.} Sous les hypoth\`eses de l'\'etape 1, l'application 
${\rm ev}^s : (\mathbb P ^1)^s \times T
\to X^s$ est donc dominante, 
autrement dit, par $s$ points g\'en\'eraux de $X$ passe 
une courbe $r$-libre d\'eformation de $f$.

\medskip

{\em Etape 3.} Montrons maintenant le point (1) du th\'eor\`eme.
Soit  $f~: \mathbb P ^1
\to X$ une courbe tr\`es libre, on \'ecrit 
$f^*T_X \simeq \oplus_{i=1}^n \mathcal O _{\mathbb P^1} (a_i)$
et soit $r:=\min(a_i) \geq 1$.
En composant $f$ \`a la source par un morphisme 
$h_{\delta}~: \mathbb P^1 \to \mathbb P^1$ de degr\'e $\delta$, on obtient 
une courbe $r \delta$-libre, \`a savoir $f\circ h_{\delta}$. 
Des \'etapes pr\'ec\'edentes, on d\'eduit
que si $r \delta + 1 \geq m$, alors par $m$ points g\'en\'eraux    
de $X$ passe une d\'eformation de $f\circ h_{\delta}$.\footnote{
On 
ne r\'esiste pas ici au commentaire suivant : la preuve qui pr\'ec\`ede montre
que si $C$ est une courbe rationnelle tr\`es 
libre, un multiple suffisamment
grand de $C$ se d\'eforme suffisamment pour passer par $m$ points
g\'en\'eraux de $X$. Ceci fonctionne bien 
car $\mathbb P ^1$ poss\`ede des 
endomorphismes de degr\'e arbitrairement grand, ceci fonctionnerait encore
pour des courbes elliptiques. Pour des courbes de genre $\geq 2$, 
ceci fonctionne encore en caract\'eristique positive seulement car on dispose
du morphisme de Frobenius. Cette remarque g\'eniale est due \`a Mori
et est le point de d\'epart de la th\'eorie de Mori.}

\medskip

{\em Etape 4.} Supposons que $X$ est rationnellement connexe.
On a vu qu'il existe alors un entier $d$ tel que 
le morphisme naturel
$$ {\rm ev}^2 : 
\PP^1 \times \PP ^1 \times {\rm Mor}_d(\PP^1,X)^{\rm red} \to X \times X$$
soit dominant.  
Quand la caract\'eristique du corps de base est suppos\'ee nulle, 
il existe $(u_1,u_2,[f]) \in \PP^1 \times \PP ^1 \times 
{\rm Mor}_d(\PP^1,X)^{\rm red}$
tel que la diff\'erentielle de ${\rm ev}^2$ en $(u_1,u_2,[f])$
soit surjective, et 
ceci \'etant une condition ouverte, on peut supposer $u_1\neq u_2$.
Un tel $(u_1,u_2,[f])$ existe aussi en caract\'eristique positive
quand $X$ est suppos\'ee 
s\'eparablement rationnellement connexe (par d\'efinition !). 
Nous allons montrer que $f$ est tr\`es libre.

A nouveau, 
la  diff\'erentielle de ${\rm ev}^2$ en $(u_1,u_2,[f])$ 
est l'application naturelle 
$$ T_{u_1}\mathbb P ^1 \oplus T_{u_2 }\mathbb P ^1 
\oplus H^0 (\mathbb P ^1,f^*T_X)
\to T_{f(u_1)}X\oplus T_{f(u_2)}X$$
qui \`a $(v_1,v_2,\sigma) \in T_{u_1}\mathbb P ^1 \oplus T_{u_2 }\mathbb P ^1 
\oplus H^0 (\mathbb P ^1,f^*T_X)$ associe 
$$((Tf)_{u_1}(v_1)+\sigma(u_1),(Tf)_{u_2}(v_2)+\sigma(u_2)).$$  
Le point cl\'e est le suivant~: 
comme $T_{\mathbb P ^1}\simeq \mathcal O _{\mathbb P^1} (2)$,
l'application 
$H^0(\mathbb P ^1,T_{\mathbb P ^1})\to
T_{u_1}\mathbb P ^1 \oplus T_{u_2 }\mathbb P ^1$ est surjective, 
donc
l'image de l'application
$H^0(\mathbb P ^1,f^*T_X) \to T_{f(u_1)}X\oplus T_{f(u_2)}X$ 
contient celle
de $(Tf)_{u_1}\oplus (Tf)_{u_2} : 
T_{u_1}\mathbb P ^1 \oplus T_{u_2 }\mathbb P ^1 \to 
T_{f(u_1)}X\oplus T_{f(u_2)}X$, d'o\`u l'on d\'eduit \'evidemment\footnote{Si
$f_1 : E_1 \to F$ et $f_2 : E_2 \to F$
sont deux applications lin\'eaires telles que
d'une part l'image de $f_1$ est contenue dans celle de $f_2$
et d'autre part $f_1 \oplus f_2 : E_1\oplus E_2 \to F$ est
surjectif, alors $f_2$ est surjectif.} 
que  
$H^0 (\mathbb P ^1,f^*T_X)
\to T_{f(u_1)}X\oplus T_{f(u_2)}X$ est surjective. 
De l\`a, si $f^*T_X \simeq \oplus_{i=1}^n \mathcal O _{\mathbb P^1} (a_i)$,
chaque 
$$H^0(\mathbb P ^1,\mathcal O _{\mathbb P^1} (a_i)) \to 
(\mathcal O _{\mathbb P^1} (a_i))_{u_1}\oplus 
(\mathcal O _{\mathbb P^1} (a_i))_{u_2}$$ est surjective, donc $a_i \geq 1$
pour tout $i$ : $f$ est tr\`es libre.\finpreuve 

\medskip 

Le lecteur d\'ebutant pourra d\'emontrer le strict analogue
du th\'eor\`eme pr\'ec\'edent pour les vari\'et\'es unir\'egl\'ees.
Le r\'esultat de l'\'etape 4 ci-dessus 
doit \^etre adapt\'e de la fa\c{c}on suivante~:
{\em si la diff\'erentielle de $ {\rm ev}^1 : 
\PP^1 \times {\rm Mor}_d(\PP^1,X)^{\rm red} \to X$
est surjective au point $(u,[f])$, alors $f$ est libre.}

\begin{theo}\label{caracunireg}
Soit $X$ une vari\'et\'e projective lisse de dimension $n$. 
\begin{enumerate}
\item 
Si $X$ contient une courbe rationnelle libre, alors
$X$ est unir\'egl\'ee.

\item En caract\'eristique z\'ero, si $X$ est unir\'egl\'ee, 
alors par un point g\'en\'eral de $X$ passe une courbe libre.
En particulier, si $X$ est unir\'egl\'ee, alors $K_X$ n'est
pas num\'eriquement effectif.
\end{enumerate}
\end{theo}

En caract\'eristique positive, les choses sont
sont bien diff\'erentes comme le montre l'exemple suivant.

\begin{ex} Soit 
$$X= \{[x_0:\cdots:x_n]\in \mathbb P^n \mid
x_0^d+\cdots+x_n^d=0\}$$ l'hypersurface de Fermat de degr\'e $d=p^r+1$
dans $\mathbb P ^n$ sur un corps alg\'ebriquement clos
de caract\'eristique positive $p$.
Alors $X$ est unirationnelle si $n\geq 3$ (donc rationnellement
connexe). Par ailleurs, si $d \geq n+1$, le fibr\'e canonique
de $X$ est num\'eriquement effectif (son intersection avec toute
courbe est $\geq 0$), en particulier, $X$ ne contient pas de 
courbes libres : la diff\'erentielle de $ {\rm ev}^1 : 
\PP^1 \times {\rm Mor}_d(\PP^1,X) \to X$ n'est donc jamais surjective.
{\em A fortiori}, $X$ n'est pas s\'eparablement rationnellement connexe.
\end{ex}
 
Pour le confort du lecteur, \'enon\c{c}ons deux corollaires,
extraits de la preuve du th\'eor\`eme \ref{caracrconn}.

\begin{cor} Soit $X$ une vari\'et\'e projective lisse
sur un corps de caract\'eristique quelconque,
\begin{enumerate}
\item $X$ est s\'eparablement rationnellement connexe si et seulement si par
tout point g\'en\'eral de $X$ passe une courbe tr\`es libre,
\item $X$ est s\'eparablement unir\'egl\'ee si et seulement si par
tout point g\'en\'eral de $X$ passe une courbe libre.
\end{enumerate}
\end{cor}

\begin{cor}
Soient $X$ une vari\'et\'e projective
et $f~: \mathbb P ^1 \to X$ une courbe rationnelle dont l'image
est contenue dans le lieu lisse de $X$.
\begin{enumerate}
\item En caract\'eristique quelconque,
\begin{enumerate}
\item $f$ est tr\`es libre si et seulement s'il 
y a $(u_1,u_2,[f]) \in \PP^1 \times \PP ^1 \times {\rm Mor}(\PP^1,X)$
tel que la diff\'erentielle de ${\rm ev}^2$ en $(u_1,u_2,[f])$
est surjective,
\item  $f$ est libre si et seulement s'il 
y a $(u,[f]) \in \PP ^1 \times {\rm Mor} (\PP^1,X)$
tel que la diff\'erentielle de ${\rm ev}^1$ en $(u,[f])$
est surjective.
\end{enumerate}
\item En caract\'eristique z\'ero,
\begin{enumerate}
\item $f$ est tr\`es libre si et seulement si les d\'eformations de
$f$ passent par deux points g\'en\'eraux de $X$.
\item $f$ est libre si et seulement si les d\'eformations de
$f$ dominent $X$.
\end{enumerate}
\end{enumerate}

\end{cor}

\subsection{Lieu des courbes libres et courbes libres
minimales}

Avec les techniques pr\'ec\'e\-dentes, on montre la proposition suivante.

\begin{prop} Soit $X$ une vari\'et\'e projective lisse
sur un corps de caract\'eristique z\'ero. Alors, il existe
une intersection d\'enombrable
d'ouverts non vides de $X$, not\'ee $X^{\rm libre}$, 
telle que toute courbe rationnelle
dont l'image rencontre $X^{\rm libre}$ est libre.
De plus, $X^{\rm libre}\neq \emptyset$ si et seulement si $X$ est 
unir\'egl\'ee. 

\end{prop}

{\bf D\'emonstration.}
Consid\'erons \`a nouveau le morphisme d'\'evaluation 
$$ {\rm ev}^1 : \PP^1 \times {\rm Mor}(\PP^1,X) \to X.$$
Pour chaque composante irr\'eductible $M_i$ de ${\rm Mor}(\PP^1,X)$, 
deux situations sont possibles :
\begin{enumerate}
\item[$\bullet$] soit $\overline{{\rm ev}^1 (\PP^1 \times M_i)}\neq X$,
on pose $U_i = X\setminus 
\overline{{\rm ev}^1 (\PP^1 \times M_i)}$,

\item[$\bullet$] soit $\overline{{\rm ev}^1 (\PP^1 \times M_i)}= X$,
comme on a suppos\'e que la caract\'eristique du corps
de base est z\'ero, il y a un ouvert non vide $U_i$ de $X$
tel que pour tout $(u,[f])\in \PP^1 \times M_i ^{\rm red}$, si $f(u)\in U_i$
alors la diff\'erentielle de ${\rm ev}^1$ en $(u,[f])$
est surjective (en particulier, si $f(u)\in U_i$, alors $f$
est libre).
\end{enumerate}
Il est alors clair que $X^{\rm libre}:= \cap_i U_i$ convient.\finpreuve

\medskip

Le lecteur familier avec le lemme de cassage\footnote{Il s'\'enonce
de la fa\c{c}on suivante.
\begin{theo} {\bf (Lemme de cassage)} Soient $X$ une vari\'et\'e projective
et $f : \mathbb P^1 \to X$ une courbe rationnelle
telle que $\dim_{[f]}({\rm Mor} (\mathbb P ^1,X,0\mapsto f(0),
\infty \mapsto f(\infty))\geq 2$. Alors le $1$-cycle
$f_*(\mathbb P ^1)$ est num\'eriquement \'equivalent \`a un $1$-cycle
connexe passant par $f(0)$ et $f(\infty)$, effectif
et non int\`egre de courbes rationnelles.
\end{theo}

Remarquons que l'on a toujours $\dim_{[f]}
{\rm Mor} (\mathbb P ^1,X,0\mapsto f(0),
\infty \mapsto f(\infty))\geq 1$ car on peut composer $f$ \`a la source
par les automorphismes de $\mathbb P^1$ fixant $0$ et $\infty$.
Ce th\'eor\`eme affirme donc que si une courbe rationnelle se 
d\'eforme ``vraiment''
en fixant deux points, alors elle d\'eg\'en\`ere en un cycle non
irr\'eductible et/ou non r\'eduit de courbes
rationnelles. Le lecteur calculera les
d\'eg\'en\'erescences de la famille des coniques planes
$$f_t([u:v])=[t(u^2-v^2):2tuv:u^2+v^2]$$ passant par les deux points
$[1:\pm i:0]$.}
l'utilisera associ\'e
aux  
techniques pr\'ec\'edentes pour d\'emontrer les deux premiers points
du th\'eor\`eme suivant, le troisi\`eme est beaucoup plus
difficile, il est d\^u \`a Cho, Miyaoka et Shepherd-Barron (voir
\cite{CMS00} ou \cite{Keb01}).

\begin{theo} Soit $X$ une vari\'et\'e projective lisse
unir\'egl\'ee, de dimension $n$, sur un corps de caract\'eristique z\'ero. 
Soit $H$ un diviseur ample sur $X$ et $d= \min H\cdot h_*(\mathbb P ^1)$
o\`u le minimum est pris sur l'ensemble des 
courbes {\bf libres}
$h~:\mathbb P ^1 \to X$.
\begin{enumerate}
\item[(1)] Soit $f~:\mathbb P ^1 \to X$ une courbe libre
telle que $d= H\cdot f_*(\mathbb P ^1)$. Alors, il existe
$s$ tel que
$$ f^*T_X \simeq  \mathcal O _{\mathbb P^1} (2)\oplus
\mathcal O _{\mathbb P^1} (1)^{\oplus s}\oplus
\mathcal O _{\mathbb P^1}^{\oplus n-1-s}.$$
\item[(2)] Soient $x\in X^{\rm libre}$, 
$f~:\mathbb P ^1 \to X$ une courbe libre
telle que $d= H\cdot f_*(\mathbb P ^1)$ et $f(0)=x$.
Soit $M_x$ la composante de ${\rm Mor}(\mathbb P ^1,X,0\mapsto x)$
passant par $[f]$. Alors $M_x$ est propre.
\item[(3)] {\bf (Cho, Miyaoka et Shepherd-Barron)}
S'il existe $f~:\mathbb P ^1 \to X$ une courbe libre
telle que $d= \min H\cdot f_*(\mathbb P ^1)$ et
$$ f^*T_X \simeq  \mathcal O _{\mathbb P^1} (2)\oplus
\mathcal O _{\mathbb P^1} (1)^{\oplus n-1},$$
alors $X\simeq \mathbb P ^n$.
\end{enumerate}
\end{theo}

Suivant les auteurs, les courbes libres $f~: \mathbb P ^1 \to X$ 
de degr\'e minimal pour une polarisation donn\'ee,
ou plus g\'en\'eralement telles que la composante de 
${\rm Mor}(\mathbb P ^1,X,0\mapsto x)$
passant par $[f]$ est propre pour $x$ (tr\`es) g\'en\'eral dans $X$
sont dites {\em minimales}.
L'\'etude des courbes libres minimales sur les vari\'et\'es
de Fano dont le rang du groupe de Picard vaut $1$ 
est aussi l'objet de toute une s\'erie de travaux dus \`a Hwang et Mok,
ainsi qu'\`a Kebekus.

\newpage

\centerline{\bf COURS 3}

\medskip

\begin{center}
\begin{minipage}{130mm}
\scriptsize

Dans ce troisi\`eme cours, 
nous introduisons la notion de connexit\'e rationnelle
par cha\^{\i}nes. Cette notion est beaucoup plus souple que
la connexit\'e rationnelle. Nous \'etudions aussi les techniques
de lissage de cha\^{\i}nes de courbes rationnelles. Ces techniques
seront utilis\'ees au cours suivant pour montrer que 
les notions de connexit\'e rationnelle par cha\^{\i}nes
et de connexit\'e rationnelle co\"{\i}ncident 
dans la cat\'egorie des vari\'et\'es projectives lisses
sur un corps de caract\'eristique z\'ero.  

\end{minipage}
\end{center}

\section{Connexit\'e rationnelle par cha\^{\i}nes. Techniques de lissage.}
Les r\'esultats de cette section sont dus \`a 
Koll\'ar, Miyaoka et Mori \cite{kmm}. 
Je me suis \`a nouveau beaucoup appuy\'e
sur \cite{Deb01} et \cite{AK03}.

\subsection{Connexit\'e rationnelle par cha\^{\i}nes}

On a vu pr\'ec\'edemment que ${\rm Mor}_d(\PP^1,X)$
n'est pas projectif en g\'en\'eral, mais seulement quasi-projectif.
Le ph\'enom\`ene de d\'eg\'en\'erescences 
qui intervient fait appara\^{\i}tre des cha\^{\i}nes
de courbes rationnelles. On introduit ici la classe des vari\'et\'es
telles que par deux points g\'en\'eraux passe une cha\^{\i}ne de courbes
rationnelles\footnote{Une cha\^{\i}ne de courbes rationnelles
est un sch\'ema propre de dimension
$1$, connexe, dont toutes les composantes irr\'eductibles sont des courbes 
rationnelles, courbe rationnelle signifiant {\bf exceptionnellement}
ici courbe dont la normalis\'ee est $\mathbb P^1$.}.

\begin{defi} Soit $X$ une vari\'et\'e (quasi-) projective ou un sch\'ema
connexe.
On dit que $X$ est rationnellement connexe par cha\^{\i}nes
s'il existe une vari\'et\'e quasi-projective
$T$ et un sous-sch\'ema $\mathcal C$ de
$T\times X$ tels que 
\begin{enumerate}
\item les fibres de la projection $p~: \mathcal C \to T$
sont propres, connexes, de dimension $1$ et toutes leurs composantes
irr\'eductibles
sont rationnelles\footnote{Courbe rationnelle
signifie {\bf (\`a nouveau exceptionnellement !) }
ici courbe dont la normalis\'ee est $\mathbb P^1$.},
\item la projection $e~: \mathcal C \times _T \mathcal C \to X \times X$
est dominante. 
\end{enumerate}
\end{defi}


Si $X$ est rationnellement connexe par cha\^{\i}nes, alors
par deux points g\'en\'eraux 
$x$ et $x'$ de $X$, il existe $t\in T$ tel que
$x$ et $x'$ appartiennent \`a $q(p^{-1}(t))$ o\`u $q$
est la projection $q~:\mathcal C \to X$, autrement dit, 
la cha\^{\i}ne de courbes rationnelles $q(p^{-1}(t))$ passe par $x$
et par $x'$. 

\medskip

{\bf Mise en garde.} Contrairement 
aux notions pr\'ec\'edemment rencontr\'ees,
cette notion n'est pas invariante par transformation birationnelle~:
si $Y$ est un c\^one sur une vari\'et\'e projective $Y_0$, $Y$
est rationnellement connexe par cha\^{\i}nes 
(passer par le sommet du c\^one). En revanche,
le ``cylindre'' $X$ obtenu en \'eclatant le sommet du c\^one
n'est pas rationnellement connexe par cha\^{\i}nes si $Y_0$
ne l'est pas.  

\begin{ex}
Sur un corps de caract\'eristique z\'ero, consid\'erons
$$ \pi : \mathcal X = \{(w,x,y,z,t)\in \mathbb 
P_{w,x,y,z} ^3\times \mathbb A_t ^1\mid
x^3+y^3+z^3=tw^3 \} \to \mathbb A_t ^1.$$
La vari\'et\'e $\mathcal X$ est lisse, \'etudions les fibres
de $\pi$. 
Si $t\neq 0$, la fibre $\pi^{-1}(t)$ est une surface cubique
lisse de $\mathbb P^3$~; elle est rationnellement connexe
par cha\^{\i}nes (deux points g\'en\'eraux peuvent \^etre reli\'es
par une cha\^{\i}ne form\'ee d'une droite et de deux coniques), 
donc rationnellement
connexe comme nous le verrons plus loin (th\'eor\`eme \ref{versus}). 
Si $t=0$, la fibre $\pi^{-1}(0)$ est un c\^one sur une courbe
elliptique~; elle est rationnellement connexe par cha\^{\i}nes
mais n'est pas rationnellement connexe. 
\end{ex}

Le fait que des cha\^{\i}nes de courbes rationnelles
ne peuvent d\'eg\'en\'erer qu'en des cha\^{\i}nes de courbes rationnelles
implique que si $X$ est rationnellement connexe par cha\^{\i}nes, alors
par deux points {\em quelconques} de $X$ passe une 
cha\^{\i}ne de courbes rationnelles.

\medskip

Evidemment, les vari\'et\'es rationnellement connexes 
sont rationnellement connexes par cha\^{\i}\-nes, le
but des lignes qui suivent est de d\'emontrer la r\'eciproque
pour les vari\'et\'es lisses en caract\'eristique z\'ero.

\subsection{Lissage de cha\^{\i}nes de courbes rationnelles}

Cette section est plus technique, le lecteur
d\'ebutant pourra prendre comme une bo\^{\i}te noire les
deux th\'eor\`emes de lissage.

\subsubsection{Lissage des arbres de courbes rationnelles}

\begin{defi} 
Un arbre de courbes rationnelles est 
un sch\'ema $C$ de dimension un, r\'eduit, connexe 
dont les composantes irr\'eductibles 
\begin{enumerate}
\item sont des courbes rationnelles lisses,
\item peuvent \^etre num\'erot\'ees en choisissant arbitrairement 
l'une d'elles comme \'etant $C_1$ et de sorte
que pour tout $i \geq 2$, $C_i$ rencontre
$C_1 \cup \cdots \cup C_{i-1}$ transversalement en un 
unique point lisse de $C_1 \cup \cdots \cup C_{i-1}$.
\end{enumerate}
Les seules singularit\'es d'un arbre de courbes rationnelles
sont donc des n\oe{}uds ordinaires, correspondant aux points
d'intersection de deux composantes. Le graphe non orient\'e 
ayant pour sommets les composantes de $C$, avec une ar\^ete entre deux 
sommets si les deux composantes correspondantes s'intersectent,
est un arbre.
\end{defi}

Soient $X$ une vari\'et\'e projective lisse, $C$ un arbre de courbes
rationnelles, $f: C \to X$
un morphisme et $p_1,\ldots,p_r$ des points lisses deux \`a deux distincts
de $C$.

\begin{defi}
On dit que $f$ est lissable en fixant $f(p_1),\ldots,f(p_r)$
s'il existe 
\begin{enumerate}
\item une courbe quasi-projective lisse $T$ avec un point distingu\'e $o\in T$,
une surface $\mathcal C$ 
et un morphisme plat $\pi~: \mathcal C \to T$
tels que $\pi^{-1}(o) = C$ et $\pi^{-1}(t)$ est une courbe rationnelle
lisse pour tout $t\in T\setminus \{o\}$,
\item $r$ sections $\sigma_i : T \to  \mathcal C$ de $\pi$
telles que $\sigma_i(o)=p_i$ pour tout $i$,
\item un morphisme $F~:\mathcal C \to X$ tel que 
$F_{|\pi ^{-1}(o)} = f$ et $F(\sigma_i (T))=f(p_i)$ pour tout $i$.
\end{enumerate}
\end{defi}

Cette d\'efinition 
appelle quelques commentaires. 

\begin{enumerate}
\item[(i)] Une surface $\mathcal C$ v\'erifiant les seuls points
(1) et (2) de la d\'efinition est facile 
\`a construire en \'eclatant $\mathbb P ^1 \times \mathbb A^1$
convenablement.
\item[(ii)] Si $f$ est lissable en fixant $f(p_1),\ldots,f(p_r)$,
pour $t\neq o$, $F_t~:\mathbb P^1 \to X$ est une courbe
rationnelle (dont l'image peut \^etre singuli\`ere) passant par  
$f(p_1),\ldots,f(p_r)$.
\item[(iii)] En g\'en\'eral, $f$ n'est pas lissable~: 
supposons que $S$ est une surface projective lisse 
et que $C\subset S$ est un arbre \`a deux composantes 
avec $C_1^2=C_2^2=-2$. Alors l'inclusion $C\subset S$
n'est pas lissable. En effet, il existe une surface singuli\`ere
$S_0$, un morphisme birationnel $h~: S \to S_0$
tels que $h(C)=x_0 \in S_0$ et $h~: S \setminus C \simeq S_0\setminus
\{x_0\}$. Si $F~:\mathcal C \to S$ 
est un lissage de $f$, notons 
$F_t$ est la restriction de $F$ \`a $\pi^{-1}(t)$.
Par le lemme de rigidit\'e, $h \circ F_t$ est constant pour
tout $t\in T$ (on peut supposer $T$ irr\'eductible affine), donc
$F_t$ est constant pour tout $t\neq o$, ce qui est absurde.
\item[(iv)] Les g\'eom\`etres italiens, mais aussi Noether ou Halphen,
se sont int\'eress\'es au lissage de courbes gauches. Le lissage 
des arbres de
courbes de $\mathbb P^3$ est \'etudi\'e
et utilis\'e par Hartshorne et Hirschowitz (voir \cite{HH83}
et le th\'eor\`eme \ref{HH} plus loin)~; Koll\'ar a depuis
donn\'e le formalisme
g\'en\'eral dans son livre \cite{Kol96}.

\end{enumerate}

Dans l'exemple (iii) ci-dessus, $C_1$ et $C_2$ ne sont pas des courbes 
libres~; bien au contraire, leur fibr\'e normal est ``tr\`es'' n\'egatif,
ce qui assure que l'on peut les contracter sur un point. Le th\'eor\`eme qui 
suit montre que la situation est bien meilleure
pour les courbes libres.

\begin{theo}
Soient $X$ une vari\'et\'e projective lisse, 
$C$ un arbre de courbes rationnelles,
$f~: C \to X$
un morphisme, $p_1, \ldots, p_r$ des points lisses distincts
de $C$ dont $r_i$ exactement sont situ\'es sur la composante $C_i$
(chaque $r_i$ est $\geq 0$~; en particulier, ils peuvent \^etre tous nuls).
Soient $f_i$, $1\leq i \leq r$ les restrictions de 
$f$ aux composantes $C_i$.
Si $f_1$ est $(r_1-1)$-libre\footnote{Dans le cas o\`u $r_1 =0$, ceci
signifie que tous les $a_i$ sont $\geq -1$ dans la d\'ecomposition
$f_1^*T_X \simeq \oplus_{i=1}^n \mathcal O _{\mathbb P^1} (a_i)$.}
et $f_i$ est $r_i$-libre pour tout $i\geq 2$, 
alors $f$ est lissable
en fixant tous les $f(p_i)$ et on peut supposer que $F_t$ est 
$(r-1)$-libre pour tout
$o \neq t \in T$. 
\end{theo}

{\bf D\'emonstration.}

{\em Etape 1.}
Soient une courbe lisse $T$ avec un point distingu\'e $o\in T$,
une surface lisse $\mathcal C$ et un morphisme $\pi~: \mathcal C \to T$
tels que $\pi^{-1}(o) = C$ et $\pi^{-1}(t)$ est une courbe rationnelle
lisse pour tout $t\in T\setminus \{o\}$. Soient aussi 
$r$ sections $\sigma_i : T \to  \mathcal C$ de $\pi$
telles que $\sigma_i(o)=p_i$. On l'a dit, ceci est 
facile \`a construire.

\medskip

{\em Etape 2.}
Si le morphisme $F$ cherch\'e existe, 
alors pour tout $t \in T$, 
$$F_t \in {\rm Mor}(\pi^{-1}(t),X,\forall i \,\, \sigma_i(t)\mapsto f(p_i)).$$
D'apr\`es Mori, les sch\'emas 
${\rm Mor}(\pi^{-1}(t),X,\forall i \,\, \sigma_i(t)\mapsto f(p_i))$
s'assemblent en un $T$-sch\'ema $$\rho : {\rm Mor}(\mathcal C,X,
\forall i \,\, \sigma_i(T)\mapsto f(p_i)) \to T$$ tel que 
${\rm Mor}(\mathcal C,X,\forall i \,\, \sigma_i(T)\mapsto f(p_i))_t \simeq 
{\rm Mor}(\pi^{-1}(t),X,
\forall i \,\,\sigma_i(t)\mapsto f(p_i))$. 

Comme dans la version absolue, il y a un crit\`ere simple
permettant de comprendre le sch\'ema $${\rm Mor}(\mathcal C,X,
\forall i \,\, \sigma_i(T)\mapsto f(p_i))$$
au voisinage d'un point $[h]$~: si $$[h]\in {\rm Mor}(\pi^{-1}(t_0),X,
\forall i\,\, \sigma_i(t_0) \mapsto f(p_i))$$ et si 
$H^1(\pi^{-1}(t_0),h^*T_X \otimes 
{\mathcal O}_{\pi^{-1}(t_0)}(-\sum_{i=1}^r \sigma_i (t_0)))=0$, alors
$\rho$ est un morphisme lisse au point $[h]$. 
En particulier
${\rm Mor}(\pi^{-1}(t_0),X,
\forall i \,\, \sigma_i(t_0) \mapsto f(p_i))$ 
est lisse au point $[h]$
de dimension $$\dim H^0(\pi^{-1}(t_0),h^*T_X \otimes 
{\mathcal O}_{\pi^{-1}(t_0)}(-\sum_{i=1}^r \sigma_i(t_0)))$$ et
${\rm Mor}(\mathcal C,X,
\forall i \,\, \sigma_i(T)\mapsto f(p_i))$
est irr\'eductible  
au point $[h]$ 
de dimension $$\dim H^0(\pi^{-1}(t_0),h^*T_X \otimes 
{\mathcal O}_{\pi^{-1}(t_0)}(-\sum_{i=1}^r \sigma_i(t_0)))+\dim T.$$

En particulier encore, {\bf et c'est ce qu'il faut retenir ici}, 
{\em si $f~: C \to X$ est le morphisme \`a lisser en fixant les $f(p_i)$
et si $H^1(C,f^*T_X \otimes 
{\mathcal O}_{C}(-\sum_i p_i))=0$, 
alors la composante 
de $${\rm Mor}(\mathcal C,X,\forall i \,\, \sigma_i(T)\mapsto f(p_i))$$ 
passant par $[f]$ domine $T$. Il existe
donc une courbe (que l'on peut supposer lisse)
$$T' \to {\rm Mor}(\mathcal C,X,\forall i \,\,\sigma_i(T)\mapsto f(p_i)
)$$ passant par $[f]$ et dominant $T$. 
Le morphisme naturel
$$\mathcal C \times_{T}T' \to X$$ fournit 
le morphisme $F~:\mathcal C \times_{T}T' \to X$ cherch\'e !}
  
\medskip

{\em Etape 3.} Terminons la preuve dans le cas o\`u
il n'y a que deux composantes :
au vu de ce qui pr\'ec\`ede, il suffit de montrer 
que $H^1(C,f^*T_X \otimes 
{\mathcal O}_{C}(-\sum_{i=1}^r p_i))=0$ (par semi-continuit\'e, on aura 
aussi $H^1(\pi^{-1}(t),F_t^*T_X \otimes 
{\mathcal O}_{\pi^{-1}(t)}(-\sum_{i=1}^r \sigma_i (t)))=0$ 
pour tout $t$ quitte \`a r\'etr\'ecir $T$,
ce qui implique que $F_t$ est $(r-1)$-libre).

Si $q = C_1 \cap C_2$, consid\'erons la suite exacte suivante :
$$ 0 \to \mathcal O _{C_2}(-q-p_{r_1+1}-\cdots -p_r) \to 
\mathcal O _C(-p_1-\cdots-p_{r})
\to \mathcal O _{C_1}(-p_1-\cdots-p_{r_1}) \to 0$$
o\`u l'on a r\'eordonn\'e les $p_i$ de sorte que les $r_1$ premiers
soient sur $C_1$.
Apr\`es tensorisation par $f^*T_X$ et passage \`a la suite exacte
longue de
cohomologie associ\'ee, on en d\'eduit 
la suite exacte suivante~:

\begin{align*}
H^1(C_2,f_2^*T_X \otimes \mathcal O _{C_2}(-q-\sum_{i=r_1+1}^{r} p_i)))
& \to H^1(C, f^*T_X \otimes \mathcal O _C(-\sum_{i=1}^r p_i)) \\
&\hspace{2.5cm} \to 
H^1(C_1,f_1^*T_X \otimes \mathcal O _{C_1}(-\sum_{i=1}^{r_1} p_i))).
\end{align*}

Comme $f_1$ est $(r_1-1)$-libre et $f_2$ est $r_2$-libre, 
les deux groupes extr\^emes sont nuls\footnote{Rappelons que
$H^1(\mathbb P ^1, \mathcal O_{\mathbb P ^1}(a))=0$ si et seulement
si $a\geq -1$.},
ce qui fournit le r\'esultat.\finpreuve

\subsubsection{Lissage des peignes}

\begin{defi}
Un peigne rationnel $C$ est un arbre de $m+1$ courbes rationnelles
lisses, avec une composante distingu\'ee $D$ (la poign\'ee) et 
$m$ dents $C_1,\ldots,C_m$ deux \`a deux disjointes, chaque dent
$C_i$ rencontrant
$D$ transversalement en un unique point $q_i:=D\cap C_i$, $i=1,\ldots,m$.
\end{defi}

Dans cette d\'efinition, il y a donc une composante privil\'egi\'ee,
\`a savoir la poign\'ee. Le th\'eor\`eme suivant
permet de lisser des peignes 
rationnels sans que la poign\'ee ne 
soit suppos\'ee libre. C'est une diff\'erence
majeure avec le paragraphe pr\'ec\'edent o\`u toutes les composantes
\'etaient suppos\'ees libres pour permettre le lissage.

\begin{theo}
Soient $X$ une vari\'et\'e projective lisse, $C$ un peigne
rationnel \`a $m$ dents et $f~: C \to X$ un morphisme.
On suppose que la restriction de $f$ \`a chaque dent est libre
et on se donne $r \geq 0$ points lisses $p_1,\ldots,p_r$ de $C$, 
situ\'es sur
la poign\'ee. 
Si
$$ m > K_X\cdot f_*D+(r-1)\dim X + \dim_{[f_{|D}]} {\rm Mor}(\mathbb P^1,
X, \forall i\, \, p_i \mapsto f(p_i)),$$
alors il existe un sous-peigne $C'$ de $C$ avec la m\^eme poign\'ee
et avec au moins une dent
tel que $f_{|C'}$ soit lissable en fixant les $p_i$\footnote{
La preuve de ce th\'eor\`eme ne permet pas de contr\^oler 
la diff\'erence $C \setminus C'$. Dans \cite{GHS03},
les auteurs d\'emontrent un th\'eor\`eme de lissage de peigne 
avec toutes les dents,
pour peu que les dents soient ``g\'en\'erales''. 
Ceci est aussi discut\'e dans le 
cours d'Olivier Wittenberg.}
.
\end{theo}

{\bf D\'emonstration.}

{\em Etape 1.} Soit $\mathcal C_m \to D \times \mathbb A ^m$
l'\'eclatement de $D \times \mathbb A ^m$  
le long des $m$ sous-vari\'et\'es disjointes $\{q_i\} \times \{y_i=0\}$
de codimension $2$
dans $D \times \mathbb A ^m$ (on note 
$E_i \simeq \mathbb P ^1 \times \mathbb A ^{m-1}  
\subset \mathcal C_m$ 
le diviseur exceptionnel
au dessus de $\{q_i\} \times \{y_i=0\}$). 
Soit $\pi :\mathcal C_m \to \mathbb A ^m$ 
la projection induite.
La fibre au-dessus de $0 \in \mathbb A ^m$ est le peigne $C$~;
plus g\'en\'eralement, la fibre de $\pi$
au-dessus de $y \in \mathbb A ^m$ est un peigne de poign\'ee $D$
dont le nombre de dents est \'egal au nombre de $y_i$ \'egaux \`a
z\'ero. Il est important de remarquer le fait facile suivant.
Pour $1\leq m'<m$, si $V_{m'}:=\{ y \in \mathbb A ^m \mid
y_1=\cdots= y_{m'}=0 \}$, alors
l'image inverse $\pi^{-1}(V_{m'})$
poss\`ede $m'+1$ composantes irr\'eductibles de dimension
$m-m'+1$ d\'ecrites ainsi : pour chaque $1 \leq i \leq m'$, 
$E_i \cap \pi^{-1}(V_{m'}) \simeq \mathbb P ^1 \times
\mathbb A ^{m-m'}$ est une composante irr\'eductible
de $\pi^{-1}(V_{m'})$, la derni\`ere est isomorphe \`a la vari\'et\'e
$\mathcal C_{m-m'}$ obtenue
en \'eclatant $D \times V_{m'} \simeq D \times \mathbb A ^{m-m'}$  
le long des $m-m'$ sous-vari\'et\'es disjointes $\{q_i\} \times \{y_i=0\}$
pour $i > m'$. La fibre de $\mathcal C_{m-m'} \to \mathbb A ^{m-m'}$  
au dessus de $0$ est le sous-peigne $C'= D \cup C_{m'+1}\cup \cdots \cup C_m$ 
de $C$.
 
Remarquons qu'il y a aussi $r$ sections
$\sigma_i : \mathbb A ^m \to \mathcal C_m$ de $\pi$ telles
que $\sigma_i(0)=p_i \in D \subset C = \pi^{-1}(0)$.

\medskip

{\em Etape 2.} A nouveau, les sch\'emas 
${\rm Mor}(\pi^{-1}(y),X,\forall i \,\, \sigma_i(y)\mapsto f(p_i))$
s'assemblent en un $ \mathbb A ^m$-sch\'ema 
$${\rm Mor}(\mathcal C_m,X,
\forall i \,\, \sigma_i( \mathbb A ^m)\mapsto f(p_i)) \to \mathbb A ^m$$ 
tel que 
${\rm Mor}(\mathcal C _m,X,\forall i \,\, 
\sigma_i(\mathbb A ^m)\mapsto f(p_i))_y \simeq 
{\rm Mor}(\pi^{-1}(y),X,
\forall i \,\,\sigma_i(y)\mapsto f(p_i))$.
 
Montrons que 
$$ \dim_{[f]} {\rm Mor}(\mathcal C_m,X,
\forall i \,\, \sigma_i( \mathbb A ^m)\mapsto f(p_i)) 
> \dim_{[f]} {\rm Mor}(C,X,\forall i \,\, p_i \mapsto f(p_i)).$$
  
Le membre de droite est facile \`a estimer : un morphisme
de $C$ dans $X$ est d\'etermin\'e par sa restriction \`a chaque composante,
et les morphismes correspondants doivent co\"{\i}ncider aux points 
d'intersections. Pour les dents, les espaces de morphismes 
sont lisses en $[f_{|C_i}]$, d'espace tangent
$H^0(C_i,f_{|C_i}^*T_X \otimes {\mathcal O}_{C_i}(-p_i))$
puisque les dents sont libres.
Or, si $a\geq 0$, $\dim H^0(\mathbb P^1,{\mathcal O}_{\mathbb P ^1}(a-1))=a$,
il vient donc :
$$\dim_{[f]} {\rm Mor}(C,X,\forall i \,\, p_i \mapsto f(p_i))
\leq \sum_{i=1}^m (-K_X \cdot f_* C_i) + 
\dim_{[f_{|D}]} {\rm Mor}(\mathbb P^1,
X, \forall i\, \, p_i \mapsto f(p_i)).$$ 

Le membre de gauche est lui minor\'e par 
$-K_X \cdot f_* C + (1-r)\dim X +m$. L'in\'egalit\'e cherch\'ee 
d\'ecoule donc de l'hypoth\`ese.

\medskip

{\em Etape 3.} Il y a donc une courbe $T$ passant
par $[f]$ dans ${\rm Mor}(\mathcal C_m,X,
\forall i \,\, \sigma_i( \mathbb A ^m)\mapsto f(p_i))$
qui ne s'envoie pas sur $0$ dans $\mathbb A ^m$. 
Si l'image de $T$ rencontre $(k^*)^m \subset \mathbb A ^m$,
c'est gagn\'e : le morphisme $f~: C \to X$ est lissable
sans avoir \`a ``enlever de dents''. Sinon, 
quitte \`a renum\'eroter les coordonn\'ees, on peut supposer 
que $y_1,\ldots,y_{m'}$ sont les coordonn\'ees qui s'annulent
sur l'image de $T$ : $\pi (T) \subset V_{m'}:=\{ y \in \mathbb A ^m \mid
y_1=\cdots= y_{m'}=0 \}$ et $\pi (T)$ rencontre l'ouvert 
$(k^*)^{m-m'} \subset V_{m'} \simeq \mathbb A ^{m-m'}$.
Comme $T$ ne s'envoie pas sur $0$,
on a $m'<m$. On a alors $$T \subset 
{\rm Mor}(\pi^{-1}(V_{m'}),X,\forall i \,\, 
\sigma_i(V_{m'})\mapsto f(p_i)).$$ Or, on a vu
\`a l'\'etape 1 que l'une des composantes de $\pi^{-1}(V_{m'})$ est isomorphe
\`a $\mathcal C _{m'}$ pour $m'<m$ si bien que $T$
fournit un lissage d'un sous-peigne $C'$ de $C$ poss\'edant
$m-m' >0$ dents.\finpreuve

\medskip

{\bf Remarques.} 
\begin{enumerate}
\item[(i)] Il est important de comprendre o\`u l'on a utilis\'e
que les dents sont libres. La preuve ci-dessus montre qu'une
in\'egalit\'e de la forme 
\begin{align*}
m & > \sum_{i=1}^m (
\dim H^0(C_i,f_{|C_i}^* T_X 
\otimes {\mathcal O}_{C_i}(-q_i))-(-K_X\cdot
f_*C_i))\\
&\hspace{1.5cm}+ K_X\cdot f_*D +
(r-1)\dim X + \dim_{[f_{|D}]} {\rm Mor}(\mathbb P^1,
X, \forall i\, \, p_i \mapsto f(p_i))
\end{align*}
suffit pour lisser un sous-peigne, sauf qu'une telle
in\'egalit\'e n'est en g\'en\'eral pas satisfaite si les quantit\'es
positives ou nulles  
$\dim H^0(C_i,f_{|C_i}^* T_X \otimes {\mathcal O}_{C_i}(-q_i))-(-K_X\cdot
f_*C_i)$ ne sont pas nulles, ce que garantit l'hypoth\`ese $C_i$ libre.
\item[(ii)] On n'a pas r\'eellement utilis\'e le fait que la poign\'ee
soit une courbe rationnelle. Le lecteur pourra prouver un 
\'enonc\'e analogue dans le cas d'une poign\'ee quelconque.
Il faut bien
s\^ur adapter la notion de lissage, la fibre g\'en\'erale du lissage 
devant \^etre une 
courbe lisse dont le genre est celui du peigne.

\end{enumerate}

\medskip

Je ne sais pas quelle est l'origine exacte de l'id\'ee illustr\'ee 
par le th\'eor\`eme pr\'ec\'edent, \`a savoir qu'on peut lisser 
toute courbe si on lui
attache suffisamment de courbes libres. Hartshorne et Hirschowitz d\'emontrent
par exemple le th\'eor\`eme suivant \cite{HH83}.

\begin{theo}\label{HH}{\bf (Hartshorne et Hirschowitz)} 
Soit $D$ une courbe lisse dans $\mathbb P ^3$
et soit $$X=D\cup L_1 \cup \cdots \cup L_m$$ un peigne de poign\'ee
$D$ dont les dents $L_i$ sont des droites.
Si $m > \dim H^1(D,N_{D/\mathbb P ^3}) +1$, alors $X$ est lissable
(en une famille de courbes lisses dont le genre d\'epend de celui
de $D$ et de $m$).
\end{theo}

\newpage

\centerline{\bf COURS 4}

\medskip

\begin{center}
\begin{minipage}{130mm}
\scriptsize

Dans ce quatri\`eme cours, 
nous montrons que 
les notions de connexit\'e rationnelle par cha\^{\i}nes
et de connexit\'e rationnelle co\"{\i}ncident 
dans la cat\'egorie des vari\'et\'es projectives lisses
sur un corps de caract\'eristique z\'ero. Plusieurs applications
sont donn\'ees, notamment \`a l'\'etude des vari\'et\'es de Fano. 

\end{minipage}
\end{center}

\medskip

\section{Connexit\'e rationnelle {\em versus} connexit\'e
rationnelle par cha\^{\i}nes. Applications.}

\subsection{Une \'equivalence remarquable}
Voici le r\'e\-sultat principal de cette partie.

\begin{theo}\label{versus} Soit $X$ une vari\'et\'e projective lisse
sur un corps de caract\'eris\-tique z\'ero.
Si $X$ est rationnellement connexe par cha\^{\i}nes,
alors par tout sous-ensemble fini de $X$
passe une courbe rationnelle tr\`es libre (d\'efinie sur une extension
non d\'enombrable du corps de base).
\end{theo}

Le corollaire suivant est tr\`es utile.

\begin{cor}\label{coversus} Soit $X$ une vari\'et\'e projective lisse
sur un corps de caract\'eristique z\'e\-ro.
Alors $X$ est rationnellement connexe si et seulement si
$X$ est rationnellement connexe par cha\^{\i}nes. 
\end{cor}

{\bf D\'emonstration du corollaire.}
Supposons $X$ rationnellement connexe par cha\^{\i}nes.
D'a\-pr\`es le th\'eor\`eme \ref{versus},
il existe un corps $K$ extension non d\'enombrable
de $k$ et une courbe rationnelle tr\`es libre 
$f: \mathbb P^1_K \to X_K$ d\'efinie sur $K$.
La courbe $f$ est d\'efinie sur une extension de type fini
de $k$, {\em i.e.} sur le corps $k(U)$ d'une $k$-vari\'et\'e $U$.
En d'autres termes, il y a une famille param\'etr\'ee par 
$U$ de $k$-courbes rationnelles
tr\`es libres. La vari\'et\'e $X$ est donc ($k$-) 
rationnellement connexe.\finpreuve

\medskip

Nous verrons plus loin que si $X$ est une vari\'et\'e de Fano non singuli\`ere
sur un corps alg\'ebriquement clos de caract\'eristique arbitraire,
alors $X$ est rationnellement connexe par cha\^{\i}nes (ce r\'esultat est
d'ailleurs plus facile en caract\'eristique positive qu'en caract\'eristique
nulle). On en d\'eduit le corollaire suivant.

\begin{cor} Soit $X$ une vari\'et\'e projective lisse
sur un corps de caract\'eristique z\'ero.
Si $X$ est de Fano, alors $X$ est rationnellement connexe.
\end{cor}

La d\'emonstration du th\'eor\`eme utilise le fait 
suivant
(voir sa preuve dans le livre de Debarre, facile une fois qu'on a \'ecrit que
$X^{\rm libre}$ est une intersection d\'enombrable d\'ecroissante d'ouverts).

\medskip

{\bf Fait.} Soit $X$ une vari\'et\'e projective lisse
sur un corps non d\'enombrable.
Soient $\pi : \mathcal C \to T$ un morphisme propre et plat de
dimension relative $1$, o\`u $T$ est irr\'eductible et $F~:\mathcal C \to X$.
Soit $o$ un point distingu\'e de $T$. Si 
l'image de $F_{o}$ rencontre $X^{\rm libre}$, alors 
il existe une famille d\'enombrable d'ouverts
non vides $T_i$ de $T$ telle que pour tout $t \in \cap_i T_i$,
l'image de $F_t~: \pi^{-1}(t) \to X$ rencontre $X^{\rm libre}$.
Autrement dit, il faut retenir le slogan~: 
{\em toute d\'eformation 
tr\`es g\'en\'erale d'une cha\^{\i}ne de courbes rationnelles
rencontrant $X^{\rm libre}$
rencontre aussi $X^{\rm libre}$.}

\medskip 

Evidemment, ce fait serait trivial si $X^{\rm libre}$
\'etait un ouvert de $X$, ce qu'il n'est pas en g\'en\'eral\footnote{Si $S$
est une surface rationnelle avec une infinit\'e d\'enombrable 
de $(-1)$-courbes, $S^{\rm libre}$ est certainement contenu
dans le compl\'ementaire de ces $(-1)$-courbes.}.
En revanche, si $X$ est unir\'egl\'ee sur un corps non d\'enombrable,
alors $X^{\rm libre}$ est dense dans $X$.

\medskip

{\bf D\'emonstration du th\'eor\`eme \ref{versus}.}
C'est tr\`es joli et \c{c}a se fait en plusieurs \'etapes.
Soit $K$ un corps alg\'ebriquement clos non d\'e\-nombrable
extension du corps de base $k$. Dans toute la suite
de cette preuve, on travaille sur $K$.

\medskip 

{\em Etape 1.} D\'emontrons le lemme pr\'eliminaire suivant.

\begin{lemm}\label{libre} Soit $X$ une vari\'et\'e projective.
Supposons qu'il existe $f~: \mathbb P^1 \to f_*(\mathbb P^1) = D  \subset X$ 
une courbe rationnelle sur $X$
et $h~: \mathbb P^1 \to X$ une courbe rationnelle libre {\em dont l'image 
rencontre $X^{\rm libre}$} telles
que $h(0)=x \in D$. Soit $y \in D$ distinct de $x$. 
Alors il existe une courbe rationnelle rencontrant $X^{\rm libre}$
et passant par
$x$ et par $y$.
\end{lemm} 

\medskip 

En effet, comme $h$ est libre, l'application d'\'evaluation 
en $0$, ${\rm ev}_0 : M \to X$ est dominante (o\`u $M$ est la composante de 
${\rm Mor}(\mathbb P^1,X)$ passant
par $[h]$ et o\`u ${\rm ev}_0([g])=g(0)$). Comme son image rencontre $D$, sa 
restriction ${\rm ev}_0: {\rm ev}_0^{-1}(D) \to D$ domine $D$. 
Par le fait ci-dessus,
on en d\'eduit que par un point tr\`es g\'en\'eral de $D$ passe une courbe
libre d\'eformation de $h$ dont l'image rencontre $X^{\rm libre}$.
Si $x_1,\ldots,x_m$ sont tr\`es g\'en\'eraux dans $D$ (donc distincts de
$x$ et $y$), il existe donc une courbe $C_i$ libre passant par $x_i$.
Pour $m$ suffisamment grand, le peigne $D\cup C_1 \cup \cdots \cup C_m$
de poign\'ee $D$ et de dents les $C_i$, 
poss\`ede un sous-peigne $C'$ lissable \`a $x$ et $y$ fix\'es. 
Comme $C'$ poss\`ede au moins une dent, $C'$ rencontre $X^{\rm libre}$,
une d\'eformation tr\`es g\'en\'erale de $C'$ rencontre $X^{\rm libre}$
donc est libre en appliquant \`a nouveau le fait ci-dessus.\finpreuve

\medskip

{\em Etape 2.} Soient $(x_1,x_2)\in X^{\rm libre} \times X$
et $C$ une cha\^{\i}ne de courbes rationnelles joignant $x_1$
\`a $x_2$. Quitte \`a supprimer et renum\'eroter des maillons, on peut supposer
que
$C=C_1\cup \cdots \cup C_k$ avec $x_1 \in C_1$, $x_2 \in C_k$
et $C_{i-1}\cap C_i \neq \emptyset$ pour $i=1,\ldots,k-1$.
Comme $x_1 \in X^{\rm libre}$, $C_1$ est libre. Choisissons un point $x$
dans $C_1\cap C_2$ et un point $y$
dans $C_2\cap C_3$ (on peut supposer que 
$y\neq x$ car si $C_1\cap C_2 = C_2\cap C_3$, on supprime alors
la courbe $C_2$ de la cha\^{\i}ne $C$).
Le lemme \ref{libre}
montre qu'il y a une courbe rencontrant $X^{\rm libre}$ et passant
par $x$ et $y$~: on a ainsi remplac\'e $C_2$ par une courbe
rencontrant $X^{\rm libre}$, 
d\'eformation d'un peigne de poign\'ee $C_2$, dans la cha\^{\i}ne $C$. 
En it\'erant ce proc\'ed\'e, on construit une cha\^{\i}ne,
passant par $x_1$ et $x_2$,
de courbes rationnelles rencontrant toutes $X^{\rm libre}$.

\medskip

{\em Etape 3.} Lissons la cha\^{\i}ne de courbes rationnelles
obtenue \`a l'\'etape pr\'ec\'edente, en fixant $x_2$~; c'est possible gr\^ace
au th\'eor\`eme de lissage des arbres de courbes rationnelles libres.

\medskip
 
{\em Etape 4.} 
Le bilan des trois \'etapes pr\'ec\'edentes est le suivant~: si $X$
est rationnellement connexe par cha\^{\i}nes, pour tout point
$x_2 \in X$, l'image de 
$$ {\rm ev}^1 : \mathbb P^1 \times 
{\rm Mor}(\mathbb P^1,X,0\mapsto x_2) \to X$$
contient $X^{\rm libre}$ dans son adh\'erence. Comme $X^{\rm libre}$
est dense, cette application est dominante et
il y a une composante $M$ de 
${\rm Mor}(\mathbb P^1,X,0\mapsto x_2) \to X$
tel que ${\rm ev}^1 : \mathbb P ^1 \times M \to X$ soit dominante.
Soit $[f] \in M$ une courbe libre. La composante $M$
est lisse au point $[f]$ et comme les d\'eformations de
$f$ fixant $x_2$ dominent $X$, on en d\'eduit que $f$ est tr\`es libre.
R\'ep\'etons cet argument qui nous a d\'ej\`a servi dans la 
preuve du th\'eor\`eme 8 et qui, je l'esp\`ere, est devenu familier au lecteur.

\medskip

En effet, il y a
(comme la caract\'eristique du corps de base est suppos\'ee nulle) 
$(u,[f]) \in \PP^1 \times M$
tel que la diff\'erentielle de ${\rm ev}^1$ en $(u,[f])$
soit surjective.
A nouveau, la diff\'erentielle de ${\rm ev}^1$ en $(u,[f])$ 
est l'application naturelle 
$$ T_{u}\mathbb P ^1 
\oplus H^0 (\mathbb P ^1,f^*T_X \otimes {\mathcal O}_{\mathbb P^1}(-1))
\to T_{f(u)}X$$
qui \`a $(v,\sigma) \in T_{u}\mathbb P ^1 
\oplus H^0 (\mathbb P ^1,f^*T_X \otimes {\mathcal O}_{\mathbb P^1}(-1))$ 
associe 
$(Tf)_{u}(v)+\sigma(u)$.  
Le point cl\'e est le suivant~: 
comme $T_{\mathbb P ^1}\simeq \mathcal O _{\mathbb P^1} (2)$,
l'application 
$H^0(\mathbb P ^1,T_{\mathbb P ^1})\to
T_{u}\mathbb P ^1$ est surjective, 
donc
l'image de l'application
$H^0(\mathbb P ^1,f^*T_X) \to T_{f(u)}X$
contient celle
de $(Tf)_{u} : 
T_{u}\mathbb P ^1 \to 
T_{f(u)}X$, d'o\`u l'on d\'eduit \'evidemment que  
$H^0 (\mathbb P ^1,f^*T_X \otimes {\mathcal O}_{\mathbb P^1}(-1))
\to T_{f(u)}X$ est surjectif. 
De l\`a, si $f^*T_X \simeq \oplus_{i=1}^n \mathcal O _{\mathbb P^1} (a_i)$,
chaque 
$$H^0(\mathbb P ^1,\mathcal O _{\mathbb P^1} (a_i-1)) \to 
\mathcal O _{\mathbb P^1} (a_i)_{u}$$ est surjectif, donc $a_i \geq 1$
pour tout $i$ : $f$ est tr\`es libre.

\medskip

{\em Etape 5.} Montrons que par deux points quelconques $x$ et $y$
passe
une courbe tr\`es libre. 

\medskip

C'est facile, il suffit de 
joindre $x$ et $y$ \`a un point $z$ tr\`es g\'en\'eral
\`a l'aide de courbes tr\`es libres, puis de lisser \`a l'aide du th\'eor\`eme
de lissage des arbres. 

\medskip

{\em Etape 6.} Montrons finalement que par tout sous-ensemble 
fini $\{x_1,\ldots,x_m\}$ de $X$
passe une courbe rationnelle tr\`es libre. 

\medskip 

On le fait par r\'ecurrence sur $m$, on vient de voir le 
cas $m=2$. Soit alors $f~: \mathbb P^1 \to X$ une courbe tr\`es libre
passant par $x_1,\ldots,x_{m-1}$ et
soit une courbe tr\`es libre passant par $x_m$ et un point auxiliaire
de $f(\mathbb P^1)$.
Quitte \`a composer $f$ \`a la source, on peut supposer que $f$
est $r$-libre avec $r$ suffisamment grand pour qu'on puisse lisser
cet arbre \`a deux composantes en fixant tous les $x_i$.\finpreuve

\medskip

En \'etant un peu plus soigneux, il est possible de borner
le degr\'e des courbes tr\`es libres ainsi construites en fonction du
nombre de maillons n\'ecessaires pour joindre deux points g\'en\'eraux de
$X$ par une cha\^{\i}ne de courbes rationnelles\footnote{Ce type de r\'esultat
permet par exemple de montrer que sur un corps non d\'enombrable, si deux
points g\'en\'eraux sont joignables par une 
cha\^{\i}ne de courbes rationnelles,
alors deux points quelconques peuvent \^etre joints
par une cha\^{\i}ne de courbes rationnelles.}. On renvoie \`a nouveau 
au livre d'O. Debarre.

\subsection{Deux applications}

\subsubsection{}
On montre que sur un corps de caract\'eristique z\'ero, la connexit\'e
rationnelle est une propri\'et\'e ouverte et ferm\'ee des vari\'et\'es
projectives lisses.

\begin{theo}\label{rcdef} On suppose que le corps de base est de
caract\'eristique z\'ero.
Soit $T$ une vari\'et\'e quasi-projective 
et soit $\pi : \mathcal X \to T$ un morphisme
projectif lisse. S'il existe $t_0 \in T$ tel
que $X_{t_0}:= \pi^{-1}(t_0)$ soit rationnellement connexe, alors
$X_{t}$ est rationnellement connexe
pour tout $t\in T$.
\end{theo} 

{\bf D\'emonstration.} 

{\em Etape 1.} La connexit\'e rationnelle est une propri\'et\'e ouverte.

\medskip  

Soit en effet $f~:\mathbb P^1 \to X_{t_0}$ une courbe tr\`es libre
de $X_{t_0}$. On voit aussi $f$ comme une courbe libre de $\mathcal X$.
Le sch\'ema ${\rm Mor}( \mathbb P^1,\mathcal X)$ est lisse
au point $[f]$, les d\'eformations de $f$ dominent $\mathcal X$
et sont verticales par le lemme de rigidit\'e : si $h$ est une d\'eformation
de $f$, son image est contenue dans une fibre $X_{t}$ et $h$ est alors
une courbe tr\`es libre de $X_t$ par un argument de semi-continuit\'e
de la cohomologie par exemple.

\medskip 

{\em Etape 2.} 
La connexit\'e rationnelle est une propri\'et\'e ferm\'ee.

\medskip 

On suppose que $T$ est une courbe et que 
$X_{t}$ est rationnellement connexe pour $t  \neq t_{\infty}$
Si $x_{\infty}$ et $y_{\infty}$ sont deux points de $X_{t_{\infty}}$,
on les approche par $x_t$ et $y_t$ dans $X_t$. Une courbe rationnelle
$C_t$ contenue dans $X_t$ et passant par $x_t$ et $y_t$ 
d\'eg\'en\`ere en une cha\^{\i}ne de courbes rationnelles 
contenue dans $X_{t_{\infty}}$ et passant par $x_{\infty}$ et $y_{\infty}$.
Autrement dit, $X_{t_{\infty}}$ est rationnellement connexe par cha\^{\i}nes,
donc rationnellement connexe par le th\'eor\`eme 
pr\'ec\'edent.\footnote{
Plus efficacement, il suffit de dire que de fa\c{c}on g\'en\'erale,
la connexit\'e
rationnelle par cha\^{\i}nes est une propri\'et\'e ferm\'ee.}\finpreuve

\subsubsection{Simple-connexit\'e}

Dans ce paragraphe, on travaille sur $\mathbb C$.
Notre objectif est de d\'emontrer que les vari\'et\'es lisses
rationnellement connexes sont simplement connexes.

La premi\`ere \'etape est de d\'emontrer que si $X$ est lisse
et rationnellement connexe par cha\^{\i}nes, alors son groupe fondamental
est fini. On utilise ici 
une jolie astuce due \`a Campana. Elle repose sur un lemme 
de g\'eom\'etrie diff\'erentielle que nous admettrons ici.

\begin{lemm}
Soit $f~: V \to W$ un morphisme dominant avec $W$ normale.
Alors l'image du morphisme induit $\pi_1(f)~: \pi_1(V) 
\to \pi_1(W)$ est d'indice
fini dans $\pi_1(W)$. En particulier, si cette image est triviale,
alors $\pi_1(W)$ est fini. 
\end{lemm} 

Soient $X$ une vari\'et\'e rationnellement connexe et
$f~: \mathbb P^1 \to X$ une courbe tr\`es libre. 
On a vu que les d\'eformations de $f$ \`a point fix\'e dominent
$X$, autrement dit, si $M\subset {\rm Mor}(\mathbb P^1,X,0\mapsto f(0))$
est la composante passant par $f$, 
l'application d'\'evaluation ${\rm ev}^1 :
\mathbb P^1 \times M \to X$ est dominante.
L'injection $\iota : \{0\}\times M \hookrightarrow \mathbb P^1 \times M$
induit un isomorphisme au niveau des groupes fondamentaux.
On en d\'eduit que l'image de $\pi_1({\rm ev}^1)~: \pi_1(\mathbb P^1 \times M)
\to \pi_1(X)$ est \'egale \`a l'image de la compos\'ee
$\pi_1({\rm ev}^1 \circ \iota)~:  \pi_1(\{0\}\times M) \to \pi_1(X)$
qui est triviale puisque ${\rm ev}^1 \circ \iota$ est constant \'egal
\`a $f(0)$ !
Par le lemme, $\pi_1 (X)$ est fini.

\begin{prop}\label{rcannul} Soit $X$ une 
vari\'et\'e projective lisse rationnellement connexe sur
$\mathbb C$.
Alors $$H^0(X,(\wedge ^m T_X ^*)^{\otimes p})=0$$ 
pour tous $m$ et $p$ strictement positifs.
En particulier, $\chi (X,\mathcal O _X)=1$.
\end{prop}

{\bf D\'emonstration.}
Si $f$ est une courbe tr\`es libre, $f^*(\wedge ^m T_X ^*)^{\otimes p})$ 
est une somme de $\mathcal O _{\mathbb P ^1}(b_i)$ avec tous les
$b_i$ strictement n\'egatifs. On a donc 
$$  H^0(\mathbb P ^1,f^*(\wedge ^m T_X ^*)^{\otimes p}) =0,$$
donc toute section de $(\wedge ^m T_X ^*)^{\otimes p}$ s'annule
sur l'image de $f$, donc est identiquement nulle
puisque les courbes tr\`es libres dominent $X$.
Enfin, on a $$\dim H^j(X,\mathcal O _X) = 
\dim H^0(X,\wedge ^j T_X ^*)=0\footnote{
L'\'egalit\'e $\dim H^j(X,\mathcal O _X) = 
\dim H^0(X,\wedge ^j T_X ^*)$
est fausse en caract\'eristique positive. En caract\'eristique z\'ero,
cette \'egalit\'e provient de la th\'eorie de Hodge.}$$
si $j>0$, donc $\chi (X,\mathcal O _X)=
\dim H^0(X,\mathcal O _X)=1$\footnote{Exercice : adapter la preuve
pr\'ec\'edente pour montrer que si $X$ est une 
vari\'et\'e projective lisse unir\'egl\'ee sur un corps
de caract\'eristique z\'ero, alors
$H^0(X,mK_X)=0$ pour tout $m>0$.}.\finpreuve

\medskip

On peut maintenant prouver le th\'eor\`eme annonc\'e.

\begin{theo}\label{simply}
Soit $X$ une vari\'et\'e lisse projective et rationnellement
connexe par cha\^{\i}nes sur $\mathbb C$. Alors $X$ est simplement connexe.
\end{theo}

{\bf D\'emonstration.} Comme $\pi_1 (X)$ est fini, 
le rev\^etement universel $\rho~: \tilde X \to X$ est encore
une vari\'et\'e projective lisse rationnellement connexe\footnote{Une courbe
rationnelle tr\`es libre sur $X$ se rel\`eve en une 
courbe rationnelle tr\`es libre
sur $\tilde X$ car $\mathbb P^1$ est simplement connexe.}.
Comme $\rho$ est \'etale, on a $\rho^* T_X= T_{\tilde X}$
et $\mathcal O _{\tilde X}=\rho^* \mathcal O _{X}$.
On a donc, si $n:=\dim X =\dim \tilde X$, d'apr\`es le
th\'eor\`eme de Riemann-Roch,
\begin{align*}
1= \chi (\tilde X,\mathcal O _{\tilde X}) = 
\deg({\rm td}(T_{\tilde X}))_n 
&= \deg({\rm ch}(\mathcal O _{\tilde X})
\cdot
{\rm td}(T_{\tilde X}))_n 
 = 
\deg({\rm ch}(\rho^* \mathcal O _{X})
\cdot
{\rm td}(\rho ^*T_{X}))_n \\ 
&= 
 \deg(\rho^* ({\rm ch}(\mathcal O _{X})
\cdot
{\rm td}(T_{X})))_n \\
& = \deg(\rho)\cdot  \deg 
({\rm ch}(\mathcal O _X)
\cdot
{\rm td}(T_{X}))_n \\
& = \deg(\rho) \cdot \chi(X,\mathcal O _X) = \deg(\rho),
\end{align*}
donc $\rho$
est un isomorphisme, $\tilde X \simeq X$, la vari\'et\'e $X$ est 
donc simplement connexe.\finpreuve

\medskip

{\bf Remarques.} 
\begin{enumerate}
\item Les deux ingr\'edients de la preuve ci-dessus
sont le fait que $\pi_1(X)$ est fini et que $\chi (X,\mathcal O _X)=1$.
Dans le cas o\`u $X$ est Fano, la premi\`ere assertion 
se d\'emontre aussi en munissant
$X$ d'une m\'etrique riemannienne \`a courbure de Ricci strictement 
positive\footnote{C'est possible et non trivial : 
ceci d\'ecoule de la r\'esolution
de la conjecture de Calabi-Yau.}
et il est bien connu que les vari\'et\'es riemanniennes compactes
\`a courbure de Ricci
positive ont un $\pi_1$ fini. La deuxi\`eme assertion
d\'ecoule imm\'ediatement du th\'eor\`eme d'annulation de Kodaira.
\item En caract\'eristique z\'ero, la proposition \ref{rcannul} reste vraie
et la d\'emonstration 
du th\'eo\-r\`e\-me \ref{simply} montre que tout rev\^etement \'etale fini
de $X$ est trivial : on dit que $X$ est {\em alg\'e\-bri\-quement simplement  
connexe}. Koll\'ar a montr\'e qu'en caract\'eristique positive,  
une vari\'et\'e projective lisse s\'eparablement rationnellement connexe
est alg\'ebriquement simplement  
connexe (on renvoie au survol \cite{Cha03} et sa bibliographie
pour une discussion
en toute carac\-t\'eristique de la simple-connexit\'e d'une vari\'et\'e 
s\'eparablement rationnellement connexe ou unirationnelle).
\end{enumerate}

\section{Connexit\'e rationnelle des vari\'et\'es de Fano.}

Les vari\'et\'es de Fano sont les vari\'et\'es projectives lisses $X$ pour 
lesquelles le fibr\'e anticanonique $-K_X = \det (T_X)$ est ample. 
Elles sont beaucoup \'etudi\'ees car ce sont des briques \'el\'ementaires du
programme de Mori (m\^eme s'il faut alors autoriser des singularit\'es). 
De fa\c{c}on g\'en\'erale, il y a un espoir 
de classifier les vari\'et\'es pour lesquelles
le fibr\'e tangent a tendance \`a \^etre positif. Un \'enonc\'e 
particuli\`erement remarquable est la 
r\'esolution par Mori d'une conjecture de Hartshorne \cite{cass}.

\begin{theo} {\bf (Mori)} Soit
$X$ une vari\'et\'e projective lisse de dimension $n$
sur un corps de caract\'eristique quelconque. 
Alors $T_X$ est ample
si et seulement si $X \simeq \mathbb P ^n$.
\end{theo}

L'hypoth\`ese qu'une vari\'et\'e 
est de Fano est beaucoup plus faible : l'hypoth\`ese
d'amplitude porte sur le d\'eterminant du fibr\'e tangent. Du point
de vue de la g\'eom\'etrie 
riemannienne, c'est la diff\'erence entre une hypoth\`ese portant sur la
courbure sectionnelle et une hypoth\`ese portant sur la
courbure de Ricci. 
Cependant, on a encore un th\'eor\`eme de 
finitude pour les vari\'et\'es de Fano,
aboutissement des travaux de Nadel, Campana, 
Koll\'ar, Matsusaka, Miyaoka et Mori.

\begin{theo}\label{limite}
Sur un corps de caract\'eristique z\'ero, pour tout $n\geq 1$,
il n'y a qu'un nombre fini
de types de d\'eformation de vari\'et\'es lisses de Fano de dimension $n$.
\end{theo}

Il n'est pas question de d\'emontrer ces deux r\'esultats ici, mentionnons
que leur preuve utilise de fa\c{c}on essentielle la g\'eom\'etrie des courbes
rationnelles dans les vari\'et\'es de Fano. Remarquons aussi que l'hypoth\`ese
$X$ Fano n'est \'evidemment pas invariante par morphisme 
birationnel\footnote{Un
peu de publicit\'e pour mes travaux : le comportement par \'eclatement
des vari\'et\'es de Fano reste un sujet actuel d'\'etude initi\'e
dans \cite{Wis91}. Lorsqu'on \'eclate
des points, on a un th\'eor\`eme de classification dont une 
cons\'equence amusante est
le th\'eor\`eme suivant \cite{BCW01}.
\begin{theo} {\bf (Bonavero, Campana et Wi\'sniewski)}
Si $X$ est une vari\'et\'e
projective lisse complexe de dimension $n \geq 3$ et s'il
existe deux points distincts $a$ et $b$ de $X$ tels que 
l'\'eclatement de $X$ de
centre $\{a,b\}$ soit une vari\'et\'e de Fano, 
alors $X$ est isomorphe \`a une quadrique
lisse de $\mathbb P ^{n+1}$.
\end{theo}}.

\begin{ex}
La seule courbe de Fano est la droite projective. En dimension deux,
les surfaces de Fano sont aussi appel\'ees surfaces de Del Pezzo,
ce sont 
$\mathbb P^2$ \'eclat\'e en au plus $8$ points en position g\'en\'erale
et $\mathbb P^1 \times \mathbb P^1$. 
En dimension trois, il y a $105$ familles classifi\'ees par Iskovskikh
quand le nombre de Picard est \'egal \`a $1$ et par Mori
et Mukai quand le nombre de Picard est $\geq 2$ (voir \cite{IsP99}
et son impressionnante bibliographie). 
\end{ex}

\subsection{Le quotient rationnel}
Si $X$ est une vari\'et\'e projective, on peut d\'efinir une relation
d'\'equivalence sur $X$ de la fa\c{c}on suivante~: deux points $x$ et
$x'$ sont \'equivalents s'il existe une cha\^{\i}ne de courbes rationnelles
passant par $x$ et $x'$. Il n'y a aucun espoir que l'espace quotient
pour cette relation d'\'equivalence puisse \^etre muni d'une structure de
vari\'et\'e alg\'ebrique (de sorte que l'application de passage au quotient
soit un morphisme alg\'ebrique). Un \'enonc\'e remarquable d\^u \`a Campana
et Koll\'ar
permet d'obtenir un ``quotient rationnel''. 
Nous \'enon\c{c}ons sans preuve ce
r\'esultat, cas particulier 
d'un \'enonc\'e beaucoup plus g\'en\'eral o\`u
Campana, en caract\'eristique z\'ero, 
consid\`ere les relations d'\'equivalences 
engendr\'ees par des familles
couvrantes de sous-vari\'et\'es (on autorise des cha\^{\i}nes de 
sous-vari\'et\'es qui ne sont pas des courbes ! ).

\begin{theo}{\bf (Campana, Koll\'ar)} 
Soit $X$ une vari\'et\'e projective normale sur un corps 
(alg\'e\-bri\-quement clos)
de 
caract\'e\-ristique 
quelconque. Alors il existe une vari\'et\'e projective normale
$R(X)$, unique \`a application birationnelle pr\`es,
(appel\'ee {\em quotient rationnel}) et
une application rationnelle $\rho : X \dashrightarrow R(X)$
dominante telles que~:
\begin{enumerate}
\item il existe un ouvert non vide $X^0$ de $X$ et un ouvert
non vide $R(X)^0$ de $R(X)$ telle que $\rho$ se restreigne
en un morphisme propre et surjectif $\rho : X^0 \to R(X)^0$ (on
dit que $\rho$ est {\em presque holomorphe}),
\item pour tout $y \in R(X)^0$, la fibre $\rho^{-1}(y)$
est rationnellement connexe par cha\^{\i}nes,
\item toute application presque holomorphe $\varphi~: X 
\dashrightarrow Y$ dont la fibre g\'en\'erale est 
rationnellement connexe par cha\^{\i}nes se factorise par $\rho$
en $\psi : Y \dashrightarrow R(X)$,
\item toute courbe rationnelle de $X$ rencontrant une fibre
tr\`es g\'en\'erale de $\rho$ est contenue dans une fibre de $\rho$.
\end{enumerate}
\end{theo}

Les trois derniers points signifient qu'une fibre tr\`es g\'en\'erale
de $\rho$ est une classe d'\'equiva\-lence.
De ce point de vue, l'\'enonc\'e est optimal : il y a des
surfaces projectives non unir\'egl\'ees contenant une infinit\'e au plus
d\'enombrable de courbes rationnelles (certaines surfaces $K3$ par exemple).

\medskip

Evidemment, si $X$ est une vari\'et\'e projective normale,
alors $X$ est rationnellement connexe par cha\^{\i}nes si
et seulement si $R(X)$ est r\'eduit \`a un point (par d\'efinition).
De m\^eme, $X$ est unir\'egl\'ee si et seulement si $\dim R(X) < \dim X$.

\medskip

Parmi les nombreuses cons\'equences du th\'eor\`eme de Graber,
Harris et Starr, on obtient les deux r\'esultats suivants, le
premier \'etait conjectur\'e par 
Koll\'ar.

\begin{theo}\label{quotientuni} 
Soit $X$ une vari\'et\'e projective lisse sur un corps de caract\'eristique 
quelconque. 
Alors le quotient rationnel $R(X)$
de $X$ n'est pas unir\'egl\'e.
\end{theo}

{\bf D\'emonstration.} C'est imm\'ediat : si $R(X)$
est unir\'egl\'e, soit $C$ une courbe rationnelle passant
par un point g\'en\'eral de $R(X)$. 
Alors, d'apr\`es le th\'eor\`eme
\ref{theoGHS}, 
$\rho : \rho^{-1}(C) \to C$ poss\`ede une section, qui est donc
une courbe rationnelle rencontrant une fibre
g\'en\'erale de $\rho$ non contenue dans une fibre de 
$\rho$, contradiction.\finpreuve

\begin{cor}\label{critererc}
Soit $X$ une vari\'et\'e projective sur un corps de caract\'eristique 
z\'ero. Alors les assertions suivantes
sont \'equivalentes~:
\begin{enumerate}
\item pour toute application rationnelle
$X \dashrightarrow Y$ dominante, $Y$ est unir\'egl\'ee ou r\'eduite
\`a un point,
\item la vari\'et\'e $X$ est rationnellement connexe.
\end{enumerate}
\end{cor} 

{\bf D\'emonstration.} (1) implique (2) : soient $\pi : \tilde X \to X$ une
r\'esolution des singularit\'es de $X$, $\rho : \tilde X \dashrightarrow
R(\tilde X)$ son quotient rationnel. 
Alors, l'application rationnelle 
$\rho \circ \pi^{-1} : X \dashrightarrow R(\tilde X)$ est dominante,
donc $R(\tilde X)$ est r\'eduit \`a un point, donc $\tilde X$
est rationnellement connexe par cha\^{\i}nes, donc
rationnellement connexe par le th\'eor\`eme \ref{versus} 
et par suite
$X$ est rationnellement connexe.
  
(2) implique (1) est \'evident.\finpreuve 

\subsection{Courbes rationnelles sur les vari\'et\'es de Fano}
Ce paragraphe est une deuxi\`eme introduction \`a la th\'eorie
de Mori, je ne donne aucune preuve. On a d\'ej\`a vu que si $K_X$ n'est 
pas nef, la vari\'et\'e $X$ poss\`ede des courbes rationnelles. Le th\'eor\`eme
suivant pr\'ecise cette id\'ee.

\begin{theo} {\bf (Miyaoka et Mukai)} Soit $X$ une vari\'et\'e projective 
sur un corps de caract\'eristique quelconque. Soit
$C$ une courbe de $X$. On suppose que $X$ est lisse le long
de $C$ et que $K_X \cdot C <0$. Soit $H$ un diviseur ample sur $X$.
Alors par tout point $x$ de $C$ passe une courbe rationnelle
$\Gamma$ telle que 
$$0< H \cdot \Gamma \leq 2\dim X \frac{H \cdot C}{-K_X \cdot C}.$$ 
\end{theo}

{\bf Remarques.}

\begin{enumerate}
\item[(i)] Si $X$ est lisse,
on peut \`a l'aide du lemme de cassage obtenir une courbe rationnelle
$\Gamma$ satisfaisant de plus
$-K_X \cdot \Gamma \leq \dim X+1$.
Il s'agit du premier pas (mais le pas essentiel)
pour d\'emontrer le th\'eor\`eme du c\^one dans le cas lisse.

\item[(ii)] Si $X$ est normale 
$\mathbb Q$-Fano\footnote{On entend par l\`a que $X$ est singuli\`ere
et qu'un multiple entier de $-K_X$ est un diviseur de Cartier ample.},
on en d\'eduit que $X$ est couverte par des courbes rationnelles
$\Gamma$ v\'erifiant $-K_X \cdot \Gamma \leq 2\dim X$.
Comme il n'y a qu'un nombre fini de composantes 
de ${\rm Mor}_d(\mathbb P^1,X)$ de degr\'e $d$ donn\'e, on en d\'eduit
que {\em toute vari\'et\'e normale projective 
$\mathbb Q$-Fano
est unir\'egl\'ee}.

\end{enumerate}

\medskip 

Ce th\'eor\`eme est, on l'a d\'ej\`a 
mentionn\'e, plus facile en caract\'eristique 
positive~: on montre qu'un grand multiple de $C$ obtenu \`a l'aide
du morphisme de Frobenius se d\'eforme en passant par $x$ et d\'eg\'en\`ere 
en une cha\^{\i}ne de courbes rationnelles. En caract\'eristique z\'ero,
on utilise le r\'esultat connu en caract\'eristique positive
avec un contr\^ole uniforme du degr\'e des courbes rationnelles produites.
Cette id\'ee g\'eniale, parmi d'autres, a valu \`a Mori la m\'edaille Fields 
et on peut mentionner
qu'on ne conna\^{\i}t pas de preuve alternative 
par des m\'ethodes transcendantes, 
ingr\'edient n\'ecessaire \`a une \'eventuelle extension du programme
de Mori aux vari\'et\'es k\"ahl\'eriennes compactes.

\medskip

Campana d'une part et Koll\'ar, Miyaoka et Mori 
d'autre part ont \'etendu le r\'esultat de  Miyaoka et Mori
de la fa\c{c}on suivante.

\begin{theo} Soit $X$ une vari\'et\'e projective normale 
$\mathbb Q$-Fano sur un corps de 
caract\'eris\-tique quelconque et $\rho : X \dashrightarrow Y$
une application presque holomorphe dominante. Soit $F$ une fibre
g\'en\'erale de $\rho$. Si $Y$ n'est pas r\'eduit \`a un point, alors
il existe une courbe rationnelle rencontrant $F$ et non contenue
dans $F$.
\end{theo} 

\subsection{Connexit\'e rationnelle des vari\'et\'es de Fano}

Il suffit d'assembler les r\'esultats pr\'ec\'edents pour obtenir le
r\'esultat g\'en\'eral suivant.

\begin{theo}\label{fanorc} Soit $X$ une vari\'et\'e projective normale 
$\mathbb Q$-Fano sur un corps de 
caract\'eris\-tique quelconque. Alors $X$ est rationnellement connexe
par cha\^{\i}nes.
\end{theo}

Le corollaire \ref{coversus} implique alors le r\'esultat
annonc\'e.

\begin{cor}
Soit $X$ une vari\'et\'e projective lisse
et de Fano sur un corps de 
caract\'eris\-tique z\'ero. Alors $X$ est rationnellement connexe.
\end{cor}

En \'etant plus soigneux, il est possible \`a nouveau de borner
le degr\'e (et le nombre de maillons) d'une cha\^{\i}ne de courbes
rationnelles permettant de joindre deux points g\'en\'eraux, et m\^eme
quelconques, de $X$. Ce type de bornes dans le cas lisse
permet aussi de montrer
le th\'eor\`eme \ref{limite}.

\newpage

\centerline{\bf COURS 5
(cette partie a \'et\'e \'ecrite avec St\'ephane DRUEL). 
}

\medskip

\begin{center}
\begin{minipage}{130mm}
\scriptsize

Dans ce dernier cours, 
nous pr\'esentons la preuve de la conjecture de connexit\'e rationnelle de 
Shokurov suivant Hacon et ${\rm M^c}$Kernan. Outre ses liens \'evidents
avec les cours pr\'ec\'edents, cette partie se veut \^etre, modestement,
une invitation \`a la g\'eom\'etrie birationnelle moderne,
celle des paires et du MMP. La litt\'erature sur le sujet 
a explos\'e ces 20 derni\`eres ann\'ees, on recommande particuli\`erement
les textes ``historiques'' \cite{KMM87} et \cite{CKM88},
les ouvrages r\'ecents \cite{Cor07}, \cite{Kol97}, \cite{KM98}, \cite{Laz04}
et \cite{Mat01} ainsi que le r\'ecent texte de synth\`ese 
\cite{Dru08}.

\end{minipage}
\end{center}

\section{La conjecture de connexit\'e rationnelle de 
Shokurov.}

Dans toute cette partie, {\bf le corps de base $K$
est alg\'ebriquement clos de caract\'eris\-tique z\'ero}.
On pr\'esente, suivant Hacon et ${\rm M^c}$Kernan, 
une vaste g\'en\'eralisation
de certains des \'enonc\'es pr\'ec\'edemment rencontr\'es
(on pense au lemme d'Abhyankar et au th\'eor\`eme \ref{fanorc}),
le th\'eor\`eme de Graber, Harris et Starr y
jouant un r\^ole essentiel.

\subsection{Les r\'esultats}

\begin{defi} Une {\em paire} $(X,\Delta)$ est la donn\'ee~:
\begin{enumerate} 
\item[(i)] d'une vari\'et\'e alg\'ebrique normale $X$,
\item[(ii)] d'un $\mathbb Q$-diviseur de Weil
$\Delta=\sum_{i=1}^N d_i \Delta_i$, o\`u les $\Delta_i$ sont 
des diviseurs de Weil irr\'eductibles et 
r\'eduits et les $d_i$ des nombres rationnels $\geq 0$,
\end{enumerate}
tels que $K_X + \Delta$ est $\mathbb Q$-Cartier.
\end{defi}

\begin{defi} Soit $(X,\Delta)$ une paire. Une {\em log-r\'esolution 
de $(X,\Delta)$} est un morphisme birationnel propre
$\pi~: Y \to X$ telle que $Y$ est lisse et  
${\rm Exc}(\pi)+\Delta '$
est un diviseur \`a croisements normaux simples
($\Delta '$ est la transform\'ee stricte de $\Delta$). 

Pour un tel $\pi$, on \'ecrit alors 
$$K_Y + \Delta '= \pi^*(K_X + \Delta) + \sum_E a_{\pi}(E,\Delta)E$$
o\`u les $a_{\pi}(E,\Delta) \in \QQ$ et o\`u la somme porte sur les 
diviseurs exceptionnels et irr\'eductibles $E$ de $\pi$. 

\end{defi}

Tout diviseur exceptionnel et irr\'eductible $E$
de $\pi$ d\'efinit une valuation divisorielle $\nu _E$ de $K (X)$,
de centre $\pi(E) \subset X$.
On montre que $a_{\pi}(E,\Delta)$ ne d\'epend que de $\nu _E$ 
et on note alors $a(E,\Delta)=a_{\pi}(E,\Delta)$. 
On \'etend $a(-,\Delta)$ \`a tout diviseur irr\'eductible non exceptionnel
$E$
en prenant l'oppos\'e de la multiplicit\'e de $E$ dans $\Delta$. 

\begin{defi} On dit que la paire 
$(X,\Delta)$ est Kawamata log terminale (klt en abr\'eg\'e) si
$$a(E,\Delta) > -1$$ pour toute valuation divisorielle $\nu_E$\footnote{Il
suffit de le v\'erifier sur une log-r\'esolution.}.
Si $(X,\Delta)$ est une paire et si $\pi$ est une log-r\'esolution
de $(X,\Delta)$, on note 
$${\rm Nklt}(X,\Delta) = 
\cup_{ \{E \mid a(E,\Delta) \leq -1 \} } \pi (E)$$
o\`u la r\'eunion porte sur tous les diviseurs irr\'eductibles.
Le lieu ${\rm Nklt}(X,\Delta)$ ne d\'epend pas de $\pi$,
il est ferm\'e et la paire $(X,\Delta)$ est klt sur
l'ouvert $X\setminus {\rm Nklt}(X,\Delta)$.
\end{defi}

\begin{defi} On dit que la paire 
$(X,\Delta)$ est log canonique (lc en abr\'eg\'e) si
$$a(E,\Delta) \geq -1$$ pour toute valuation divisorielle $\nu_E$\footnote{Il
suffit de le v\'erifier sur une log-r\'esolution.}.
\end{defi}

Le cadre des paires est maintenant le cadre naturel dans lequel 
prend place le programme de Mori (on sait m\^eme maintenant
qu'il s'av\`ere n\'ecessaire de travailler avec des 
$\mathbb R$-diviseurs). Dire qu'une paire $(X,\Delta)$ est klt
signifie 
que $X$ et $\Delta$ sont peu singuli\`eres. Ces singularit\'es 
sont in\'evitables lorsqu'on applique le programme de Mori. 
Il est imm\'ediat (voir la preuve du lemme \ref{perturb} ci-dessous)
et important de constater qu'une petite perturbation d'une 
paire klt est encore une paire klt : si $D$ est un $\mathbb Q$-diviseur de 
Cartier
effectif et si $(X,\Delta)$ est klt, alors $(X,\Delta + \varepsilon D)$
est klt pour tout $\varepsilon >0$ suffisamment petit.

\begin{ex}
\begin{enumerate}

\item Les paires klt sont lc.

\item Soient $S$ une surface lisse et $C\subset S$ une 
courbe lisse de genre $g$ telle que $C^2=-n<0$.
Soit $\pi : S \to S_0$ le morphisme birationnel
consistant \`a contracter $C$ sur un point.
Alors $\pi$ est une r\'esolution de la paire $(S_0,0)$
et  
$$ K_S = \pi^* K_{S_0}+\left(-1+\frac{2-2g}{n}\right)C.$$
La paire $(S_0,0)$ est donc klt si et seulement si $g=0$,
est lc si et seulement si $g\leq 1$. 

\item Si $X$ est lisse et $\Delta=\sum_{i=1}^N d_i \Delta_i$
est un diviseur \`a croisements normaux simples, alors  
$${\rm Nklt}(X,\Delta)= \cup_{ \{ i \mid d_i \geq 1 \} } \Delta_i$$
et $(X,\Delta)$
est klt si et seulement si $d_i<1$ pour tout $i$, est lc 
si et seulement si $d_i\leq 1$ pour tout $i$.

\end{enumerate}
\end{ex}

\medskip

Les r\'esultats suivants sont dus \`a Hacon et ${\rm M^c}$Kernan
(on renvoie \`a l'article original \cite{HM07} pour des \'enonc\'es un peu 
plus g\'en\'eraux).

\begin{theo}{\bf (Hacon et ${\rm \bf M^c}$Kernan)}\label{HMmain}
Soient $(X,\Delta)$ une paire klt et $f~: X \to S$ un morphisme
propre tels que $-K_X$ est $f$-gros et $-(K_X+\Delta)$
est $f$-nef\footnote{Un $\mathbb Q$-diviseur est dit $f$-gros
s'il est $\mathbb Q$-lin\'eairement \'equivalent \`a $A+B$ o\`u $A$ est
$f$-ample et
$B$ est un $\mathbb Q$-diviseur de Weil
effectif. Un diviseur $\mathbb Q$-Cartier est
$f$-nef s'il est nef en restriction \`a {\bf toute} 
fibre de $f$.}.
Soient $Y$ une vari\'et\'e alg\'ebrique normale,
$g~: Y \to X$ un morphisme propre birationnel et
$\pi = f \circ g : Y \to S$. Alors toute composante connexe de
toute fibre de
$\pi$ est rationnellement connexe par cha\^{\i}nes.
\end{theo}

{\bf Mise en garde.} Dans cet \'enonc\'e, les composantes
connexes des fibres de $\pi$
ne sont en g\'en\'eral pas irr\'eductibles. Dire qu'une composante
connexe $W$ d'une fibre 
de $\pi$ est rationnellement
connexe par cha\^{\i}nes implique que deux points quelconques
de $W$ peuvent \^etre joints par une cha\^{\i}ne de courbes
rationnelles contenues dans $W$. Ceci n'implique pas
que chaque composante irr\'eductible de $W$ est elle-m\^eme rationnellement
connexe par cha\^{\i}nes. Soient $X$ une vari\'et\'e lisse de dimension
$3$, $B_x(X)$ l'\'eclat\'e de $X$ en $x$ et $Y$ l'\'eclat\'e
de $B_x(X)$ le long d'une courbe elliptique contenue
dans le diviseur exceptionnel de $B_x(X) \to X$. Soit $\pi : Y \to X$
la compos\'ee. Alors $\pi^{-1}(x)$ est connexe,
a deux composantes irr\'eductibles,
est rationnellement connexe par cha\^{\i}nes mais la composante correspondant
au diviseur exceptionnel du deuxi\`eme \'eclatement n'est pas rationnellement
connexe par cha\^{\i}nes : c'est un fibr\'e en $\mathbb P ^1$ sur une courbe
elliptique.

\medskip

Les deux \'enonc\'es suivants \'etaient conjectur\'es par Shokurov.
Il s'agit du th\'eor\`eme \ref{HMmain} dans les cas extr\^emes o\`u
$S$ est un point et $S=X$.

Dans le cas o\`u $S$ est r\'eduit \`a un point, on obtient 
une g\'en\'eralisation du th\'eor\`eme \ref{fanorc} due \`a Zhang.

\begin{cor}\label{zhang}{\bf (Hacon et ${\rm \bf M^c}$Kernan, Zhang)}
Soit $(X,\Delta)$ une paire klt, avec $X$ projective, 
telle que $-(K_X+\Delta)$
est gros et nef. Alors $X$ est rationnellement 
connexe.
\end{cor}

Dans le cas o\`u $S=X$, on obtient encore le corollaire suivant, 
g\'en\'eralisation
des propositions \ref{prol} et \ref{term} et du corollaire 
\ref{prolon}~; le point (1) est la conjecture de connexit\'e rationnelle
de Shokurov.

\begin{cor}\label{cozhang}{\bf (Hacon et ${\rm \bf M^c}$Kernan)} 
Soit $(X,\Delta)$ une paire klt.
\begin{enumerate} 
\item Soient $Y$ une vari\'et\'e alg\'ebrique normale et $\pi~: Y \to X$ 
un morphisme birationnel propre.
Alors toute fibre\footnote{Inutile ici de se restreindre aux composantes
connexes, le th\'eor\`eme principal de Zariski affirme que si
$\pi : Y \to X$ est propre birationnel avec $X$ normale, alors les
fibres de $\pi$ sont connexes.} de
$\pi$ est rationnellement connexe par cha\^{\i}nes.
\item Soient $Z$ une vari\'et\'e alg\'ebrique normale
et $\varphi : X \dashrightarrow Z$ 
une application rationnelle propre. Si $Z$ ne contient pas de courbes 
rationnelles, alors $\varphi$ se prolonge en une application
r\'eguli\`ere $\varphi : X \to Z$.
\end{enumerate}
\end{cor}

Une derni\`ere cons\'equence est la connexit\'e rationnelle par cha\^{\i}nes
des fibres des contractions de Mori fournies par le th\'eor\`eme du 
c\^one.

\begin{cor}
Soient $X$ une vari\'et\'e projective \`a singularit\'es 
terminales
et $f: X \to Z$ un morphisme sur une vari\'et\'e projective $Z$.
Soient $R_i:= {\mathbb R}^+ [\Gamma_i ]$ une ar\^ete du c\^one 
$\overline{{\rm NE}}(X/Z)$ telle que $-K_X \cdot \Gamma _i >0$
et $c_i : X/Z \to X_i /Z$ la contraction associ\'ee. 
Alors toute fibre de $c_i$ est rationnellement connexe par cha\^{\i}nes.
\end{cor}

Mentionnons aussi que Broustet et Pacienza ont \'etendu
depuis les r\'esultats de Hacon et ${\rm M^c}$Kernan au cas o\`u 
$-(K_X+\Delta)$
n'est plus suppos\'e $f$-nef \cite{BP09}.

\subsection{Les grandes lignes de la preuve du th\'eor\`eme \ref{HMmain}}
{\em On donne dans ce paragraphe les grandes \'etapes, 
tous les d\'etails sont donn\'es
dans la suite de ce cours}.

Soient $(X,\Delta)$ une paire klt et $f~: X \to S$ un morphisme
propre tels que $-K_X$ est $f$-gros et $-(K_X+\Delta)$
est $f$-nef. Soient $Y$ une vari\'et\'e alg\'ebrique normale,
$g~: Y \to X$ un morphisme propre birationnel et
$\pi = f \circ g : Y \to S$.
Supposons pour simplifier que $K_X+\Delta$
est lin\'eairement \'equivalent \`a $0$ et que $\Delta = A +B$
avec $A$ ample aussi g\'en\'eral que souhait\'e, $B$ effectif
et $(X,B)$ klt. 

\subsubsection{Le cas o\`u $S$ est r\'eduit \`a un point.}
Il s'agit alors de montrer que $Y$ 
est rationnellement connexe par cha\^{\i}nes. Supposons \`a nouveau
pour simplifier
que $Y$ est lisse et qu'il y a un morphisme r\'egulier
$t : Y \to R(Y)$ de $Y$ sur son quotient rationnel $R(Y)$.
Comme $R(Y)$ n'est pas 
unir\'egl\'e (th\'eor\`eme \ref{quotientuni}),
le th\'eor\`eme
\ref{carunireg} affirme que $K_{R(Y)}$ est pseudo-effectif,
et par suite que $$\kappa(R(Y),K_{R(Y)}+H)=\dim R(Y)$$
pour tout $\mathbb Q$-diviseur ample $H$ sur $R(Y)$.

On peut alors \'ecrire\footnote{Les participants
de l'Ecole d'Et\'e 2007 de l'Institut Fourier \`a Grenoble
reconnaitront la d\'ecomposition
``dagger'' dont Alessio Corti est fan.} 
$$ K_{Y}+ \Omega = g^*(K_{X}+\Delta)+ D$$
o\`u $\Omega$ et $D$ sont effectifs sans composante commune,
la paire $(Y,\Omega)$
est klt et $D$ est $g$-exceptionnel.
On peut aussi supposer que 
$\Omega = g^*A+\bar B$ avec $(Y,\bar B)$ klt.
Comme $D$ est $g$-exceptionnel,
$$ \kappa (Y, K_{Y}+ \Omega ) = 
\kappa (X,K_X+ \Delta) =0.$$

En utilisant la positivit\'e de 
$\Omega$, on montre alors, c'est l'\'etape cl\'e, qu'il
existe un $\mathbb Q$-diviseur ample $H$ sur $R(Y)$
tel que 
$$0= \kappa (Y, K_{Y}+ \Omega ) \geq \kappa(R(Y),K_{R(Y)}+H)=\dim R(Y)$$
(et ainsi $R(Y)$ est r\'eduit \`a un point, {\em i.e.}
$Y$ est rationnellement connexe par cha\^{\i}nes). Cette in\'egalit\'e
(de type conjecture $C_{n,m}$) est obtenue \`a l'aide d'un th\'eor\`eme
de positivit\'e d'images directes (d\^u \`a Campana) pour les faisceaux
du type $t_*\mathcal O _Y(m(K_{Y/R(Y)}+C))$, o\`u
$C$ est un diviseur effectif dont la restriction \`a la fibre g\'en\'erale de
$t$ est lc.

\subsubsection{Le cas o\`u $\dim S >0$.}
Fixons $s\in S$ et supposons que $S$ est affine. Supposons aussi
pour simplifier que $Y$ est lisse et que $F := \pi^{-1}(s)$
est un diviseur \`a croisements normaux simples.

On montre alors que si $k$ est le nombre de composantes irr\'eductibles 
de $F=\pi^{-1}(s)$, 
il existe une num\'erotation des composantes 
irr\'eductibles $F_1,\ldots,F_k$ de $F$ et pour tout $i$
une paire $(Y,\Theta _{i})$ tels que
\begin{enumerate}
\item $\Theta _{i} = \bar A_0 + B_i$ o\`u $\bar A_0$ est
un $\mathbb Q$-diviseur ample et effectif
et $B_i$ est un $\mathbb Q$-diviseur effectif,
\item $K_Y +\Theta _{i} \sim \bar E_i $ o\`u $\bar E_i$ est effectif,
$g$-exceptionnel et n'a pas de composante commune avec $\Theta_i$,
\item ${\rm Nklt}(Y,\Theta _{i})=F_1\cup \cdots \cup F_i$ et
le coefficient de $F_i$ dans $\Theta _{i}$ vaut $1$,
celui des $F_1,\ldots,F_{i-1}$ est $>1$. En particulier, 
$\Theta_{1}=\bar A_0+F_1+C_1$ o\`u
$(Y,C_{1})$ est klt et $C_{1}$ ne contient pas $F_1$ dans son 
support.
\end{enumerate} 

On montre ensuite que 
\begin{enumerate}
\item La composante $F_1$ est rationnellement connexe,
\item pour tout $i\geq 2$, $F_i$ est rationnellement connexe
par cha\^{\i}nes modulo $$F_i\cap {\rm Nklt}(Y,\Theta_{i-1})
= F_i\cap(F_1\cup \ldots \cup F_{i-1}),$$ ce qui signifie 
que pour tout $x \in F_i$, il existe une cha\^{\i}ne de courbes
rationnelles joignant $x$ \`a un point de 
$F_i\cap (F_1\cup \ldots \cup F_{i-1})$.
\end{enumerate}

Pour montrer que $F_1$ est rationnellement connexe,
on consid\`ere comme pr\'ec\'edemment son quotient
rationnel $t : F_1 \to R(F_1)$ (que l'on suppose r\'egulier 
pour simplifier) puis l'on montre qu'il existe
un $\mathbb Q$-diviseur ample $H$ sur $R(F_1)$
tel que
$$  \kappa(F_1,K_{F_1}+ \{ \Theta_1\}_{|F_1})
\geq \kappa(R(F_1),K_{R(F_1)}+H)$$
\`a l'aide \`a nouveau du th\'eor\`eme de positivit\'e d'images directes
de Campana.

Pour conclure, il suffit de montrer que  
$\kappa(F_1,K_{F_1}+ \{ \Theta_1\}_{|F_1})=0$. 
Or $K_{F_1}+ \{ \Theta_1\}_{|F_1}= (K_Y+\Theta_1)_{|F_1}$
si bien qu'il s'agit alors d'appliquer un th\'eor\`eme d'extension
de sections fourni lui-aussi par Hacon et ${\rm M^c}$Kernan.
 
\subsection{Un peu de technologie des paires}

L'avantage de travailler avec des paires r\'eside dans leur grande
souplesse, on peut les perturber  et ``suivre leur positivit\'e'' 
dans les log-r\'esolutions
successives. On donne ici un lemme qui illustre bien ce
principe.

\begin{lemm}\label{perturb} Soient $(X,\Delta)$ une paire klt. 
On suppose
que $\Delta$ est gros. Alors il existe un $\mathbb Q$-diviseur
ample effectif $A$ et un $\mathbb Q$-diviseur effectif $\Gamma_1$ tels que
les paires $(X,\Gamma_1+A)$ et $(X,\Gamma_1)$ sont klt et
$\Delta \sim \Gamma_1+A$\footnote{Dans tout ce texte, $\sim$
est l'\'equivalence lin\'eaire des $\mathbb Q$-diviseurs~: deux 
$\mathbb Q$-diviseurs $M$ et $N$ sont lin\'eairement \'equivalents 
s'il existe un entier positif $m$ tels que $mM$ et $mN$
sont deux diviseurs entiers lin\'eairement \'equivalents.}.
\end{lemm}   

Ce lemme permet donc de remplacer la paire klt $(X,\Delta)$
par une paire $(X,\Gamma)$ elle aussi klt, avec $\Gamma \sim \Delta$
et $\Gamma = {\rm ample} + {\rm effectif}$. L'hypoth\`ese que $\Delta$
est gros est \'evidemment essentielle ici.

\medskip

{\bf D\'emonstration.}
Comme $\Delta$ est gros, il existe des $\mathbb Q$-diviseurs
$H$ et $B$ tels que $H$ est ample, $B$ est effectif et
$\Delta \sim H+B$. Soit $H_m \in |mH|$ un \'el\'ement g\'en\'eral
et $A_m := (1/m)H_m$, o\`u $m$ est choisi suffisam\-ment grand 
et suffisamment divisible pour que $mH$ soit un diviseur (entier) 
tr\`es ample sur $X$. 
On pose alors 
$$ \Gamma := (1-\varepsilon)\Delta + \varepsilon A_m + \varepsilon B
= \Gamma_1 + A \sim \Delta$$
o\`u $A:= \varepsilon A_m$ et $\Gamma_1 := 
(1-\varepsilon)\Delta + \varepsilon B$.
Par construction, $K_X+A_m+B \sim K_X+\Delta$ est $\mathbb Q$-Cartier,
donc $K_X+\Gamma \sim (1-\varepsilon)(K_X+\Delta)+\varepsilon(K_X+A_m+B)$
l'est aussi, de m\^eme que $K_X+\Gamma_1$.

Soit $\pi~: Y \to X$ une log-r\'esolution $\Delta$ telle
que ${\rm Exc}(\pi)+\Delta'+B'$ soit un diviseur \`a croisements
normaux simples. 
Si $H_m$ est suffisamment g\'en\'eral, la transform\'ee stricte
de $H_m$ est \'egale \`a sa transform\'ee totale $\pi^* H_m$.  

Le lemme d\'ecoule alors des calculs suivants, la somme
porte \`a chaque fois sur les diviseurs $\pi$-exceptionnels
et les $'$ d\'esignent les transform\'ees strictes.
\begin{align*} 
K_Y+A_m '+B' & = \pi^*(K_X+A_m+B) + \sum_E r_E E\\
K_Y + \Delta' & = \pi^*(K_X+\Delta)+ \sum_E a(E,\Delta)E\\
K_Y + \Gamma_1' &=
\pi^*(K_X +\Gamma_1)
+ \sum_E((1-\varepsilon)a(E,\Delta)+ \varepsilon r_E)E\\
K_Y+ \Gamma' &=
\pi^*(K_X +\Gamma)
+ \sum_E((1-\varepsilon)a(E,\Delta)+ \varepsilon r_E)E.
\end{align*}
Comme $a(E,\Delta)>-1$ pour tout $E$, les quantit\'es
$(1-\varepsilon)a(E,\Delta)+ \varepsilon r_E$
sont aussi $>-1$ pour $\varepsilon >0$ suffisamment petit.
De m\^eme, les multiplicit\'es de $\Gamma$
et $\Gamma_1$ le long de leurs composantes irr\'eductibles
sont $<1$ si $\varepsilon >0$ est suffisamment petit.\finpreuve

\subsection{Le th\'eor\`eme de positivit\'e d'images directes de Campana.}

L'un des ingr\'edients essentiels, non encore introduit dans
ce cours, de la preuve du th\'eor\`eme
de Hacon et ${\rm \bf M^c}$Kernan est un th\'eor\`eme de positivit\'e 
d'images directes. 
Pour l'\'enoncer, nous avons besoin d'une d\'efinition due \`a Viehweg. 

\begin{defi}
Un faisceau coh\'erent sans torsion
$\mathcal E$ sur une vari\'et\'e projective $V$ est faiblement positif
si, $V_0$ d\'esignant le plus grand ouvert sur lequel $\mathcal E$ 
est localement libre, il existe un ouvert non vide $U \subset V_0$
tel que pour tout diviseur ample $H$ sur $V$, pour tout entier $a>0$,
il existe $b>0$ tel que les sections globales de 
$\mathcal F^{ab}:=  ({\rm Sym}^{ab}\mathcal E _{|V_0})\otimes 
H_{|V_0}^{\otimes b}$ engendrent $\mathcal F^{ab}_{|U}$. 
\end{defi}

{\bf Remarque.} 
On met en garde le lecteur avec le fait que le faisceau nul 
est faiblement positif~!!

\medskip

Le th\'eor\`eme de positivit\'e d'images directes
d\^u \`a Campana \cite{campab}, et faisant suite aux travaux de Viehweg,
s'\'enonce de la fa\c{c}on suivante.

\begin{theo}{\bf (Campana)}\label{Campa}
Soit $f~: V' \to V$ un morphisme \`a fibres connexes
entre vari\'et\'es projectives lisses et soit $C$ un $\mathbb Q$-diviseur
effectif. On suppose que 
la restriction de $C$
\`a une fibre g\'en\'erale de $f$
est lc.
Alors le faisceau $f_*\mathcal O _{V'}(m(K_{V'/V}+C))$
est faiblement positif pour tout entier $m$ tel que  
$m(K_{V'/V}+C)$ est entier.
\end{theo}

On l'a dit, le faisceau nul est faiblement positif. Sous des hypoth\`eses
de positivit\'e de la fibre, on garantit que le faisceau images directes
est non trivial.

\begin{theo}\label{Debar}  
Soit $f~: V' \to V$ un morphisme \`a fibres connexes
entre vari\'et\'es projectives lisses, soit $C$ un $\mathbb Q$-diviseur
effectif et soit $W$ la fibre g\'en\'erale de $f$.

On suppose que 
la paire $(W,C_{|W})$
est lc et que
$$ \kappa(W,m_0(K_{V'/V}+C)_{|W}) \geq 0$$
pour un entier $m_0$ tel que $m_0(K_{V'/V}+C)$ est entier.
 
Alors, pour tout diviseur ample $H$ sur $V$, il y a un entier 
$b>0$ tel que 
$$ H^0(V',bm_0(K_{V'/V}+C+f^*H))\neq 0.$$
\end{theo}

{\bf D\'emonstration.} R\'ep\'etons, pour le confort 
du lecteur,
cet argument d\^u \`a Campana et Viehweg et d\'etaill\'e
dans \cite{Deb06}. 

Le lemme d'aplatissement de Raynaud combin\'e au th\'eor\`eme
de d\'esingularisation d'Hiro\-naka permet de construire un diagramme
 
\centerline{
\xymatrix{V'_1 \ar[r]^ {\tau '} \ar[d]_{f_1} & {V'}\ar[d]^{f}
\\
V_1 \ar[r]^{\tau} & V \\
}
}
\noindent o\`u $V'_1$ et $V_1$ sont des vari\'et\'es projectives lisses,
$\tau '$ et $\tau$ sont des modifications 
telles que d'une part $f_1({\rm Exc}(\tau '))\subset {\rm Exc}(\tau)$
et d'autre part
{\em toute hypersurface $S\subset V'_1$
dont l'image par $f_1$ est de codimension $\geq 2$ dans
$V$ est $\tau'$-exceptionnelle.}  
 
Comme $(\tau ') ^* C$ est lc,
le th\'eor\`eme \ref{Campa}
assure que $$(f_1)_* {\mathcal O}_{V'_1}(m(K_{V'_1/V_1}+ 
(\tau') ^* C+f_1^*K_{V_1/V}))$$ 
est faiblement positif pour tout
$m$ suffisamment divisible.
Ce m\^eme faisceau est aussi non nul car
$$\tau_* (f_1)_* {\mathcal O}_{V'_1}(m_0(K_{V'_1 /V_1}+ 
(\tau') ^* C+f_1^*K_{V_1/V}))= 
f _*{\mathcal O}_{V'}
(m_0(K_{V'/V}+ C))\neq 0$$
puisque la fibre du faisceau $f _*{\mathcal O}_{V'}
(m_0(K_{V'/V}+ C))$ en un point g\'en\'eral de $V$ vaut
$$H^0(W,(m_0(K_{V'/V}+ C))_{|W}),$$ suppos\'e non nul. 

Comme $H$ est ample, $\tau ^* (m_0 H)$ est lin\'eairement
\'equivalent \`a $H_{V_1}+ E_{V_1}$ avec $H_{V_1}$ 
ample sur $V_1$ et $E_{V_1}$ effectif.
Par d\'efinition de la faible positivit\'e du faisceau non nul 
$f _*{\mathcal O}_{V'}
(m_0(K_{V'/V}+ C))$, il existe
un entier $b$ et un ouvert $V_0$ de $V_1$, dont le compl\'ementaire
dans $V_1$ est de codimension $\geq 2$, tels que 
$$ H^0(V_0,(f_1)_* {\mathcal O}_{V'_1}(bm_0(K_{V'_1/V_1}+ 
(\tau') ^* C+f_1^*K_{V_1/V}))\otimes H_{V_1} 
^{\otimes b}) \neq 0.$$
{\em A fortiori}, puisque $\tau^* (m_0 H) \sim H_{V_1}+E_{V_1}$,
\begin{align*}
H^0(f_1^{-1}(V_0),bm_0 (K_{V'_1 /V_1}+ 
(\tau') ^* C + f_1^* K_{V_1/V}+f_1^* \tau ^*H))&=\\
H^0(V_0,(f_1)_* {\mathcal O}_{V'_1}(bm_0(K_{V'_1 /V_1}+ 
(\tau') ^* C+ f_1^* K_{V_1/V}))\otimes \tau ^* H  
^{\otimes b}) \neq 0.
\end{align*}
Or, 
$$ K_{V'_1 /V_1}+ (\tau') ^* C + f_1^* K_{V_1/V}+f_1^* \tau ^*H
= K_{V'_1 /V'}+ 
(\tau') ^* ( K_{V'/V} + C + f ^*H).$$
Comme $K_{V'_1 /V'}$ est $\tau '$-exceptionnel, de m\^eme que le 
compl\'ementaire dans $V'_1$ de $f_1^{-1}(V_0)$, on en d\'eduit que 
$$0 \neq H^0(f_1^{-1}(V_0),bm_0(K_{V'_1 /V_1}+ 
(\tau') ^* C + f_1^* K_{V_1/V}+f_1^* \tau ^*H))=
H^0(V',bm_0(K_{V' /V} + C + f ^*H)).$$
\finpreuve

On en d\'eduit, suivant la terminologie
introduite par Campana, une version ``orbifold''
de la conjecture $C_{n,m}$.

\begin{cor}\label{cnm}
Soit $f~: V' \to V$ un morphisme \`a fibres connexes
entre vari\'et\'es projectives lisses, soit $C$ un $\mathbb Q$-diviseur
effectif et soit $W$ la fibre g\'en\'erale de $f$.
On suppose que 
la paire $(W,C_{|W})$
est lc et que
$$ \kappa(W,(K_{V'/V}+C)_{|W}) \geq 0.$$
 
Alors, pour tout diviseur ample $H$ sur $V$,
$$ \kappa(V', K_{V'}+C+ 2f^*H) \geq \kappa(V,K_V+H).$$
\end{cor}

{\bf D\'emonstration.} 
Comme il existe un entier $b>0$ tel que
$$ H^0(V',b(K_{V'/V}+ C + f^* H))\neq 0,$$
on en d\'eduit que pour tout $m$ suffisamment divisible, 
$H^0(V,mb(K_{V}+H))$ s'injecte 
dans 
$$ H^0(V',mb(K_{V'/V}+C+ 
f^* H + f^*
(K_{V}+H))) \simeq 
H^0(V',mb(K_{V'}+C + 2f^* H)).$$
{\em A fortiori}, 
$$\dim H^0(V,mb(K_{V}+H)) \leq \dim H^0(V',mb(K_{V'}+C + 2f^* H)).$$
\finpreuve

\subsection{Le r\'esultat de Zhang.} 
On traite ici le cas o\`u $\dim S =0$. Rappelons
l'\'enonc\'e.

\medskip

{\bf Corollaire \ref{zhang}.} {\bf (Zhang)}
{\em
Soit $(X,\Delta)$ une paire klt, avec $X$ projective, 
telle que $-(K_X+\Delta)$
est gros et nef. Alors $X$ est rationnellement 
connexe.
}

\medskip

{\bf D\'emonstration.}

Soit $g~:Y \to X$ une log-r\'esolution de la paire
$(X,\Delta)$ fix\'ee dans toute la suite. On va montrer
que $Y$ est rationnellement connexe. Comme $Y$ est lisse,
il suffit de montrer que $Y$ est rationnellement connexe
par cha\^{\i}nes (corollaire \ref{coversus}), 
autrement dit que son quotient rationnel $R(Y)$ est
r\'eduit \`a un point. 

\medskip

Comme $-(K_X+\Delta)$
est gros et nef, le ``base-point-free theorem''\footnote{
En voici l'\'enonc\'e, dans le cadre relatif utile pour la suite
du cours. On l'applique ici dans le cas o\`u $Z$
est r\'eduit \`a un point. 

\begin{theo} {\bf (``Base-point-free theorem'')}
Soient $(X,\Gamma)$ une paire klt, $\pi : X \to Z$ un
morphisme projectif et $D$ un diviseur
de Cartier $\pi$-nef sur $X$. On suppose qu'il existe
un rationnel $a>0$ tel que $aD-(K_X+\Gamma)$ est $\pi$-nef
et $\pi$-gros. Alors le syst\`eme lin\'eaire $mD$
est $\pi$-globalement engendr\'e
pour tout $m$ suffisamment divisible.
\end{theo}

Dans l'\'enonc\'e de ce th\'eor\`eme, ``$mD$ est $\pi$-globalement engendr\'e''
signifie qu'il existe un recouvrement de
$Z$ par des ouverts affines tel que $mD$ est sans point base
en restriction \`a $\pi^{-1}(U)$ pour tout ouvert $U$ du recouvrement.}
affirme que $-m(K_X+\Delta)$ est sans point base
pour $m$ assez grand et assez divisible.  
Soit alors $D\in |-m(K_X+\Delta)|$ un \'el\'ement g\'en\'eral pour $m$ 
suffisamment grand et divisible. 
La paire
$(X,\Delta+(1/m)D)$ est toujours klt et si $\Delta_1:=\Delta+(1/m)D$,
alors $\Delta_1$ est gros et $K_X + \Delta_1 \sim 0$.
Comme $\Delta_1$ est gros, le lemme \ref{perturb}
montre qu'il existe $\Delta_2 = \varepsilon A_1 + B_1$,
$A_1$ ample et $B_1$ effectif, $A_1$ suffisamment g\'en\'eral
pour que sa transform\'ee totale sous $g$ soit \'egale \`a sa transform\'ee
stricte sous $g$, tels que $(X,B_1)$ est klt, $(X,\Delta_2)$ est klt 
et
$K_X + \Delta_2 \sim 0$.

\medskip

{\em Dor\'enavant, on peut donc supposer que 
$(X,\Delta)$ est une paire klt, que $K_X+\Delta$
est lin\'eai\-rement \'equivalent \`a $0$ et que $\Delta = A +B$
avec $A$ ample aussi g\'en\'eral que souhait\'e, $B$ effectif
et $(X,B)$ klt.}

\medskip

{\em Etape 2.}
Soient $t~: Y \dashrightarrow R(Y)$ le quotient rationnel
de $Y$, $\varphi : \tilde Y \to Y$ une suite d'\'eclatements
de centres lisses qui l\`eve les ind\'eterminations de
$t$ et $\tilde t = t \circ \varphi : \tilde Y \to R(Y)$ le
morphisme associ\'e.

On \'ecrit \`a nouveau la d\'ecomposition 
$$ K_{\tilde Y}+ \Omega = (g \circ \varphi)^*(K_{X}+\Delta)+ D$$
o\`u $\Omega$ et $D$ sont effectifs sans composante commune,
la paire $(\tilde Y,\Omega)$
est klt, $(g \circ \varphi)_*(\Omega)=\Delta$ 
et $D$ est $(g\circ \varphi)$-exceptionnel.
Enfin, $\Omega = (g\circ \varphi)^*A+\bar B$ avec $(\tilde Y,\bar B)$ klt.
Comme $D$ est $(g\circ \varphi)$-exceptionnel,
$$ \kappa (\tilde Y, K_{\tilde Y}+ \Omega ) = 
\kappa (X,K_X+ \Delta) =0.$$
  
\medskip

{\em Etape 3.} D\'emontrons le lemme suivant.

\begin{lemm} \label{kodaira}
Il existe un $\mathbb Q$-diviseur $H$ ample sur $R(Y)$
tel que $(g\circ \varphi)^*A \sim 2\tilde t ^* H +\Gamma$
o\`u $\Gamma$ est effectif et la paire $(\tilde Y,\bar B+\Gamma)$
est klt.
\end{lemm}
 
\medskip

En effet, il existe un diviseur effectif $E^*$ tel que
$(g\circ \varphi)^*A - \varepsilon E^*$ est ample sur $\tilde Y$ pour tout
$\varepsilon >0$ suffisamment petit.
Si $H_1$ est tr\`es ample
sur $R(Y)$, soit $A_m$ un membre g\'en\'eral de 
$|m((g\circ \varphi)^*A - \varepsilon E^*)-\tilde t^*H_1|$.
On a alors 
$$ (g\circ \varphi)^*A \sim \frac{1}{m}A_m +\varepsilon E^*+
\frac{1}{m}\tilde t^*H_1.$$
Si $m$ est suffisamment grand et $\varepsilon$ suffisamment petit,
la paire $(\tilde Y,(1/m)A_m +\varepsilon E^*+\bar B)$ est klt puisque
$(\tilde Y,\bar B)$ l'est 
et $H=(1/2m)H_1$ convient. Ceci termine ce lemme.\finpreuve

\medskip

{\em Etape 4.} 
La paire $\bar B+\Gamma$ est klt, donc sa restriction \`a la fibre 
g\'en\'erale $\tilde t ^{-1}(u)$ l'est aussi. 
De plus, 
$$K_{\tilde Y/R(Y)}+\bar B+\Gamma \sim K_{\tilde Y}-\tilde t^*K_{R(Y)}
+ \Omega - (g\circ \varphi)^*A+(g\circ \varphi)^*A
-2\tilde t ^* H$$
donc
$$(K_{\tilde Y/R(Y)}+\bar B+\Gamma)_{|\tilde t ^{-1}(u)}
\sim (K_{\tilde Y}+ \Omega)_{|\tilde t ^{-1}(u)} \sim D_{|\tilde t ^{-1}(u)}$$
est effectif
donc $$ H^0(\tilde t ^{-1}(u),m(K_{\tilde Y/R(Y)}+\bar B+
\Gamma)_{|\tilde t ^{-1}(u)})\neq 0.$$

Comme $(g\circ \varphi)^*A \sim 2\tilde t^* H + \Gamma$
et $\Omega = (g\circ \varphi)^*A+\bar B$, on en d\'eduit par le Corollaire
\ref{cnm}
que 
\begin{align*}
0 = \kappa(\tilde Y,K_{\tilde Y}+\Omega)&  
= \kappa(\tilde Y,K_{\tilde Y}+ \bar B + \Gamma + 2\tilde t ^* H)\\
& \geq \kappa(R(Y),K_{R(Y)}+H).
\end{align*}

Or $R(Y)$ n'est pas 
unir\'egl\'e (th\'eor\`eme \ref{quotientuni}),
le th\'eor\`eme
\ref{carunireg} affirme donc que $K_{R(Y)}$ est pseudo-effectif,
et par suite que $$\kappa(R(Y),K_{R(Y)}+H)=\dim R(Y).$$ 
Contradiction sauf si $R(Y)$ est un point !\finpreuve

\subsection{Le th\'eor\`eme de Hacon et ${\rm \bf M^c}$Kernan.} 
On traite ici le cas o\`u $\dim S > 0$. Rappelons
l'\'enonc\'e.

\medskip

{\bf Th\'eor\`eme \ref{HMmain}.} {\bf (Hacon et ${\rm \bf M^c}$Kernan)}
{\em
Soient $(X,\Delta)$ une paire klt et $f~: X \to S$ un morphisme
propre tels que $-K_X$ est $f$-gros et $-(K_X+\Delta)$
est $f$-nef.
Soient $Y$ une vari\'et\'e alg\'ebrique normale,
$g~: Y \to X$ un morphisme propre birationnel et
$\pi = f \circ g : Y \to S$. Alors toute composante connexe de
toute fibre de
$\pi$ est rationnellement connexe par cha\^{\i}nes.
}

\subsubsection{Quelques r\'eductions ``faciles''}
Comme l'\'enonc\'e est local, on peut supposer que $S$ est affine,
et \`a l'aide 
du ``base-point-free theorem'' dans sa version relative, on montre,
exactement comme dans la preuve du th\'eor\`eme de Zhang\footnote{
Dans notre situation, $S$ est affine, et quitte \`a r\'eduire $S$,
on peut 
\'ecrire $-K_X \sim A+B$ avec $A$ ample et $B$ effectif,
et $\Gamma = (1-\varepsilon)\Delta + \varepsilon B$.
La paire $(X,\Gamma)$ est klt pour $\varepsilon >0$ petit
et $\varepsilon(-(K_X+\Delta))-(K_X+\Gamma) \sim -(K_X+\Delta) + 
\varepsilon A $ est nef et gros, donc $|-m(K_X+\Delta)|$ est
sans point base pour $m$ suffisamment divisible.
}, 
qu'il existe $\Delta_2 = \varepsilon A_1 + B_1$,
$A_1$ ample et $B_1$ effectif, $A_1$ suffisamment g\'en\'eral
pour que sa transform\'ee totale sous $g$ soit \'egale \`a sa transform\'ee
stricte sous $g$, tels que $(X,B_1)$ est klt, $(X,\Delta_2)$ est klt
et
$K_X + \Delta_2 \sim 0$.

\medskip

{\em Dor\'enavant, on suppose donc que 
$(X,\Delta)$ est une paire klt, que $f~: X \to S$ est un morphisme
tel que $K_X+\Delta$
est lin\'eairement \'equivalent \`a $0$ et que $\Delta = A +B$
avec $A$ ample aussi g\'en\'eral que souhait\'e, $B$ effectif
et $(X,B)$ klt. On supposera aussi, quitte \`a faire une factorisation de
Stein, que $f$ est \`a fibres connexes et que $S$ est affine normale.}

\subsubsection{On r\'ealise la fibre comme lieu non klt d'une paire
bien choisie}

Fixons $s \in S$. Quitte \`a \'eclater encore $Y$, on peut supposer
que $g$ est une log-r\'esolution
de $(X,\Delta)$ telle
que $$\pi^{-1}(s) + {\rm Exc}(g) + \Delta'$$  
est un diviseur \`a croisements normaux simples\footnote{
$\pi^{-1}(s)$ et ${\rm Exc}(g)$ ont bien s\^ur des composantes communes.}
~: si $\sigma : Y'\to Y$
est une suite d'\'eclatements (donc \`a fibres connexes) 
et si toute composante connexe
de $(\pi \circ \sigma)^{-1}(s)$ est rationnellement connexe par 
cha\^{\i}nes, alors toute composante connexe de $\pi^{-1}(s)$ est 
aussi rationnellement connexe par 
cha\^{\i}nes. 

\medskip

Soit $F:=\pi^{-1}(s)$ dont on rappelle qu'on vient de supposer que
c'est un diviseur \`a croisements normaux simples.

\begin{prop} 
Si $k$ est le nombre de composantes irr\'eductibles 
de $F=\pi^{-1}(s)$, 
il existe une num\'erotation des composantes 
irr\'eductibles $F_1,\ldots,F_k$ de $F$ et pour tout $i$
compris entre $1$ et $k$
une paire $(Y,\Theta _{i})$ tels que
\begin{enumerate}
\item $\Theta _{i} = \bar A_0 + B_i$ o\`u $\bar A_0$ est
un $\mathbb Q$-diviseur ample et effectif
et $B_i$ est un $\mathbb Q$-diviseur effectif,
\item $K_Y +\Theta _{i} \sim \bar E_i $ o\`u $\bar E_i$ est effectif,
$g$-exceptionnel et n'a pas de composante commune avec $\Theta_i$,
\item 
${\rm Nklt}(Y,\Theta _{i})=F_1\cup \cdots \cup F_i$ et
le coefficient de $F_i$ dans $\Theta _{i}$ vaut $1$,
celui des $F_1,\ldots,F_{i-1}$ est $>1$, celui des $F_{i+1},\ldots,F_k$
est $<1$. En particulier, $\Theta_{1}=\bar A_0+F_1+C_1$ o\`u
$(Y,C_{1})$ est klt et $C_{1}$ ne contient pas $F_1$ dans son 
support\footnote{
On conseille au lecteur l'exercice suivant : soit 
$\pi : \tilde S \to S$ l'\'eclatement de deux points infiniment
voisins dans une surface lisse $S$. On note $E_1$ et $E_2$ les deux courbes 
exceptionnelles, $p$ le point de $S$ tel que
$\pi^{-1}(p) = E_1 \cup E_2$.  
Soient $C_m$ une courbe poss\'edant une
singularit\'e nodale d'ordre $m$ en $p$, $\Delta =0$ sur $S$. 
Alors, $\pi$ est une log-r\'esolution de $\Delta$ 
telle que $\pi^{-1}(p) \cup C_m'$ est \`a croisements normaux simples
et on a ${\rm Nklt}(\tilde S,\Theta_2)=E_1 \cup E_2$
et ${\rm Nklt}(\tilde S,\Theta_1)=E_1$ o\`u 
$\Theta_1=(2/m)C_m '+E_1$, $\Theta_2=(3/m)C_m '+2E_1+E_2$
et $K_{\tilde S}+\Theta_i \sim_{\pi} 0$.}.
\end{enumerate} 
\end{prop}

{\bf D\'emonstration.}

{\em Etape 1. Construction d'un diviseur $\Theta_0$ auxiliaire.}

On \'ecrit la d\'ecomposition~:
$$ K_Y + \Theta = g^*(K_X+\Delta) + E \sim
E $$ 
o\`u $\Theta $ et $E$ sont effectifs, 
sans composante commune,
$E$ est $g$-exceptionnel et $g_* \Theta =\Delta$.
Comme les composantes irr\'eductibles $F_j$ de $F=\pi^{-1}(s)$ peuvent \^etre
$g$-exceptionnelles, on \'ecrit $E=\sum_{j}e_j F_j + \bar E$
o\`u les $e_j$ sont $\geq 0$ et le support de $\bar E$
ne contient pas de $F_j$. 

Evidemment, $(Y,\Theta )$ est klt,
$\Theta \geq \Delta' =g^*A+B'$
(comme toujours, un $'$ d\'esigne la transform\'ee stricte)
car on a suppos\'e que $A$ est suffisamment g\'en\'eral.
Finalement, $\Theta = g^*A+C$
o\`u $C$ est effectif et $(Y,C)$ est klt. 
On note au passage que le support de 
$C+\textup{Exc}(g)$ est \`a croisements normaux simples.
Soit $N$ un $\mathbb Q$-diviseur effectif et $g$-exceptionnel tel que 
$g^*A-N$ soit ample. On remarque que pour tout $0<\varepsilon \le 1$, 
$g^*A-\varepsilon N=(1-\varepsilon)g^*A+\varepsilon (g^*A-N)$ est encore ample.
On v\'erifie facilement que pour $0<\varepsilon_0 \le 1$ convenable, la paire 
$(Y,C+\varepsilon_0 N)$ est klt. 
On choisit $A_0\sim g^*A-\varepsilon_0 N$ suffisamment g\'en\'eral et 
on pose $B_0=C+\varepsilon N$ et $\Theta_0 = A_0 + B_0$~; 
$(Y,\Theta_0)$ et $(Y,B_0)$ sont des paires klt et les supports 
de $\Theta_0$ et $B_0$ sont \`a 
croisements normaux simples. On peut \'egalement supposer ici
que $B_0$ 
et $E$ n'ont pas de composante commune, quitte \`a remplacer
$B_0$, $\Theta_0$ et $E$ par les diviseurs 
obtenus en simplifiant les composantes communes.



\medskip

{\em Etape 2. Num\'erotation des composantes.}

Soit $l$ un entier, soient
$G_1, G_2,\dots,G_l$ des diviseurs amples g\'en\'eraux sur $S$ passant
par $s$ et soit $G=G_1+\cdots+G_l$. 
Ecrivons $\pi^*G= \sum g_j F_j + \bar G$ o\`u
$\bar G$ ne contient pas de $F_j$ dans son support. 
Si $t_0$ est un rationnel $>0$
et suffisamment petit {\em ind\'ependamment de $l$}, 
alors la paire $(Y,t_0 \bar G)$ est klt. Choisissons 
maintenant $l$ de sorte que $t_0 g_j -e_j >1$ pour tout $j$ compris
entre $1$ et $k$.

Ecrivons pour $0 < t \leq t_0$~:
$$ K_Y + \Theta_0 + t \pi^*G - \sum_{j} e_j F_j = 
K_Y + A_0 + B_0 + \sum_j (t g_j -e_j) F_j +
t \bar G.$$

Choisissons des rationnels $\varepsilon_1,\ldots,\varepsilon_k$
strictements positifs de sorte que $A_0+\sum_j \varepsilon_j F_j$
soit ample et de sorte que les rationnels $\dfrac{1+\varepsilon_j +e_j}{g_j}$
soient deux \`a deux distincts et $<t_0$. On renum\'erote alors
les composantes afin que 
$$ \dfrac{1+\varepsilon_1 +e_1}{g_1}< \dfrac{1+\varepsilon_2 +e_2}{g_2}<\cdots
< \dfrac{1+\varepsilon_k +e_k}{g_k}$$
et on pose $t_i= \dfrac{1+\varepsilon_i +e_i}{g_i}$.

\medskip

{\em Etape 3. Construction des $\Theta_i$.}

L'identit\'e ci-dessus se r\'e-\'ecrit alors 
$$ K_Y + \Theta_0 + t_i \pi^*G - \sum_{j} e_j F_j = 
K_Y + (A_0+ \sum_{j=1}^k \varepsilon_j F_j ) + B_0 + \sum_{j=1}^k (t_i g_j -e_j
-\varepsilon_j) F_j +
t_i \bar G.$$  

Soit alors, comme dans la preuve du lemme \ref{perturb},
un $\mathbb Q$-diviseur $\bar A_0$ g\'en\'eral et lin\'eairement \'equivalent 
\`a $A_0+\sum_{j=1}^k \varepsilon_j F_j$.
On pose alors 
$$\Theta _{i}= \bar A_0 + B_0 + \sum_{j=1}^k \max (t_i g_j- 
e_j - \varepsilon_j,0) F_j + t_i \bar G  \; \mbox{ et }\,
B_i= B_0 + \sum_{j=1}^k \max(t_i g_j- 
e_j - \varepsilon_j,0)F_j + t_i \bar G.$$

Par construction, on a
$$  K_Y + \Theta _{i} \sim K_Y + \Theta  - \sum_{j=1}^ke_j F_j 
+ \sum_{j=1}^k \big ( \max (t_i g_j- 
e_j - \varepsilon_j,0) - (t_ig_j-e_j- \varepsilon_j) \big) F_j $$
d'o\`u 
$$  K_Y + \Theta _{i} \sim \bar E + \sum_{j=1}^k \big ( \max (t_i g_j- 
e_j - \varepsilon_j,0) - (t_ig_j-e_j- \varepsilon_j ) \big) F_j 
=: \bar E_i.$$
La paire
$(Y,\Theta _{i})$ convient : en effet, si $t_i g_j- 
e_j - \varepsilon_j < 0$ alors $e_j >0$ donc
$F_j$ est $g$-exceptionnelle. 
\finpreuve

\subsection{Deuxi\`eme \'etape de la preuve 
du th\'eor\`eme \ref{HMmain} : 
extension de sections}

Si $\Gamma$ est un $\mathbb Q$-diviseur effectif, on note 
$\{ \Gamma \}$ sa partie fractionnaire.  



\begin{prop}\label{kodfibre} 
Avec les notations pr\'ec\'edentes,
pour tout $1\leq i \leq k$
$$ \kappa(F_i,K_{F_i}+\{\Theta_{i}\}_{|F_i}) \le 0.$$
\end{prop}

La preuve de cette proposition repose sur un th\'eor\`eme 
d'extension de sections d\^u \`a Hacon et ${\rm M^c}$Kernan
(\cite{HM06}, Corollary 3.17).
Il est l'aboutis\-se\-ment d'une
s\'erie de travaux initi\'es par Siu, Tsuji et Takayama. 

\begin{theo} {\bf (Hacon et ${\rm \bf M^c}$Kernan)} 
Soient $\pi~: Y \to S$ un morphisme projectif avec $S$ affine,
$Y$ lisse de dimension $d$. Soit $H$ un diviseur ample sur $Y$.
Soient $V$ une hypersurface lisse de $Y$
et $\Gamma$ un $\mathbb Q$-diviseur effectif 
dont le support est \`a croisements normaux simples
tel que $\Gamma = \{\Gamma\} +V$. Soit enfin $C$ 
un $\mathbb Q$-diviseur effectif
dont le support ne contient pas $V$.
On suppose 
\begin{enumerate}
\item que $K_V+(\Gamma-V)_{|V}$ est pseudo-effectif
\item et qu'il existe un $\mathbb Q$-diviseur effectif $G\sim K_Y+\Gamma+C$ ne
contenant aucune intersection (non vide) de composantes de $\Gamma$.
\end{enumerate}
Alors, il existe un entier $k \ge 1$ tel que pour tout $m\ge 1$, 
l'image
de l'application naturelle
$$H^0(Y,mk(K_Y+\Gamma+C)+k(d+1)H) 
\to H^0(V,mk(K_Y+\Gamma+C)_{|V}+k(d+1)H_{|V}) $$
contient
$H^0(V,mk(K_Y+\Gamma)_{|V}+kH_{|V})$
o\`u l'inclusion 
$$H^0(V,mk(K_Y+\Gamma)_{|V}+kH_{|V})\subset 
H^0(V,mk(K_Y+\Gamma+C)_{|V}+k(d+1)H_{|V})$$
est induite par le diviseur $kmC_{|V}+kdH_{|V}$.
\end{theo}


Montrons comment on utilise ce r\'esultat pour d\'emontrer
la proposition \ref{kodfibre}.

\medskip

{\bf D\'emonstration de la proposition \ref{kodfibre}.}
On le fait pour la ``premi\`ere composante'', \`a savoir $F_1$.
On applique le r\'esultat ci-dessus
\`a $V=F_1$,
$\Gamma = \Theta_{1} = \bar A_0 +C_{1}+F_1$ et $C=0$.
On peut toujours supposer qu'un multiple entier non nul 
de $K_V+(\Gamma-V)_{|V}$ est effectif (et donc
certainement pseudo-effectif) sinon il n'y a rien \`a d\'emontrer.
La condition (1) est donc satisfaite.

De m\^eme,
$K_Y+\Theta_{1} \sim \bar E_1 $, donc le lieu 
base de $K_Y+\Theta_{1}$ est contenu dans $\bar E_1$ (rappelons que
$S$ est affine donc ce qui vient de $S$ peut \^etre suppos\'e sans
point base). Or $\pi^{-1}(s) + {\rm Exc}(g) + \Delta'$  
est un diviseur \`a croisements normaux simples
et $\Theta_{1}$ et $\bar E_1$ n'ont pas de composantes communes,
une intersection de composantes de $\Theta_{1}$
ne peut donc pas \^etre incluse dans $\bar E_1$~: la condition
(2) est ainsi aussi satisfaite.

Le th\'eor\`eme pr\'ec\'edent affirme alors
qu'il existe un entier $k \ge 1$ tel que pour tout $m\ge 1$, l'image
de l'application naturelle
$$H^0(Y,mk(K_Y+\Theta_{1})+k(d+1)H) \to 
H^0(F_1,mk(K_{F_1}+\{\Theta_{1}\}_{|F_1})+k(d+1)H_{|F_1})
$$
contient
$H^0(F_1,mk(K_{F_1}+\{\Theta_{1}\}_{|F_1})+kH_{|F_1})$
o\`u l'inclusion 
$$H^0(F_1,mk(K_{F_1}+\{\Theta_{1}\}_{|F_1})+kH_{|F_1})\subset
H^0(F_1,mk(K_{F_1}+\{\Theta_{1}\}_{|F_1})+k(d+1)H_{|F_1})$$
est induite par le diviseur $kdH_{|F_1}$ et 
o\`u $H$ est un diviseur ample sur $Y$ et $d$ est la dimension de $Y$.

On en d\'eduit que 
\begin{align*} 
\dim H^0(F_1,m(K_{F_1}+
\{\Theta_{1}\}_{|F_1}))  
& \leq \dim H^0(F_1,mk(K_{F_1}+\{\Theta_{1}\}_{|F_1})+kH_{|F_1})\\
& \le {\rm rg} (\pi_* \mathcal O_Y (mk(K_Y+\Theta_{1})+k(d+1)H))_s\\
& = {\rm rg} (\pi_* \mathcal O_Y (mk(\bar E_1)+k(d+1)H))_s.\\
\end{align*}

Soit $D$ un diviseur sur $X$ tel que $g^*D-k(d+1)H$ soit effectif.
On a finalement

\begin{align*} 
{\rm rg} (\pi_* \mathcal O_Y (mk(\bar E_1)+k(d+1)H))_s
& \le {\rm rg} (\pi_* \mathcal O_Y (mk(\bar E_1)+g^*D))_s\\
& = {\rm rg}(f_*(g_* \mathcal O_Y (mk \bar E_1+g^*D)))_s\\
& = {\rm rg} (f_* \mathcal O_X (D))_s\\
\end{align*}
puisque 
$g_* \mathcal O_Y (m \bar E_1)= \mathcal O_X$.
On a bien
$ \kappa(F_1,K_{F_1}+\{\Theta_{1}\}_{|F_1}) \le 0.$

Pour la composante $F_i$, il s'agit du m\^eme
raisonnement avec $V:=F_i$, $\Gamma := \{\Theta_{i}\}+F_i$ et 
$C=\Theta_i-\Gamma$.\finpreuve

\subsection{Derni\`ere \'etape de la preuve du th\'eor\`eme \ref{HMmain}}

Les notations et la situation sont celles de l'\'etape pr\'ec\'edente.
La proposition suivante termine la preuve du 
th\'eor\`eme \ref{HMmain}.

\begin{prop}\label{Frc}
\begin{enumerate}
\item La composante $F_1$ est rationnellement connexe,
\item pour tout $i\geq 2$, $F_i$ est rationnellement connexe
par cha\^{\i}nes modulo $$F_i\cap {\rm Nklt}(Y,\Theta_{i-1})
= F_i\cap(F_1\cup \ldots \cup F_{i-1}),$$ ce qui signifie 
que pour tout $x \in F_i$, il existe une cha\^{\i}ne de courbes
rationnelles joignant $x$ \`a un point de 
$F_i\cap (F_1\cup \ldots \cup F_{i-1})$.
\end{enumerate}
\end{prop}

{\bf D\'emonstration.}
On commence par d\'emontrer le premier point. 
Comme $F_1$ est lisse,
il suffit de montrer que $F_1$ est rationnellement connexe
par cha\^{\i}nes (corollaire \ref{coversus}), 
autrement dit que son quotient rationnel $R:=R(F_1)$ est
r\'eduit \`a un point. 

Soit donc $t: F_1 \dashrightarrow R$ le quotient
rationnel de $F_1$. On peut supposer que $R$ est lisse.
La paire klt $(F_1,\{\Theta_{1}\}_{|F_1})$ 
va jouer un r\^ole cl\'e. Rappelons la liste de ses propri\'et\'es :
\begin{enumerate}
\item la dimension de Kodaira de $K_{F_1}+\{\Theta_{1}\}_{|F_1}$ vaut $0$
ou $-\infty$ : ceci d\'ecoule imm\'edia\-tement de la proposition
\ref{kodfibre},
\item $\{\Theta _{1}\}_{|F_1} = (\bar A_0) _{|F_1} +(C_{1})_{|F_1}$ 
o\`u $\bar A_0$ est ample
sur $Y$, la paire $(F_{1},(C_{1})_{|F_1})$ est klt
et 
$$K_{F_1}+\{\Theta_{1}\}_{|F_1} \sim (\bar E_1)_{|F_1}$$
est effectif. 
\item soient $\varphi~: \tilde F_1 \to F_1$ une suite d'\'eclatements
qui l\`eve les ind\'eterminations de $t$ et $$\tilde t = t\circ  
\varphi : \tilde F_1 \to R$$ le morphisme compos\'e.
Alors la dimension de Kodaira de la restriction de 
$$\varphi^*(K_{F_1}+\{\Theta_{1}\}_{|F_1})$$ 
\`a la fibre g\'en\'erale de $\tilde t$
est $\geq 0$ (car sup\'erieure ou \'egale \`a la 
dimension de Kodaira de la restriction de 
$\varphi^* (\bar E_1 ) _{|F_1}$ 
\`a la fibre g\'en\'erale de $\tilde t$ !).
\end{enumerate} 

On utilise alors \`a nouveau le th\'eor\`eme de positivit\'e d'images directes
de Campana.
Soient donc, comme ci-dessus, 
$\varphi~: \tilde F_1 \to F_1$ une suite d'\'eclatements
qui l\`eve les ind\'e\-termi\-na\-tions de $t$ et $\tilde t = t\circ  
\varphi : \tilde F_1 \to R$ le morphisme compos\'e.
On \'ecrit comme d'habitude la d\'ecomposition 
$$ K_{\tilde F_1}+ \Omega = \varphi^*(K_{F_1}+\{\Theta_{1}\}_{|F_1})+ D$$
o\`u $\Omega$ et $D$ sont effectifs sans composantes communes,
$\varphi_*(\Omega)=\{\Theta_{1}\}_{|F_1}$ et $D$ est $\varphi$-exceptionnel. 
Les propri\'et\'es pr\'ec\'edemment list\'ees pour 
$(F_1,\{\Theta_{1}\}_{|F_1})$ 
se propagent \`a $(\tilde F_1,\Omega)$~:

\begin{enumerate}
\item la dimension de Kodaira de $K_{\tilde F_1}+\Omega$ vaut $0$
ou $-\infty$,
\item la dimension de Kodaira de la restriction de 
$K_{\tilde F_1}+\Omega$
\`a la fibre g\'en\'erale de $\tilde t$
est $\geq 0$,
\item $\Omega = \varphi ^* (\bar A_0)_{|F_1} + \bar B$ o\`u 
$(\tilde F_{1},\bar B)$ est klt.
\end{enumerate}

\medskip

D'apr\`es le lemme \ref{kodaira}, 
il existe un diviseur $H$ ample sur $R$
tel que $\varphi ^* (\bar A_0)_{|F_1} \sim 2\tilde t ^* H +\Gamma$
o\`u $\Gamma$ est un diviseur effectif sur 
$\tilde F_1$ et la paire $(\tilde F_1,\bar B+\Gamma)$
est klt.

La paire $\bar B+\Gamma$ est klt
et 
$$K_{\tilde F_1/R}+\bar B+\Gamma \sim K_{\tilde F_1}-\tilde t^*K_{R}
+ \Omega - \varphi ^* (\bar A_0)_{|F_1}+ \varphi ^* (\bar A_0)_{|F_1}
-2\tilde t ^* H$$
donc
$$(K_{\tilde F_1/R}+\bar B+\Gamma)_{|\tilde t ^{-1}(u)}
\sim (K_{\tilde F_1}
+ \Omega)_{|\tilde t ^{-1}(u)}$$
est de dimension de Kodaira $\geq 0$.

Comme $\varphi ^* (\bar A_0)_{|F_1} \sim 2\tilde t^* H + \Gamma$
et $\Omega = \varphi ^* (\bar A_0)_{|F_1} + \bar B$,
on en d\'eduit alors
que 
\begin{align*}
0 = \kappa(\tilde F_1,K_{\tilde F_1}+\Omega)&  
= \kappa(\tilde F_1,K_{\tilde F_1}+ \bar B + \Gamma + 2\tilde t ^* H)\\
& \geq \kappa(R,K_{R}+H).
\end{align*}

Or $R$ n'est pas 
unir\'egl\'e (th\'eor\`eme \ref{quotientuni}),  
le th\'eor\`eme
\ref{carunireg} affirme donc que $K_{R}$ est pseudo-effectif,
et par suite que $$\kappa(R,K_{R}+H)=\dim R.$$ 
Contradiction sauf si $R$ est un point ! 
C'est magnifique et \c{c}a termine la preuve du premier point.

Pour montrer que $F_i$ est rationnellement connexe
par cha\^{\i}nes modulo $$F_i\cap {\rm Nklt}(Y,\Theta_{i-1})
= F_i\cap(F_1\cup \ldots \cup F_{i-1}),$$
il suffit d'observer que si $t~: F_i \dashrightarrow R(F_i)$
est le quotient rationnel de $F_i$, alors 
soit $$ F_i\cap {\rm Nklt}(Y,\Theta_{i-1})
= F_i\cap(F_1\cup \ldots \cup F_{i-1})$$ 
domine $R(F_i)$, 
soit le m\^eme raisonnement que pr\'ec\'edemment montre
que $R(F_i)$ est r\'eduit \`a un point.\finpreuve

\end{document}